\pdfoutput=1
\RequirePackage{silence}
\documentclass[12pt,a4paper,svgnames]{amsart}
\usepackage[hmarginratio={1:1},vmarginratio={1:1},lmargin=60.0pt,tmargin=58.0pt]{geometry}

\usepackage{booktabs}
\allowdisplaybreaks

\usepackage[T1]{fontenc}
\usepackage{latexsym,exscale,amsfonts,amssymb,mathtools}
\usepackage{amsmath,amsthm,amsfonts,amssymb,amscd,textcomp,bbm}
\usepackage{mathrsfs,stackrel}
\usepackage{etoolbox}
\usepackage{tikz,tikz-cd}
\usetikzlibrary{matrix,arrows.meta,decorations.markings,shapes.multipart,decorations.pathreplacing,decorations.pathmorphing}
\usetikzlibrary{shadings}

\pgfdeclareverticalshading{rainbow}{100bp}
{color(0bp)=(red); color(25bp)=(red); color(35bp)=(yellow);
color(45bp)=(green); color(55bp)=(cyan); color(65bp)=(blue);
color(75bp)=(violet); color(100bp)=(violet)}

\usepackage{enumitem}
\setlist[enumerate]{itemsep=0.15cm,label=\emph{\upshape(\alph*)}}
\setlist[enumerate,2]{itemsep=0.15cm,label=\emph{\upshape(\roman*)}}
\usepackage{aliascnt}
\usepackage{ytableau}
\usepackage{xparse}
\usepackage{cite}

\usepackage{dynkin-diagrams}

\usepackage{array}
\newcolumntype{C}{>{$}c<{$}}

\usepackage{xcolor,colortbl}

\definecolor{mygray}{gray}{0.6}
\definecolor{mygraydark}{gray}{0.4}
\definecolor{mygraylight}{gray}{0.85}
\definecolor{spinach}{RGB}{46,139,87}
\definecolor{tomato}{RGB}{255,99,71}
\definecolor{orchid}{RGB}{143,40,194}
\definecolor{neon}{RGB}{77,77,255}
\definecolor{pumpkin}{RGB}{224,180,80}
\definecolor{citron}{RGB}{190,180,90}

\definecolor{lava}{RGB}{207,16,32}
\definecolor{cream}{RGB}{255,253,208}
\definecolor{verdigris}{RGB}{67,179,174}
\definecolor{Black}{RGB}{0,0,0}
\definecolor{mydarkblue}{RGB}{10,10,170}
\definecolor{darkspinach}{RGB}{20,70,20}
\definecolor{darktomato}{RGB}{155,40,30}
\definecolor{darkorchid}{RGB}{50,10,100}
\definecolor{darklava}{RGB}{150,8,16}

\usepackage{enumitem}
\setlist[enumerate]{itemsep=0.15cm,label=\emph{\upshape(\alph*)}}
\setlist[enumerate,2]{itemsep=0.15cm,label=\emph{\upshape(\roman*)}}
\setlist[enumerate,3]{itemsep=0.15cm,label=\emph{\upshape(\Alph*)}}

\let\emph\relax
\DeclareTextFontCommand{\emph}{\bfseries\em}

\newcommand{\acts}{\centerdot}

\renewcommand{\dots}{\text{...}}
\renewcommand{\vdots}{\rotatebox{90}{\text{...}}}

\newcommand{\placeholder}{{}_{-}}
\newcommand{\mystrut}{\rule[-0.2\baselineskip]{0pt}{0.9\baselineskip}}

\usepackage{todonotes}

\newcommand{\fieldsymbol}{\scalebox{0.4}{$\heartsuit$}}

\newcommand{\C}{\mathbb{C}}

\newcommand{\N}{\mathbb{Z}_{\geq 0}}
\newcommand{\Z}{\mathbb{Z}}
\newcommand{\Aa}{\mathbb{A}}
\newcommand{\Q}{\mathbb{Q}}
\newcommand{\Ff}{\mathbb{F}}
\newcommand{\A}{\mathbb{A}_{\fieldsymbol}}
\newcommand{\K}{\mathbb{Q}_{\fieldsymbol}}
\newcommand{\F}{\mathbb{F}_{\fieldsymbol}}

\newcommand{\vcirc}{\circ}
\newcommand{\hcirc}{\otimes}
\newcommand{\wtlnit}{\mathbbm{1}}

\newcommand{\setstuff}[1]{\mathrm{#1}}
\newcommand{\catstuff}[1]{\mathbf{#1}}
\newcommand{\functorstuff}[1]{\mathcal{#1}}

\newcommand{\obstuff}[1]{\mathtt{#1}}
\newcommand{\morstuff}[1]{\mathrm{#1}}

\newcommand{\idmor}{\morstuff{1}}

\newcommand{\End}{\setstuff{End}}
\newcommand{\Hom}{\setstuff{Hom}}

\newcommand{\Kar}[1][1]{{#1}^{\oplus,\subset_{\oplus}}}

\newcommand{\sbinom}[2]{{\textstyle\binom{#1}{#2}}}

\newcommand{\tr}{\mathrm{tr}}
\newcommand{\ev}{\mathrm{ev}}
\newcommand{\coev}{\mathrm{coev}}

\newcommand{\bigboxtimes}{\scalebox{1.5}{$\boxtimes$}}

\newcommand{\semisimple}[1]{\scalebox{0.9}{$\overline{#1}$}}

\newcommand{\pivo}{\ast}

\newcommand{\rank}{\mathrm{rank}}

\newcommand{\lsym}{N}
\newcommand{\llsym}{n}
\newcommand{\rsym}{m}

\newcommand{\slgroupf}[1][\C]{\mathrm{SL}_{\lsym}(#1)}
\newcommand{\glgroupf}[1][\C]{\mathrm{GL}_{\lsym}(#1)}
\newcommand{\glgroupff}[1][\C]{\mathrm{GL}_{\rsym}(#1)}
\newcommand{\ogroupf}[1][\C]{\mathrm{O}_{\lsym}(#1)}
\newcommand{\spgroupf}[1][\C]{\mathrm{SP}_{2\lsym}(#1)}
\newcommand{\sogroupf}[1][\C]{\mathrm{SO}_{3}(#1)}

\newcommand{\ogroup}[1][\lsym]{\mathrm{O}(#1)}
\newcommand{\sogroup}[1][\lsym]{\mathrm{SO}(#1)}
\newcommand{\glgroup}[1][\lsym]{\mathrm{GL}(#1)}
\newcommand{\sgroup}[1][n]{\mathrm{S}(#1)}

\newcommand{\weyl}[1][\lambda]{\Delta_{\A}(#1)}
\newcommand{\coweyl}[1][\lambda]{\nabla_{\A}(#1)}
\newcommand{\aweyl}[1][\lambda]{\Delta_{\Aa}(#1)}
\newcommand{\acoweyl}[1][\lambda]{\nabla_{\Aa}(#1)}
\newcommand{\fweyl}[1][\lambda]{\Delta_{\F}(#1)}
\newcommand{\fcoweyl}[1][\lambda]{\nabla_{\F}(#1)}
\newcommand{\ffweyl}[1][\lambda]{\Delta_{\Ff}(#1)}
\newcommand{\ffcoweyl}[1][\lambda]{\nabla_{\Ff}(#1)}

\newcommand{\qweyl}[1][\lambda]{\Delta_{\Q}(#1)}

\newcommand{\fsimple}[1][\lambda]{\setstuff{L}_{\F}(#1)}
\newcommand{\ffsimple}[1][\lambda]{\setstuff{L}_{\Ff}(#1)}
\newcommand{\ksimple}[1][\lambda]{\setstuff{L}_{\K}(#1)}
\newcommand{\qsimple}[1][\lambda]{\setstuff{L}_{\Q}(#1)}

\newcommand{\ftilting}[1][\lambda]{\setstuff{T}_{\F}(#1)}
\newcommand{\fftilting}[1][\lambda]{\setstuff{T}_{\Ff}(#1)}

\newcommand{\ext}[1][\ast]{\Lambda^{#1}}

\newcommand{\NBoxm}{\square_{\lsym{\times}\rsym}}

\newcommand{\ohor}{<_{h}}
\newcommand{\overt}{<_{v}}

\newcommand{\hreading}{\phi_{h}}
\newcommand{\vreading}{\phi_{v}}

\newcommand{\web}[1][\A]{\catstuff{Web}_{#1}\big(\mathrm{O}(\lsym)\big)}
\newcommand{\webb}[1][\A]{\catstuff{Web}_{#1}\big(\glgroup\big)}
\newcommand{\webF}[1][\F]{\catstuff{Web}_{#1}\big(\mathrm{O}(\lsym)\big)}

\newcommand{\sweb}[1][\A]{\catstuff{Web}_{#1}^{S}\big(\mathrm{O}(\lsym)\big)}

\newcommand{\Web}[1][\A]{\catstuff{Web}_{#1}(\mathrm{O}(\lsym))}

\newcommand{\brauer}[1][\A]{\catstuff{Br}_{#1}(N)}
\newcommand{\brauertwo}[1][\A]{\catstuff{Br}}

\newcommand{\ifunctor}[1][\A]{\functorstuff{E}_{#1}^{\mathrm{GL}\to\mathrm{O}}}

\newcommand{\pfunctor}[1][\F]{\functorstuff{P}_{#1}^{\ogroup}}
\newcommand{\ppfunctor}[1][\F]{\functorstuff{P}_{#1}^{\mathrm{GL}(\lsym)}}
\newcommand{\sfunctor}[1][\F]{\functorstuff{S}_{#1}^{\gso}}

\newcommand{\glm}[1][\rsym]{\mathfrak{gl}_{#1}}

\newcommand{\ahowe}[1][\A]{\functorstuff{H}_{#1}^{gl}(\lsym,\rsym)}
\newcommand{\dhowe}[1][\A]{\functorstuff{H}_{#1}^{so}(\lsym,\rsym)}
\newcommand{\ddhowe}[1][\A]{\functorstuff{HS}_{#1}^{so}(\lsym,\rsym)}

\newcommand{\hideal}{\catstuff{I}_{>\lsym}}

\newcommand{\sand}[2]{\setstuff{X}_{#1}^{#2}}
\newcommand{\sandtwo}[2]{\setstuff{O}_{#1}^{#2}}

\newcommand{\Fund}{\catstuff{Fund}_{\A}\big(\mathrm{O}(\lsym)\big)}
\newcommand{\FundF}{\catstuff{Fund}_{\F}\big(\mathrm{O}(\lsym)\big)}
\newcommand{\FundFF}{\catstuff{Fund}_{\F}\big(\mathrm{GL}(\lsym)\big)}
\newcommand{\Rep}{\catstuff{Rep}}
\newcommand{\Tilt}{\catstuff{Tilt}}
\newcommand{\Ver}{\catstuff{Ver}}
\newcommand{\Vero}[1][\F]{\catstuff{Ver}_{#1}\big(\mathrm{O}(\lsym)\big)}

\newcommand{\Wcat}{\catstuff{Fil}^{\Delta}}
\newcommand{\DWcat}{\catstuff{Fil}^{\nabla}}

\newcommand{\id}{\mathrm{id}}
\newcommand{\Ya}{\mathbf{Y}}
\newcommand{\Pa}{\reflectbox{\rotatebox[origin=c]{180}{\textbf{Y}}}}
\newcommand{\cu}{\mathbf{U}}
\newcommand{\ca}{\reflectbox{\rotatebox[origin=c]{180}{\textbf{U}}}}
\newcommand{\one}{\mathbf{1}}

\newcommand{\TiltF}{\catstuff{Tilt}_{\F}\big(\mathrm{O}(\lsym)\big)}
\newcommand{\TiltFi}{\catstuff{Tilt}_{\F}\big(\mathrm{O}(\lsym_{i})\big)}
\newcommand{\wt}{\mathrm{wt}\,}
\newcommand{\Comp}{\Pi_{\rsym}}
\newcommand{\CompN}{\Pi_{\rsym}^{\leq\lsym}}
\newcommand{\DCompN}{\Pi_{\rsym,+}^{\leq\lsym}}
\newcommand{\ACompN}{\Pi_{\rsym,-}^{\leq\lsym}}
\newcommand{\Par}{\Lambda_{+}}
\newcommand{\ParN}{\Lambda_{+}^{\ogroup}}
\newcommand{\ParNm}{\Lambda_{+,\leq\rsym}^{\ogroup}}
\newcommand{\led}{\trianglelefteq}
\newcommand{\ltd}{\triangleleft}
\newcommand{\ltO}{<_{\ogroup}}
\newcommand{\leO}{\leq_{\ogroup}}
\newcommand{\legso}{\leq_{\gso}}

\newcommand{\ct}{\mathrm{ct}}
\newcommand{\IndN}{\functorstuff{I}^{\mathrm{O}}_{\mathrm{SO}}}
\newcommand{\ResN}{\functorstuff{R}_{\mathrm{O}}^{\mathrm{SO}}}
\newcommand{\hs}{\mathrm{hs}}
\newcommand{\Ext}{\mathrm{Ext}}
\newcommand{\gso}{\mathfrak{so}_{2\rsym}}
\newcommand{\lso}{\mathfrak{so}_{\lsym}}

\newcommand{\XX}{\mathbf{X}}
\newcommand{\wtk}{\obstuff{K}}
\newcommand{\wtl}{\obstuff{L}}

\newcommand{\wta}{\obstuff{a}}
\newcommand{\wtb}{\obstuff{b}}

\newcommand{\xp}{\mathcal{Y}}

\newcommand{\cbasis}{\setstuff{B}}

\newcommand{\uplain}[1][\A]{\setstuff{U}_{#1}}
\newcommand{\uplainsecond}[1][\A]{\setstuff{V}_{#1}}
\newcommand{\udot}[1][\A]{\dot{\setstuff{U}}_{#1}}
\newcommand{\udotsecond}[1][\A]{\dot{\setstuff{V}}_{#1}}

\newcommand{\splain}[1][\A]{\setstuff{S}_{#1}^{\leq\lsym}}
\newcommand{\spi}[1][\A]{\setstuff{S}_{#1}^{\leq\lsym}}

\usepackage[all]{xy}
\usepackage{tikz}
\usetikzlibrary{cd}
\usetikzlibrary{decorations}
\usetikzlibrary{decorations.markings}
\usetikzlibrary{decorations.pathreplacing}
\usetikzlibrary{decorations.pathmorphing}
\usetikzlibrary{arrows.meta,shapes,positioning,matrix,calc}
\usetikzlibrary{shapes.callouts}
\tikzset{anchorbase/.style={baseline={([yshift=-0.5ex]current bounding box.center)}},
tinynodes/.style={font=\tiny,text height=0.25ex,text depth=0.05ex},
smallnodes/.style={font=\scriptsize,text height=0.75ex,text depth=0.15ex},
usual/.style={line width=2.0,color=black},
usual0/.style={line width=2.0,color=magenta},
usual1/.style={line width=2.0,color=orchid},
usual2/.style={line width=2.0,color=blue},
mor/.style={line width=0.75,color=black,fill=cream},
}

\numberwithin{equation}{subsection}

\newtheorem{theorem}{Theorem}[section]

\makeatletter
\newcommand{\newthmwithalias}[2]{%
\newaliascnt{#1}{theorem}%
\newtheorem{#1}[#1]{#2}%
\aliascntresetthe{#1}%
}
\makeatother

\newthmwithalias{lemma}{Lemma}
\newthmwithalias{definition}{Definition}
\newthmwithalias{notation}{Notation}
\newthmwithalias{corollary}{Corollary}
\newthmwithalias{proposition}{Proposition}
\newthmwithalias{conjecture}{Conjecture}
\newthmwithalias{example}{Example}
\newthmwithalias{remark}{Remark}
\newthmwithalias{question}{Question}
\newthmwithalias{assumption}{Assumption}

\AtEndEnvironment{definition}{\null\hfill$\Diamond$}
\AtEndEnvironment{notation}{\null\hfill$\Diamond$}
\AtEndEnvironment{example}{\null\hfill$\Diamond$}
\AtEndEnvironment{remark}{\null\hfill$\Diamond$}
\AtEndEnvironment{corollary}{\null\hfill$\square$}

\usepackage[hypertexnames=false]{hyperref}
\usepackage{bookmark}
\hypersetup{
pdftoolbar=true,
pdfmenubar=true,
pdffitwindow=false,
pdfstartview={FitH},
pdftitle={Orthogonal webs and semisimplification},
pdfauthor={Elijah Bodish and Daniel Tubbenhauer},
pdfsubject={},
pdfcreator={Elijah Bodish and Daniel Tubbenhauer},
pdfproducer={Elijah Bodish and Daniel Tubbenhauer},
pdfkeywords={},
pdfnewwindow=true,
colorlinks=true,
linkcolor=mydarkblue,
citecolor=teal,
filecolor=magenta,
urlcolor=orchid,
linkbordercolor=lava,
citebordercolor=teal,
urlbordercolor=orchid,
linktocpage=true
}

\usepackage[nameinlink, noabbrev]{cleveref}

\crefname{theorem}{Theorem}{Theorems}
\crefname{lemma}{Lemma}{Lemmas}
\crefname{definition}{Definition}{Definitions}
\crefname{notation}{Notation}{Notations}
\crefname{corollary}{Corollary}{Corollaries}
\crefname{proposition}{Proposition}{Propositions}
\crefname{conjecture}{Conjecture}{Conjectures}
\crefname{example}{Example}{Examples}
\crefname{remark}{Remark}{Remarks}
\crefname{question}{Question}{Questions}
\crefname{assumption}{Assumption}{Assumptions}
\crefname{equation}{Equation}{Equations}
\crefname{section}{Section}{Sections}%
\crefname{subsection}{Subsection}{Subsections}%
\crefname{subsubsection}{Subsubsection}{Subsubsections}%

\newcommand{\nnfootnote}[1]{%
\begin{NoHyper}
\renewcommand\thefootnote{}\footnote{#1}%
\addtocounter{footnote}{-1}%
\end{NoHyper}
}

\begin{document}
\title[Orthogonal webs and semisimplification]{Orthogonal webs and semisimplification}
\author[E. Bodish and D. Tubbenhauer]{Elijah Bodish and Daniel Tubbenhauer}

\address{E.B.: MIT, Department of Mathematics, Building 2, Office 2-178, Cambridge, MA 02139, United States, https://orcid.org/0000-0003-1499-5136}
\email{ebodish@mit.edu}

\address{D.T.: The University of Sydney, School of Mathematics and Statistics F07, Office Carslaw 827, NSW 2006, Australia, \href{http://www.dtubbenhauer.com}{www.dtubbenhauer.com}, https://orcid.org/0000-0001-7265-5047}
\email{daniel.tubbenhauer@sydney.edu.au}

\begin{abstract}
We define a diagrammatic category 
that is equivalent 
to tilting representations for the orthogonal group. Our construction 
works in characteristic not equal to two.
We also describe the semisimplification of this category.
\end{abstract}

\nnfootnote{\textit{Mathematics Subject Classification 2020.} Primary: 18M05, 20G05; Secondary: 18M30, 20J15.}
\nnfootnote{\textit{Keywords.} Representations of linear algebraic groups, webs, Howe duality, diagrammatic presentation, positive characteristic, semisimplification.}

\addtocontents{toc}{\protect\setcounter{tocdepth}{1}}

\maketitle

\tableofcontents

%%%%%%%%%%%%%%%%

\section{Introduction}\label{S:Intro}

%%%%%%%%%%%%%%%%

In this paper we define a diagrammatic category of \emph{orthogonal webs}
that is equivalent to tilting representations of the orthogonal group.

%%%%%%%%%%%%%%%%

\subsection{Motivation}\label{SS:IntroductionOne}

%%%%%%%%%%%%%%%%

Diagrammatic methods and diagrammatic presentations in representation theory have been studied for over a century. 

An early example is Schur's celebrated \emph{Schur--Weyl
duality} \cite{Schur} relating representations of $\slgroupf$ (or $\glgroupf$) with the symmetric group. The diagrammatic description of a permutation gives rise to diagrammatic methods for $\slgroupf$-representations. In modern terms, we say that there is a functor from a generators and relations (i.e. diagrammatic) monoidal category to the monoidal category of representations of $\slgroupf$, and Schur--Weyl duality says this functor is full. 

The functor is not faithful, but there is an explicit diagrammatic formulation of the kernel which involves the antisymmetrizer on $N$ strands, see e.g. \cite{Ha-murphy-generalized-tl}. Identification of this kernel yields a diagrammatic presentation of the monoidal category of representations of $\slgroupf$. However, since the kernel is described by a complicated sum of diagrams, this presentation is not very aesthetically pleasing from a diagrammatic perspective. In the special case $\lsym = 2$, the kernel has a particularly simple description, as shown by 
Rumer--Teller--Weyl \cite{RuTeWe-sl2}, which resulted in the Temperley--Lieb calculus, a presentation which does look nice diagrammatically. 

Brauer extended Schur's result to $\ogroupf$ (and $\spgroupf$) relating them via \emph{Brauer duality} to Brauer's diagrammatic algebra \cite{Br-brauer-algebra-original}. Brauer's duality, combined with an identification of the kernel in terms of Brauer's diagrammatic algebra, yields a diagrammatic presentation for the monoidal category of representations of $\ogroupf$. Similarly as for $\slgroupf$, the kernel is well-known, see for example \cite{LeZh-brauer-invariant-theory}, but does not seem to admit a diagrammatic description without complicated sums when $\lsym>2$.

With the birth of quantum topology in the 1980s many more 
diagrammatic presentations were found and these often include the identification of the kernels. For example, Yamada's 
presentation for $\sogroupf$-representations \cite{Ya-invariant-graphs}, following ideas from \cite{TeLi-the-tl-paper}, which employs the diagrammatics of \emph{webs} (also known as birdtracks \cite{Cv-bridtracks}, spiders \cite{Ku-spiders-rank-2} etc.). Yamada's presentation is in the spirit 
of Rumer--Teller--Weyl and all relations are beautifully included in the diagrammatics. These webs were then extended to many other settings, and have 
been at the heart of diagrammatic representation theory every since.

The strategy 
that we run in this paper is inspired by the observation of Cautis--Kamnitzer--Morrison 
\cite{CaKaMo-webs-skew-howe} that the kernel under 
\emph{Howe's duality} for $\glgroupf$-$\glgroupff$ 
\cite{Ho-perspectives-invariant-theory} 
has again a pleasing diagrammatic interpretation in terms of webs. 
They masterfully used Howe's duality to define a presentation for $\slgroupf$-representations which generalizes Rumer--Teller--Weyl's presentation. Since the kernel under Schur--Weyl duality for $\lsym>2$ is diagrammatically rather ugly, such a nice presentation did not seem possible
from Schur--Weyl duality itself.

Another upshot of Cautis--Kamnitzer--Morrison's 
description is that it works over any field when one slightly 
modifies the target category to be \emph{tilting representations}. 
Taking the prime characteristic version 
of Howe's $\glgroupf$-$\glgroupff$ duality from \cite{AdRy-tilting-howe-positive-char}
and running the Cau\-tis--Kamnitzer--Morrison strategy, one gets a diagrammatic 
category equivalent to tilting representations of $\slgroupf[\overline{\mathbb{F}_{p}}]$ for any prime $p$. The technicalities however are much more involved, see \cite{El-ladders-clasps} for details.

\begin{remark}
In characteristic zero all representations are tilting, so 
the setting in the above paragraph generalizes Cautis--Kamnitzer--Morrison's result. The same is true in the orthogonal world and we 
will state all of our results for tilting representations.
\end{remark}

This is the starting point of our paper. We begin by fixing $p=0$ or any prime $p\ne 2$. Let $\F$ be an infinite field of characteristic $p$ containing $\sqrt{-1}$. Our main results are:

\begin{enumerate}[label=(\Alph*)]

\item We give a diagrammatic presentation of tilting $\ogroupf[\F]$-representations using 
\emph{orthogonal webs}. This extends the result of Sartori \cite{SaTu-bcd-webs} to prime characteristic.

\item As an application, we give an orthogonal version of the main result of \cite{BrEnAiEtOs-semisimple-tilting}, i.e. we describe the semisimplification of tilting $\ogroupf[\F]$-representations. 
Here $p$ does not need to be bigger than $\lsym$.

\end{enumerate}

\noindent A key ingredient in our proofs is Howe's orthogonal duality in prime characteristic \cite{AdRy-tilting-howe-positive-char}.

Before coming to the main body of the paper we now explain our results in more detail.

%%%%%%%%%%%%%%%%

\subsection{What we prove}\label{SS:IntroductionWhat}

%%%%%%%%%%%%%%%%

A \emph{closed orthogonal pre-web} is a trivalent graph with edges labeled 
with integers $\{1,\dots,\lsym\}$ 
such that we have $k$, $l$ and $k+l$ around 
every trivalent vertex. A \emph{closed orthogonal web} 
is a planar embedding of a closed orthogonal pre-web such that each 
point of intersection is a crossing in the usual sense.
As usual in diagrammatic algebra, cutting these graphs open and putting 
them into a strip with bottom and top boundary points gives a way to define morphisms, called \emph{orthogonal webs}, in a monoidal category. 
Here are two examples, the left one being closed:
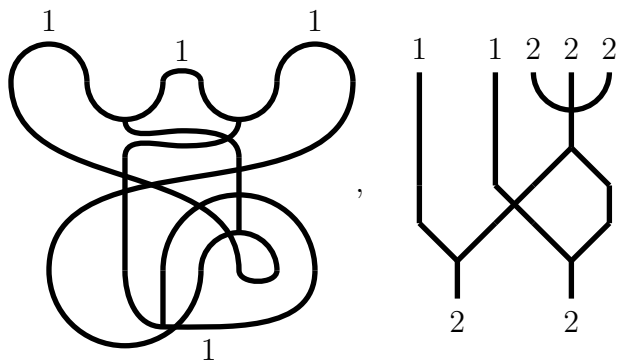
\begin{gather*}
\begin{tikzpicture}[anchorbase,scale=1]
\draw[usual] (0,0) to[out=90,in=180] (0.5,0.5);
\draw[usual] (1,0) to[out=90,in=0] (0.5,0.5);
\draw[usual] (0.5,0.5) to (0.5,1) to (0.5,1.5);
\draw[usual] (-0.5,0) to[out=90,in=180] (0.5,1) to[out=0,in=90] (1.5,0);
\draw[usual] (-1,0) to (-1,1.5);
\draw[usual] (0.5,1.5) to[out=90,in=0] (-0.25,1.85) to[out=180,in=270] (-1,2);
\draw[usual] (-1,1.5) to[out=90,in=180] (-0.25,1.65) to[out=0,in=270] (0.5,2);
\draw[usual] (-0.5,2.5) to[out=270,in=0] (-1,2) to[out=180,in=270] (-1.5,2.5);
\draw[usual] (0,2.5) to[out=270,in=180] (0.5,2) to[out=0,in=270] (1,2.5);
\draw[usual] (-1.5,2.5) to[out=90,in=0] (-2,3)node[above]{$1$} to[out=180,in=90] (-2.5,2.5) to[out=270,in=90] (0.5,0) to [out=270,in=270] (1,0);
\draw[usual] (-0.5,2.5) to[out=90,in=90] (0,2.5)node[above,xshift=-0.25cm,yshift=0.1cm]{$1$};
\draw[usual] (1,2.5) to[out=90,in=180] (1.5,3)node[above]{$1$} to[out=0,in=90] (2,2.5);
\draw[usual] (0,0) to[out=270,in=0] (-1,-1) to[out=180,in=270] (-2,0) to[out=90,in=270] (2,2.5);
\draw[usual] (-1,0) to[out=270,in=180] (-0.5,-0.75)node[below,xshift=0.6cm]{$1$} to[out=0,in=270] (1.5,0);
\draw[usual] (-0.5,0) to (-0.5,-0.75);
\end{tikzpicture}
,\quad
\begin{tikzpicture}[anchorbase,scale=1]
\draw[usual] (-1.5,0) to (-1.5,1.5)node[above]{$1$};
\draw[usual] (-0.5,0) to (-0.5,1.5)node[above]{$1$};
\draw[usual] (0,0) to (0.5,0.5);
\draw[usual] (1,0) to (0.5,0.5);
\draw[usual] (0.5,0.5) to (0.5,1) to (0.5,1.5)node[above]{$2$};
\draw[usual] (0,1.5)node[above]{$2$} to[out=270,in=180] (0.5,1) to[out=0,in=270] (1,1.5) node[above]{$2$};
\draw[usual] (-1.5,-0.5) to (-1.5,0);
\draw[usual] (-0.5,-0.5) to (0,0);
\draw[usual] (0,-0.5) to (-0.5,0);
\draw[usual] (1,-0.5) to (1,0);
\draw[usual] (-1,-1) to (-1.5,-0.5);
\draw[usual] (-1,-1) to (-0.5,-0.5);
\draw[usual] (-1,-1.5)node[below]{$2$} to (-1,-1);
\draw[usual] (0.5,-1) to (0,-0.5);
\draw[usual] (0.5,-1) to (1,-0.5);
\draw[usual] (0.5,-1.5)node[below]{$2$} to (0.5,-1);
\end{tikzpicture}
\end{gather*}

\begin{notation}
If in this or other illustrations an edge is not labeled, then 
its label is determined by other labels and we omitted it to avoid clutter.
\end{notation}

Orthogonal webs form a combinatorial and topological category
akin to the Temperley--Lieb calculus and Yamada's webs.
Let $\F$ be as above. Enriching orthogonal 
webs $\F$-linearly and imposing relations of the form
\begin{gather*}
\scalebox{0.96}{$\begin{tikzpicture}[anchorbase,scale=1]
\draw[usual] (0,0)node[below]{${>}\lsym$} to (0,1.5)node[above]{${>}\lsym$};
\end{tikzpicture}
=0
,\quad
\begin{tikzpicture}[anchorbase,scale=1,rounded corners]
\draw[usual] (0.5,0.35) to (0,0.75)node[left]{$k$} to (0.5,1.15);
\draw[usual] (0.5,0.35) to (1,0.75)node[right]{$l$} to (0.5,1.15);
\draw[usual] (0.5,1.15) to (0.5,1.5)node[above]{$k{+}l$};
\draw[usual] (0.5,0.35) to (0.5,0)node[below]{$k{+}l$};
\end{tikzpicture}
=\binom{k+l}{k}
\cdot
\begin{tikzpicture}[anchorbase,scale=1]
\draw[usual] (0.5,0)node[below]{$k{+}l$} to (0.5,1.5)node[above]{$k{+}l$};
\end{tikzpicture}
,\quad
\begin{tikzpicture}[anchorbase,scale=1]
\draw[usual] (0,0)node[left]{$k$} to[out=270,in=180] (0.5,-0.5) to[out=0,in=270] (1,0);
\draw[usual] (0,0) to[out=90,in=180] (0.5,0.5) to[out=0,in=90] (1,0);
\end{tikzpicture}	
=\binom{\lsym}{k}
,\quad
\begin{tikzpicture}[anchorbase,scale=1,rounded corners]
\draw[usual] (0.5,0.35) to (0,0.75)node[left]{$k{-}a$} to (0.5,1.15);
\draw[usual] (0.5,0.35) to (1,0.75)node[right]{$a$} to (0.5,1.15);
\draw[usual] (0.5,1.15) to (0.5,1.5)node[above]{$l$};
\draw[usual] (0.5,0.35) to (0.5,0)node[below]{$k$};
\end{tikzpicture}
=0
\text{ for $k>l$}
,$}
\end{gather*}
together with associativity and coassociativity and some additional relations, gives a 
symmetric ribbon $\F$-linear category $\web[\F]$, which is the main object under study in this paper. Precisely, let $\Rep_{\F}\big(\ogroup[\lsym]\big)$
denote the symmetric ribbon $\F$-linear category of finite dimensional $\ogroupf[\F]$-representations.
We show:

\begin{theorem}\label{T:Main}
There is a fully faithful symmetric ribbon $\F$-linear functor
\begin{gather*}
\web[\F]\to\Rep_{\F}\big(\ogroup[\lsym]\big),
\end{gather*}
sending $k$ to the $k$th exterior power of the vector 
$\ogroupf[\F]$-representation.
This functor induces an equivalence of symmetric ribbon (additive) $\F$-linear categories between:
\begin{enumerate}[label=(\roman*)]

\item The additive idempotent closure of $\web[\F]$.

\item The category of all tilting $\ogroupf[\F]$-representations.

\end{enumerate}
\end{theorem}

With respect to \cref{T:Main} we note:
\begin{enumerate}[label=$\blacktriangleright$]

\item The result generalizes \cite[Section 3]{SaTu-bcd-webs} to positive characteristic and we get a nice diagrammatic calculus.

\item The characteristic of our ground field is essentially arbitrary, except that we cannot work in characteristic two (a not surprising restriction for the orthogonal group).

\item The points (i) and (ii) are the expected
caveats in diagrammatic representation theory. Here we stress that an $\F$-linear equivalence of additive $\F$-linear categories is automatically additive as well, and we will omit the ``additive'' below.

\end{enumerate}

Let $p=\infty$ in case the characteristic of $\F$ is zero.
For a number $\lsym\in\N$ let $\lsym_{i}$ be the \emph{$p$-adic digits}, i.e. numbers $\lsym_{i}\in\{0,\dots,p-1\}$ such that
$N=\sum_{i\in\N}\lsym_{i}p^{i}$.
Finally, we write $\TiltF$ for the category of 
$\ogroupf[\F]$-tilting representations, and we use overline to indicate semisimplifications.
Having established \cref{T:Main}, we then prove:

\begin{theorem}\label{T:O-ss}
There is an equivalence of symmetric ribbon $\F$-linear categories
\begin{gather*}
\semisimple{\TiltF}\cong
\bigboxtimes_{i\in\N}\,\semisimple{\TiltFi},
\end{gather*}
where $\boxtimes$ is Deligne's tensor product.
\end{theorem}

The categories $\semisimple{\TiltFi}$ appearing in 
\cref{T:O-ss} are \emph{Verlinde categories}, 
so \cref{T:O-ss} expresses $\semisimple{\TiltF}$ 
in terms of well-known categories.

\begin{remark}
We expect similar results for the symplectic group instead 
of the orthogonal group. However, the exterior powers are more complicated in this case, e.g. even in characteristic zero they are not simple. As a result of this complication, one needs new arguments to identify the semisimplification in terms of well-known categories.
\end{remark}

%%%%%%%%%%%%%%%%

\subsection{How we prove this}\label{SS:IntroductionHow}

Our proof of \cref{T:Main} is, in broad outline, similar to the proof in Cautis--Kamnitzer--Mor\-rison's paper on type A webs \cite{CaKaMo-webs-skew-howe}. We will now sketch our approach and also indicate what the new steps are.

\begin{remark}
To keep things simple in this introduction, we are not precise with the ground rings and fields: there are no serious issues, but some care needs to be taken. We will give the precise statements and details in the main body of the paper.
\end{remark}

\subsubsection{Webs in type A in characteristic zero (known)}\label{SSS:intro-websAQ}

We first recall a simplified outline of the Cautis--Kamnitzer--Morrison proof when working over $\Q$ (so $q=1$, or as we say later, \emph{dequantized}). (Technical side note: Their paper considers $\mathfrak{sl}_N$ and not $\mathfrak{gl}_N$ representations, but it is easy to adapt their arguments to work for $\mathfrak{gl}_N$.)

They consider the following objects connected with the exterior powers $\Lambda^k:=\Lambda^k(\Q^N)$.
\begin{enumerate}[label=\arabic*.]
\item A generators and relations algebra, denoted $\udot[\Q]^{N}(\mathfrak{gl}_m)$, which is finite dimensional, semi\-simple, and is equipped with an isomorphism
\[
\Phi_m:\udot[\Q]^{N}(\mathfrak{gl}_{m})\xrightarrow{\cong} \End_{\glgroupf[\Q]}\big(\bigoplus_{\substack{(k_1, \dots, k_m) \\ 0\le k_{i}\le N}}\Lambda^{k_1}\otimes \dots \otimes \Lambda^{k_m}\big).
\]

\item A full monoidal subcategory of $\Rep(\glgroupf[\Q])$, denoted $\catstuff{Fund}\big(\glgroupf[\Q]\big)$, with objects $\Lambda^{k_1}\otimes \dots \otimes \Lambda^{k_m}$.
\item A generators and relations monoidal category, denoted $\catstuff{Web}\big(\glgroupf[\Q]\big)$, which has objects $(k_1, \dots, k_m)$, and is equipped with a monoidal functor 
\[
\Gamma:\catstuff{Web}\big(\glgroupf[\Q]\big)\rightarrow \catstuff{Fund}\big(\glgroupf[\Q]\big).
\]
\end{enumerate}
Given the above, their proof that $\Gamma$ is fully faithful proceeds as follows. 
\begin{enumerate}[label=(\roman*)]
\item Using the generators and relations for the algebra $\udot[\Q]^{N}(\mathfrak{gl}_m)$, they construct an algebra homomorphism 
\[
\Psi_m:\udot[\Q]^{N}(\mathfrak{gl}_m)\longrightarrow \End_{\catstuff{Web}(\glgroupf[\Q])}\big(\bigoplus_{\substack{(k_1, \dots, k_m) \\ 0\le k_{i}\le N}} (k_1, \dots, k_m)\big).
\]
Everything is set up so that $\Gamma\circ \Psi_m = \Phi_m$. 

\item One of the key features of $\Psi_m$ is that the Chevalley generators $e_{i}, f_{i}$, $i=1, \dots, m-1$, and their divided powers $e_{i}^{(a)}$ and $f_{i}^{(a)}$, $i=1, \dots, m-1$, map to \emph{ladder} shaped diagrams in the diagrammatic category $\catstuff{Web}\big(\glgroupf[\Q]\big)$. This allows them to give a topological argument, using what they call \emph{ladderization of webs}, to prove that $\Psi_m$ is surjective.

\item The proof that $\Gamma$ is fuly faithful uses nothing more than that for all $m\ge 1$, $\Phi_m$ is an isomorphism, $\Psi_m$ is surjective, and $\Gamma\circ \Psi_m = \Phi_m$. 
\end{enumerate}

\subsubsection{Webs in type A integrally (known, with a mild tweak)}\label{SSS:intro-websAZ}

Amazingly, this is not the only proof that $\Gamma$ is fully faithful. Another proof, which avoids using $\udot[\Q]^{N}(\mathfrak{gl}_m)$ entirely, was found by Elias \cite{El-ladders-clasps}. 
Elias' result also generalized the result in \cite{CaKaMo-webs-skew-howe} from $\Q$ to $\Z$ in the following sense.

There are analogs of the source and target of $\Gamma$ which can be defined over $\Z$, denoted by $\catstuff{Web}_{\Z}\big(\glgroup\big)$ and $\catstuff{Fund}_{\Z}\big(\glgroup\big)$. The former is easy to define, since any $\Q$-linear diagrammatic category such that the coefficients of the relations are in $\Z$ can be defined over $\Z$. The definition of $\catstuff{Fund}_{\Z}\big(\glgroup\big)$ is a bit more complicated, but it uses only standard technology from modular representation theory of reductive algebraic groups, like tilting modules. In \cite{El-ladders-clasps}, Elias carefully studies these categories and actually proves that $\Gamma$ is fully faithful over $\Z$. 

It is a consequence of well-known results about tilting modules that $\catstuff{Fund}_{\Z}\big(\glgroup\big)$ has torsion free homomorphism spaces which are the same rank as the dimension of the analogous homomorphism space over $\Q$. In contrast, for $\catstuff{Web}_{\Z}\big(\glgroup\big)$ it is not at all clear from the definition that the homomorphism spaces are torsion free. However, it follows from $\Gamma$ being an equivalence over $\Z$ that the homomorphism spaces in $\catstuff{Web}_{\Z}\big(\glgroup\big)$ are torsion free and the same rank as the analogous homomorphism spaces over $\Q$.  

Elias' proof is technically difficult, so let us sketch a potential alternative proof.

\begin{remark}\label{R:New}
To the best of our knowledge, the arguments below are new, but were known to experts. Since we do not provide complete details, the present paper still relies on Elias' theorem that $\Gamma$ is fully faithful over $\Z$.
\end{remark}

Cautis--Kamnitzer--Morrison already knew that there was an integral version of the algebra $\udot[\Q]^{N}(\mathfrak{gl}_m)$ \cite[Remark, Section 4.1]{CaKaMo-webs-skew-howe}, which we denote by $\udot[\Z]^{N}(\mathfrak{gl}_m)$. A key difference between $\udot[\Q]^{N}(\mathfrak{gl}_m)$ and $\udot[\Z]^{N}(\mathfrak{gl}_m)$ is that the latter is defined using \emph{higher Serre relations}, while the former only needs Serre relations and the higher Serre relations are consequences. It is natural to try to adapt Cautis--Kamnitzer--Morrison's proof to work over $\Z$.

Defining the $\Z$ analog of $\Phi_m$, and proving that it is an isomorphism is a classic result \cite{Do-tilting-alg-groups}. See \cite[Proposition 4.13]{BrEnAiEtOs-semisimple-tilting} for a modern discussion of this result, which also explains the connection to diagrammatic algebra. Moreover, the ladderization of webs argument \cite[Theorem 5.3.1]{CaKaMo-webs-skew-howe} also works over $\Z$, without any changes. Thus, the proof reduces to constructing a $\Z$ analog of $\Psi_m$. 

We know from working over $\Q$ that, in order to have $\Gamma\circ \Psi_m = \Phi_m$, the generators $e_{i}^{(a)}$ and $f_{i}^{(a)}$ have to go to specific ladder shaped diagrams. So we just need to check these diagrams satisfy the defining relations of $\udot[\Z](\mathfrak{gl}_m)$. However, it is highly nontrivial to directly verify the higher Serre relations for ladder webs. To the best of our knowledge, this calculation does not appear in the literature. Filling in this gap, would give an alternative proof of Elias' result. 

On the other hand, Elias's result that $\Gamma$ is an equivalence over $\Z$ can be used to prove that the higher Serre relations hold for the ladder web diagrams in $\catstuff{Web}_{\Z}\big(\glgroup\big)$, see \cref{L:BackgroundWebsHigherSerre}.

\subsubsection{Webs for orthogonal groups in characteristic zero (known before)}\label{SSS:intro-websOQ}

Following the breakthrough work of Cautis--Kamnitzer--Morrison on webs for $\glgroup$, people immediately began hunting for the generalization to orthogonal groups. The (arguably) first result in this direction is the paper \cite{SaTu-bcd-webs}, which although quite subtle in the quantum case gives the following results over $\Q$. 

\cite{SaTu-bcd-webs} continues to consider exterior powers $\Lambda^k:=\Lambda^k(\Q^N)$, but now as representations of $\ogroupf[\Q]\subset \glgroupf[\Q]$, and studies the following objects related to these exterior powers. 

\begin{enumerate}[label=\arabic*.]
\item A generators and relations algebra, denoted $\udot[\Q]^{N}(\mathfrak{so}_{2m})$, which is finite dimensional, semi\-simple, and is equipped with an isomorphism
\[
\Phi_m:\udot[\Q]^{N}(\mathfrak{so}_{2m})\xrightarrow{\cong} \End_{\ogroupf[\Q]}\big(\bigoplus_{\substack{(k_1, \dots, k_m) \\ 0\le k_{i}\le N}}\Lambda^{k_1}\otimes \dots \otimes \Lambda^{k_m}\big).
\]
The Chevalley generators for this algebra are $e_{i}, f_{i}$, $i=1, \dots, m-1$, and $e_m,f_m$. The algebra $\udot[\Q]^{N}(\mathfrak{gl}_{m})$ maps to $\udot[\Q]^{N}(\mathfrak{so}_{2m})$, sending the generators to Chevalley generators for $i=1, \dots, m-1$. 

\item A full monoidal subcategory of $\Rep\big(\ogroupf[\Q]\big)$, denoted $\catstuff{Fund}\big(\ogroupf[\Q]\big)$, with objects $\Lambda^{k_1}\otimes \dots \otimes \Lambda^{k_m}$.
\item A generators and relations monoidal category, denoted $\catstuff{Web}\big(\ogroupf[\Q]\big)$, which has objects $(k_1, \dots, k_m)$, and is equipped with a monoidal functor 
\[
\Gamma:\catstuff{Web}\big(\ogroupf[\Q]\big)\rightarrow \catstuff{Fund}\big(\ogroupf[\Q]\big).
\]
Restriction induces a monoidal functor $\catstuff{Fund}\big(\glgroupf[\Q]\big)\rightarrow \catstuff{Fund}\big(\ogroupf[\Q]\big)$. This is paralleled in the definition of orthogonal webs \cite[Section 3]{SaTu-bcd-webs}, which are webs for $\glgroupf[\Q]$ along with cups and caps realizing the symmetric form on $\Q^N$ preserved by $\ogroupf[\Q]$. In particular, there is a monoidal functor $\catstuff{Web}\big(\glgroupf[\Q]\big)\rightarrow \catstuff{Web}\big(\ogroupf[\Q]\big)$, and we refer to webs in the image as \emph{type A webs}. (It is not clear from the definition that this functor is faithful, but one can show that it is.)
\end{enumerate}

The proof that $\Gamma$ is fully faithful follows the same outline.

\begin{enumerate}[label=(\roman*)]

\item Using the generators and relations for the algebra $\udot[\Q]^{N}(\mathfrak{so}_{2m})$, \cite{SaTu-bcd-webs} constructs an algebra homomorphism (denoted $\Upsilon$ in \cite[Section 6B]{SaTu-bcd-webs})
\[
\Psi_m:\udot[\Q]^{N}(\mathfrak{so}_{2m})\longrightarrow \End_{\catstuff{Web}(\ogroupf[\Q])}\big(\bigoplus_{\substack{(k_1, \dots, k_m) \\ 0\le k_{i}\le N}} (k_1, \dots, k_m)\big),
\]
such that $\Gamma\circ \Psi_m = \Phi_m$. 

\item This $\Psi_m$ still maps Chevalley generators $e_{i}, f_{i}$, and their divided powers $e_{i}^{(a)}$ and $f_{i}^{(a)}$, to the same ladder shaped type A web diagrams. Moreover, the Chevalley generators $e_m, f_m$, and their divided powers, also map to ladder shaped diagrams. A \emph{ladderization of webs} argument then proves $\Psi_m$ is surjective.

\item The proof that $\Gamma$ is fully faithful uses nothing more than that, for all $m\ge 1$, $\Phi_m$ is an isomorphism, $\Psi_m$ is surjective, and $\Gamma\circ \Psi_m = \Phi_m$. 

\end{enumerate}

To construct $\Psi_{m}$, one needs to verify that the defining relations for $\udot[\Q]^{N}(\mathfrak{so}_{2m})$ are satisfied by the ladder webs corresponding to Chevalley generators. Of course, since the type A webs in $\catstuff{Web}\big(\ogroupf[\Q]\big)$ satisfy type A web relations, and the Chevalley generators $e_{i},f_{i}\in \udot[\Q]^{N}(\mathfrak{so}_{2m})$, for $i=1, \dots, m-1$, satisfy the same relations as in $\udot[\Q]^{N}(\mathfrak{gl}_{m})$, it suffices to check relations involving $e_m,f_m$. The wonderful observation of Sartori, \emph{Sartori's trick}, is that this can be done effortlessly \cite[Section 3C]{SaTu-bcd-webs}, thanks to the topological intuition provided by the diagrammatic description of $\catstuff{Web}\big(\ogroupf[\Q]\big)$.

\subsubsection{Webs for orthogonal groups integrally (new)}\label{SSS:intro-websOZ}

In \cref{S:Diagrams}, we define the orthogonal web category over $\Z$ and establish its connection to the familiar orthogonal web category in \cref{SSS:intro-websOQ} after base changing to $\Q$. We prove \cref{T:Main} by combining various results to establish $\Z$ analogs of everything above. 

\begin{remark}
Up to this point we have just been discussing motivation for the present work (although the ideas discussed in \cref{R:New} were new). For the rest of the introduction, the results are new (excluding the discussion about semisimplification in type A \cref{SSS:intro-sstilt}). The new results are nontrivial, but can still be accomplished efficiently, and can be grouped into two classes.

On the one hand, the technical results like (a) and (c) below do not to our knowledge appear in the literature. For comparison, in type A, the analogous technical results were well-established prior to the study of the corresponding web categories.

On the other hand, many of the arguments in the present paper involve generators and relations checks, or ladderization arguments, the key steps of which have appeared in prior papers. Thus, our main task is to carefully ``piece things together.'' Indeed, a central and new observation of this paper is that this approach works.
\end{remark}

\begin{enumerate}[label=\arabic*.]

\item In \cref{S:Background}, we combine results of \cite{Ta-hyperalgebra-presentation} and \cite{Lu-intro-quantum-groups} to establish that $\udot[\Z]^{N}(\mathfrak{g})$ has a generators and relations presentation.

\item In \cref{S:RepsHowe}, we construct a $\Z$ analog of $\Phi_m:\udot[\Z](\mathfrak{so}_{2m})\rightarrow \catstuff{Fund}_{\Z}\big(\ogroup\big)$ and establish that it is an isomorphism. The existence of such an isomorphism was proven in \cite{AdRy-tilting-howe-positive-char}, but in order to guarantee that $\Phi_m$ is compatible with our web diagrammatics, i.e. $\Gamma\circ \Psi_m = \Phi_m$, we provide more details. 

\item In \cref{SS:RepsTiltingO}, we provide background on tilting modules for the disconnected group $\ogroup$, over $\F$. This is inspired by \cite{AcHaRi-disconnected-groups}. We use these results to define and understand $\catstuff{Fund}_{\F}\big(\ogroup\big)$.

\item In \cref{SS:FunctorDefiningMaps}, we construct $\Gamma:\catstuff{Web}_{\Z}\big(\ogroup\big)\rightarrow \catstuff{Fund}_{\Z}\big(\ogroup\big)$.
We then use Elias' results on type A webs over $\Z$ to prove that the relations for the divided power Chevalley generators in $\udot[\Z]^{N}(\mathfrak{gl}_m)$ are satisfied by the ladder webs in $\catstuff{Web}_{\Z}\big(\glgroup\big)$.

\item In \cref{SS:DiagramOHowe}, we establish the existence of $\Psi_m:\udot[\Z](\mathfrak{so}_{2m})\rightarrow \catstuff{Web}_{\Z}\big(\ogroup\big)$, a $\Z$ analog of $\Psi_m$ from \cref{SSS:intro-websOQ}, by checking that the relations satisfied by divided power Chevalley generators in the presentation of $\udot[\A](\mathfrak{so}_{2m})$ are satisfied by the ladder webs in $\catstuff{Web}_{\Z}\big(\ogroup\big)$. This uses that we know the relations for $e_{i}^{(a)},f_{i}^{(a)}$ are satisfied by type A ladder diagrams, and the observation that Sartori's trick to establish relations between orthogonal ladder diagrams work in an identical manner over $\Z$. 

\item The exact same ladderization argument for orthogonal webs over $\Q$ still works over $\Z$, this establishing surjectivity of $\Psi_m$.

\end{enumerate}

\begin{remark}
Everything above about type A webs has a quantum analog \cite{CaKaMo-webs-skew-howe, El-ladders-clasps, LaTu-gln-webs, Br-q-Schur-category}. However, even in characteristic zero, the quantum analog of \cite{SaTu-bcd-webs} does not describe the monoidal category of representations of quantum $\ogroupf[\Q]$. This is because the orthogonal quantum skew Howe duality in \cite{SaTu-bcd-webs} is for the pair $U_q^{\prime}(O_N)\text{-}U_q(\mathfrak{so}_{2m})$, where $U_q^{\prime}(O_N)$ is not a Drinfeld--Jimbo quantum group. The nonstandard quantum group $U_q^{\prime}(O_N)$ is not a Hopf algebra, but is a quantum symmetric pair, i.e. a coideal algebra inside the Drinfeld--Jimbo quantum group $U_q\big(\glgroup\big)$. 

Later, the web category for the Drinfeld--Jimbo quantum group $U_q\big(\ogroup\big)$ was constructed in \cite{BoWu-type-O-webs}, which proves an equivalence between quantum analogs of $\catstuff{Web}\big(\ogroupf[\Q]\big)$ and also $\catstuff{Fund}\big(\ogroupf[\Q]\big)$. This works only in characteristic zero when $q$ is generic and does not use Howe duality, since the $q$-Howe duality for the pair $U_q(O_N)\text{-}U_q^{\prime}(\mathfrak{so}_{2m})$ has not been studied.
\end{remark}

%%%%%%%%%%%%%%%%%%%%%%%%%%%%%
\subsubsection{Semisimplification in type A (known)}\label{SSS:intro-sstilt}
%%%%%%%%%%%%%%%%%%%%%%%%%%%%%

Let $p$ be a prime. The results above about type A webs over $\Z$ can be base changed to give an equivalence
\[
\catstuff{Web}_{\overline{\mathbb{F}}_{p}}\big(\glgroup\big)\rightarrow \catstuff{Fund}_{\overline{\mathbb{F}}_{p}}\big(\glgroup\big).
\]
Taking the idempotent closure of $\catstuff{Fund}_{\overline{\mathbb{F}}_{p}}\big(\glgroup\big)$, i.e. including all direct summands as objects, yields the category of tilting modules $\catstuff{Tilt}_{\overline{\mathbb{F}}_{p}}\big(\glgroup\big)$. 

As before, we count $p=0$ as $p=\infty$ for the remainder of the introduction.
The category of tilting modules has a unique semisimple monoidal quotient by the negligible ideal. The irreducible objects in the semisimple quotient correspond to the indecomposable tilting modules with nonzero dimension modulo $p$. When $p\ge \lsym$, this category is well-studied under the name \emph{Verlinde category}.

A consequence of main Theorem in \cite{BrEnAiEtOs-semisimple-tilting} is that for any $p$, this semisimplified tilting module category for $\glgroup$ is equivalent to a tensor product of Verlinde categories for $\glgroup[N_{i}]$ where $N_{i}$ is the $p^{i}$th term in the $p$-adic decomposition of $N$. Having established the connection between type A webs and tilting modules, their argument can be summarized as follows.
\begin{enumerate}[label=\arabic*.]

\item The tensor product of Verlinde categories for $\glgroup[N_{i}]$ is identified with a quotient of the colored oriented Brauer category \cite[Lemma 3.3]{BrEnAiEtOs-semisimple-tilting}. The colored oriented Brauer category also maps to the semisimplification of $\catstuff{Web}_{\overline{\mathbb{F}}_{p}}\big(\glgroup\big)$ \cite[Equation 3.4]{BrEnAiEtOs-semisimple-tilting}. 

\item A general argument about semisimplification \cite[Lemma 2.6]{BrEnAiEtOs-semisimple-tilting}, along with some facts about exterior powers in characteristic $p$ \cite[Lemma 3.4]{BrEnAiEtOs-semisimple-tilting}, reduces the cl\-aimed equivalence to showing that the functor from the colored oriented Brauer category to the semisimplification of the type A web category is full.

\item Motivated by a connection between endomorphism algebras of permutation modules for symmetric groups and type A webs, they find a spanning set of what they call chicken foot diagrams in $\catstuff{Web}_{\overline{\mathbb{F}}_{p}}\big(\glgroup\big)$ \cite[Section 4]{BrEnAiEtOs-semisimple-tilting}, which we call \emph{few-to-many-to-few (fmf)} diagrams.

\item They prove fullness by showing that the only fmf diagrams which can survive in the the semisimple quotient are in the image of the functor from the colored oriented Brauer category \cite[Lemma 4.16 and Theorem 4.17]{BrEnAiEtOs-semisimple-tilting}.

\end{enumerate}

%%%%%%%%%%%%%%%%%%%%%%%%%%%%%
\subsubsection{Semisimplification for orthogonal groups (new)}\label{SSS:intro-sstilt2}
%%%%%%%%%%%%%%%%%%%%%%%%%%%%%

When $p>2$, there is an analog of tilting modules for $\ogroup[N]$, see \cref{S:BackgroundOrthogonal} (for the potentially first exposition of this). As before, when $p>N$, the semisimplification of this category is a well-understood Verlinde category for $\ogroup$. For $p<N$, the semisimplification was not studied previously in the literature. 

Upon establishing the connection between orthogonal webs and tilting modules in \cref{T:Main}, we spend \cref{S:Semisimple} proving \cref{T:O-ss}, that the semisimplification of tilting modules of $\catstuff{Tilt}_{\overline{\mathbb{F}}_{p}}\big(\ogroup\big)$ is a tensor product of Verlinde categories for $\ogroup[N_{i}]$, where again $N_{i}$ is the $p^{i}$th term in the p-adic expansion of $N$. Our arguments mirror the argument in \cite{BrEnAiEtOs-semisimple-tilting}.

\begin{enumerate}[label=\arabic*.]

\item We recall the colored Brauer category (no orientation) in \cref{SS:ColorBrauer}. We relate its semi\-simple quotient to the tensor product of $\ogroup[N_{i}]$ Verlinde categories in \cref{L:ss-coloredBrauer}, and to $\catstuff{Web}_{\overline{\mathbb{F}}_{p}}\big(\ogroup\big)$, see \cref{P:DiagramsBrauer} and \cref{L:dim-Lambdapi-mod-p}. 

\item In \cref{SS:SemisimpleBrauer}, we use the exact same general arguments about semisimplification and exterior powers in characteristic $p$ to reduce the desired equivalence to showing that the functor from the colored Brauer category to the semisimplification of the orthogonal web category is full.  

\item Although there is no longer a connection between something like permutation modules of Brauer diagrams and orthogonal webs, we draw inspiration from the type A chicken foot diagrams and in \cref{SS:DiagramsAtInfinity}, we find an analogous spanning set of fmf diagrams in $\catstuff{Web}_{\overline{\mathbb{F}}_{p}}\big(\ogroup\big)$ in \cref{P:DiagramsAtInfinity}.

\item We deduce fullness by using \cref{L:merge-split-negligible} to argue that the only fmf diagrams which can be nonzero in the quotient are in the image of the functor from the colored Brauer category.

\end{enumerate}

All of this is new and requires some careful analysis of tilting modules for orthogonal groups which appears in \cref{SS:RepsTiltingO}.
\newline
%%%%%%%%%%%%%%%%%
%
%%%%%%%%%%%%%%%%%
%% Acknowledgments
%%%%%%%%%%%%%%%%%
%
\noindent\textbf{Acknowledgments.}
We thank Jon Brundan, Kevin Coulembier, Ben Elias, Pavel Etingof, 
Harry Geranios, Jonathan Gruber and Victor Ostrik 
for very helpful discussions,
and the referee for useful comments.
D.T. likes to thank the only good thing going on in their life, a Thanksgiving leftover sandwich. E.B. was supported by the National Science Foundation’s M.S.P.R.F.-2202897. D.T. was sponsored by the ARC Future Fellowship FT230100489.

%%%%%%%%%%%%%%%%

\section{Idempotented algebras and type A webs}\label{S:Background}

%%%%%%%%%%%%%%%%

In this section, we present a mix of well-known material and new, though expected (and well-known to experts), observations such as a presentation of $\udot[\Z]^{N}(\mathfrak{g})$ and \cref{L:DiagramsIntegralTypeA}.

%%%%%%%%%%%%%%%%

\subsection{Idempotent version of divided powers form}\label{SS:BackgroundIdemForm}

%%%%%%%%%%%%%%%%

\begin{notation}\label{N:GeneralA}
In this section (and in some sections below) we let $\Aa=\Z$ and not $\A=\Z[\frac{1}{2},\sqrt{-1}]$
as in e.g. \cref{S:Diagrams}. Note that the field of fractions of $\Aa$ is $\Q$. We can go from $\Aa$ to $\A$ by scalar extension, and this is how we use the results in this section in the rest of the paper.
\end{notation}

Let $\mathfrak{g}$ be a semisimple Lie algebra. 
Write $\Z\Phi$ for the associated root 
lattice and $X$ for the integral weight lattice. 
Fix a choice of simple roots $\Delta$. We 
consider the \emph{universal enveloping algebra} 
$\uplain[\Q]=\uplain[\Q](\mathfrak{g})$, viewed as a 
$\Q$-algebra, which has a 
generators and relations presentation 
called the \emph{Serre presentation}. 
In this section we recall 
a presentation of the idempotented $\Aa$-form 
$\udot[\Aa]$ of $\uplain[\Q]$ which is similar in spirit.

\begin{remark}\label{R:BackgroundKnown}
The main result of this section, 
the presentation in \cref{L:BackgroundDivPowerTwo}, 
is well-known but difficult to find explicitly spelled out.
So we decided to give details how to get it.
\end{remark}

\begin{definition}\label{D:BackgroundDivPower}
The \emph{divided power algebra} $\uplain[\Aa]=\uplain[\Aa](\mathfrak{g})$
(for $\uplain[\Q]$)
is the $\Aa$-subalgebra of $\uplain[\Q]$ 
generated by $e_{\alpha}^{(a)}$, 
$f_{\alpha}^{(b)}$, and $\binom{h_{\alpha}}{c}$ 
for all $\alpha\in\Delta$, and $a,b,c\in\N$.
\end{definition}

Fix formal variables $t$ and $u$ and define elements in $\uplain[\Aa]{}[[t,u]]$ by
\begin{gather*}
e_{\alpha}(u)=\sum_{i\geq 0}e_{\alpha}^{(i)}u^{i},\quad
f_{\alpha}(u)=\sum_{i\geq 0}f_{\alpha}^{(i)}u^{i},\quad
h_{\alpha}(u)=\sum_{i\geq 0}\sbinom{h_{\alpha}}{i}u^{i},
\end{gather*}
and similarly with $t$ in place of $u$.

\begin{lemma}\label{L:BackgroundDivPower}
As an $\Aa$-algebra $\uplain[\Aa]$ has a presentation
with generators $e_{\alpha}^{(a)}$, 
$f_{\alpha}^{(b)}$, and $\binom{h_{\alpha}}{c}$ 
for all $\alpha\in\Delta$, and 
$a,b,c\in\N$, and relations
\begin{gather}\label{Eq:BackgroundTanTrivial}
e_{\alpha}^{(0)}=f_{\alpha}^{(0)} 
=\sbinom{h_{\alpha}}{0}=1
\quad\text{for all $\alpha\in\Delta$},
\\[0.05cm]
\label{Eq:BackgroundTanH}
\left\{
\begin{aligned}
h_{\alpha}(t)h_{\alpha}(u)&=
h_{\alpha}(t+u+tu) 
\quad\text{for all $\alpha\in\Delta$},
\\[0.05cm]
h_{\alpha}(t)h_{\beta}(u)&=
h_{\beta}(u)h_{\alpha}(t)
\quad\text{for all $\alpha,\beta\in\Delta$}, 
\\[0.05cm]
h_{\alpha}(t)e_{\beta}(u)&=
e_{\beta}\big((1+t)^{\alpha^{\vee}(\beta)}u\big)
h_{\alpha}(t)
\quad\text{for all $\alpha,\beta\in\Delta$},
\\[0.05cm]
h_{\alpha}(t)f_{\beta}(u)&=
f_{\beta}\big((1+t)^{-\alpha^{\vee}(\beta)}u\big)h_{\alpha}(t)
\quad\text{for all $\alpha,\beta\in\Delta$},
\end{aligned} 
\right.
\\[0.05cm]
\label{Eq:BackgroundTanEE}
\left\{
e_{\alpha}(t)e_{\alpha}(u)=e_{\alpha}(t+u),
\quad f_{\alpha}(t)f_{\alpha}(u)=f_{\alpha}(t+u)
\quad\text{for all $\alpha\in\Delta$},
\right.
\\[0.05cm]
\label{Eq:BackgroundTanEFOne}
\left\{
e_{\alpha}(t)f_{\alpha}(u)=
f_{\alpha}\left(\tfrac{u}{1+tu}\right)h_{\alpha}(tu)
e_{\alpha}\left(\tfrac{t}{1+tu}\right)
\quad\text{for all $\alpha\in\Delta$},
\right.
\\[0.05cm]
\label{Eq:BackgroundTanEFTwo}
\left\{
e_{\alpha}(t)f_{\beta}(u)=f_{\beta}(u)e_{\alpha}(t)
\quad\text{for all $\alpha\neq\beta\in\Delta$},
\right.
\\[0.05cm]
\label{Eq:BackgroundTanSerre}
\left\{
\begin{aligned}
\sum_{p+r=n}(-1)^{r}e_{\alpha}^{(p)}e_{\beta}^{(m)}e_{\alpha}^{(r)}&=0 \quad\text{whenever $n>-m\cdot\alpha_{i}^{\vee}(\beta)$ 
and $\alpha\neq\beta \in\Delta$},
\\[0.05cm]
\sum_{p+r=n}(-1)^{r}f_{\alpha}^{(p)}f_{\beta}^{(m)}f_{\alpha}^{(r)}&=0 \quad\text{whenever $n>-m\cdot\alpha_{i}^{\vee}(\beta)$ 
and $\alpha\neq\beta\in\Delta$}.
\end{aligned}
\right.
\end{gather}
The relations are interpreted as equalities in 
$\uplain[\Aa]{}[[t,u]]$.
\end{lemma}

\begin{proof}
This is \cite[Corollary 5.2]{Ta-hyperalgebra-presentation}.
\end{proof}

The rather trivial relation \cref{Eq:BackgroundTanTrivial} 
is listed for completeness only.

We define the idempotented (universal enveloping) algebra as in 
\cite[Section 23.1]{Lu-intro-quantum-groups}. Note that this idempotented algebra is not unital.

\begin{definition}\label{D:BackgroundIdempotented}
Given $\wtk,\wtl\in X$, set
\begin{gather*}
{}^{\wtl}\uplain[\Q]^{\wtk}
=\uplain[\Q]\Big/\Big(\sum_{\alpha\in\Delta}\big(h_{\alpha}-\alpha^{\vee}(\wtl)\big)
\uplain[\Q] 
+\sum_{\alpha\in\Delta}\uplain[\Q]\big(h_{\alpha}-\alpha^{\vee}(\wtk)\big)\Big),
\end{gather*}
and let $\one_{\wtk}$ be the image 
of $1$ under the canonical projection $\uplain[\Q]\twoheadrightarrow {^{\wtk}\uplain[\Q]^{\wtk}}$. Define the 
\emph{idempotented enveloping algebra} as 
\begin{gather*}
\udot[\Q]=\bigoplus_{\wtl,\wtk\in X}{^{\wtl}\uplain[\Q]^{\wtk}}.
\end{gather*}
We consider $\udot[\Q]$ as an idempotented $\Q$-algebra with multiplication 
inherited from $\uplain[\Q]$.
\end{definition}

Analogous to \cref{D:BackgroundDivPower} we define:

\begin{definition}\label{D:BackgroundIdempotentedTwo}
The \emph{divided power idempotented algebra} $\udot[\Aa]$ (for $\udot[\Q]$)
is the $\Aa$-subalge\-bra of 
$\udot[\Q]$ generated by 
$e_{\alpha}^{(a)}\one_{\wtk}$, 
$f_{\alpha}^{(a)}\one_{\wtk}$ for all 
$\alpha\in\Delta$, $a\in\N$, and $\wtk\in X$. 
\end{definition}

Let and $\mathfrak{n}=\mathfrak{n}_{+}$ and $\mathfrak{n}_{-}$ be the spans of the positive and negative root spaces.
Let $\cbasis$ denote the 
\emph{canonical basis} for 
$\uplain[\Q]^{+}=\uplain[\Q](\mathfrak{n})$, see e.g. \cite[Part IV]{Lu-intro-quantum-groups}. 
There is an isomorphism 
$\uplain[\Q]^{+}\cong\uplain[\Q]^{-}=\uplain[\Q](\mathfrak{n}_{-})$, 
denoted $u^{+}=u\mapsto u^{-}$, which is 
determined by $e_{\alpha}\mapsto f_{\alpha}$ 
for all $\alpha\in\Delta$. 

\begin{lemma}\label{L:BackgroundIdempotentedTwo}
We have
\begin{gather*}
\udot[\Aa]=
\bigoplus_{b_{1},b_{2}\in\cbasis,\wtk\in X}\
\Aa\cdot b_{1}^{-}\one_{\wtk}b_{2}^{+},
\end{gather*}
and the structure constants for the basis $\cbasis$ are in $\Aa$, i.e. $x\cdot y\in\Aa\cdot\cbasis$ for $x,y\in\cbasis$. 
\end{lemma}

\begin{proof}
This is \cite[Section 23.2.2]{Lu-intro-quantum-groups}.
\end{proof}

We have seen two procedures for extening enveloping algebras: adding idempotents and adding divided powers. The divided power idempotented algebra is defined by adding divided powers to the idempotented algebra. Next, we study adding divided powers first, and then adding idempotents. 

\begin{definition}\label{D:BackgroundIdemRels}
Given $\wtk,\wtl\in X$, set
\begin{gather*}
^{\wtl}{\uplainsecond[\Aa]^{\wtk}}=\uplain[\Aa]\Big/
\Big(\sum_{\alpha\in\Delta,i\in\N}
\big(\sbinom{h_{\alpha}}{i}-\sbinom{\alpha^{\vee}(\wtl)}{i}\big)\uplain[\Aa]
+ 
\sum_{\alpha\in\Delta,i\in\N}
\uplain[\Aa]\big(\sbinom{h_{\alpha}}{i}-\sbinom{\alpha^{\vee}(\wtk)}{i}\big)
\Big)
\end{gather*}
and let $\one_{\wtk}$ be the image of 
$1$ under the canonical projection 
$\uplainsecond[\Aa]\twoheadrightarrow{^{\wtk}\uplainsecond[\Aa]^{\wtk}}$. 
Define the \emph{idempotented divided power algebra} to be
\begin{gather*}
\udotsecond[\Aa]=\bigoplus_{\wtk,\wtl\in X}{}^{\wtl}\uplainsecond[\Aa]^{\wtk}.
\end{gather*}
Similarly as before, 
we consider $\udotsecond[\Aa]$ as an idempotented $\Aa$-algebra with multiplication 
inherited from $\uplainsecond[\Aa]$.
\end{definition}

\begin{lemma}\label{L:BackgroundDivPowerTwo}
As an $\Aa$-algebra $\udotsecond[\Aa]$ has a presentation as an idempotented agebra with generators $e_{\alpha}^{(a)}\one_{\wtk}$ and 
$f^{(a)}_{\alpha}\one_{\wtk}$, for 
all $\wtk\in X$, $\alpha\in\Delta$, and 
$a,b\in\N$, and relations
\begin{gather}
\label{Eq:BackgroundH}
\left\{
\one_{\wtk}\one_{\wtl}
=\delta_{\wtk,\wtl}\one_{\wtk}, 
\quad 
e_{\alpha}^{(0)}\one_{\wtk}=\one_{\wtk}, 
\quad 
f_{\alpha}^{(0)}\one_{\wtk}=\one_{\wtk},
\right.
\\[0.05cm]
\label{Eq:BackgroundEE}
\left\{
e_{\alpha}^{(a)}e_{\alpha}^{(b)}\one_{\wtk}= \sbinom{a+b}{a}e_{\alpha}^{(a+b)}\one_{\wtk}, 
\quad 
f_{\alpha}^{(a)}f_{\alpha}^{(b)}\one_{\wtk}= \sbinom{a+b}{a}f_{\alpha}^{(a+b)}\one_{\wtk}
\quad\text{for all $\alpha\in\Delta$},
\right.
\\[0.05cm]
\label{Eq:BackgroundEFOne}
\left\{
e_{\alpha}^{(a)}f_{\alpha}^{(b)}\one_{\wtk} = \sum_{x\in\N}\sbinom{\alpha^{\vee}(\wtk) +a-b}{x}f_{\alpha}^{(b-x)}e_{\alpha}^{(a-x)}
\text{for all $\alpha\in\Delta$},
\right.	
\\[0.05cm]
\label{Eq:BackgroundEFTwo}
\left\{
e_{\alpha}^{(a)}f_{\beta}^{(b)}\one_{\wtk}= f_{\beta}^{(b)}e_{\alpha}^{(a)}\one_{\wtk}
\quad\text{for all $\alpha\ne\beta\in\Delta$},
\right.
\\[0.05cm]
\label{Eq:BackgroundSerre}
\left\{
\begin{aligned}
\sum_{p+r=n}(-1)^{r}e_{\alpha}^{(p)}
e_{\beta}^{(m)}e_{\alpha}^{(r)}\one_{\wtk}
&=0\quad\text{whenever $n>-m\cdot\alpha_{i}^{\vee}(\beta)$ and 
$\alpha\neq\beta\in\Delta$},
\\
\sum_{p+r=n}(-1)^{r}f_{\alpha}^{(p)}
f_{\beta}^{(m)}f_{\alpha}^{(r)}\one_{\wtk} 
&=0 \quad\text{whenever $n>-m\cdot\alpha_{i}^{\vee}(\beta)$ and 
$\alpha\neq\beta\in\Delta$}.
\end{aligned}
\right.
\end{gather}
\end{lemma}

\begin{proof}
This follows from \cref{L:BackgroundDivPower}. The match between the relations is 
\cref{Eq:BackgroundTanH} 
$\leftrightsquigarrow$ \cref{Eq:BackgroundH},
\cref{Eq:BackgroundTanEE} 
$\leftrightsquigarrow$ \cref{Eq:BackgroundEE},
\cref{Eq:BackgroundTanEFOne} 
$\leftrightsquigarrow$ \cref{Eq:BackgroundEFOne},
\cref{Eq:BackgroundTanEFTwo} 
$\leftrightsquigarrow$ \cref{Eq:BackgroundEFTwo},
\cref{Eq:BackgroundTanSerre} 
$\leftrightsquigarrow$ \cref{Eq:BackgroundSerre}.
\end{proof}

The following gives the promised presentation:

\begin{proposition}\label{P:BackgroundGenRels}
There is an isomorphism of $\Aa$-algebras 
$\eta\colon\udotsecond[\Aa]\rightarrow\udot[\Aa]$ such that
$e_{\alpha}^{(a)}\one_{\wtk}\mapsto e_{\alpha}^{(a)}\one_{\wtk}$ and $f_{\alpha}^{(a)}\one_{\wtk}\mapsto f_{\alpha}^{(a)}\one_{\wtk}$.
\end{proposition}

\begin{proof}
Thanks to \cref{L:BackgroundDivPowerTwo}, $\udotsecond[\Aa]$ has a generators and relations presentation, so the existence of $\eta$ is 
a (omitted) relations check. Surjectivity is then immediate. We show $\eta$ is injective. Let 
$\setstuff{PBW}$ be the $\Aa$-basis 
of $\udot[\Aa]$ from \cite[Section 23.2.1]{Lu-intro-quantum-groups}.
By the same arguments as in \cite{Lu-intro-quantum-groups},
the preimage $\eta^{-1}(\setstuff{PBW})$ is an $\Aa$-spanning set of 
$\udotsecond[\Aa]$. Since $\setstuff{PBW}$ is an $\Aa$-basis, it follows that $\eta^{-1}(\setstuff{PBW})$ is $\Aa$-linearly
independent. Thus, $\eta$ sends an $\Aa$-basis 
to an $\Aa$-basis and is therefore an $\Aa$-isomorphism.
\end{proof}

\begin{notation}\label{N:BackgroundGenRels}
We use $\eta$ from \cref{P:BackgroundGenRels} to identify $\udotsecond[\Aa]$ with $\udot[\Aa]$. Thus, we are justified to only refer to $\udot[\Aa]$.
\end{notation}

For certain results about Howe duality below, we need 
to pass from a semisimple Lie algebra to a reductive Lie algebra, namely the general linear group. 
For $\udot[\Aa](\mathfrak{gl}_{\lsym})$, this combines the definition of $\udot[q](\mathfrak{gl}_{\lsym})$ in
\cite[Section 4.1]{CaKaMo-webs-skew-howe} with the higher Serre relations \eqref{Eq:BackgroundSerre}. A more sophisticated approach could be developed using \cite[Theorem 5.1]{Ta-hyperalgebra-presentation}, but we do not pursue it in this paper.

%%%%%%%%%%%%%%%%

\subsection{A quick recap on type A webs}\label{SS:BackgroundRecapGLWebs}

%%%%%%%%%%%%%%%%

The following is a condensed reformulation of the webs from \cite{CaKaMo-webs-skew-howe}.

\begin{notation}\label{N:DiagramsBasics}
We specify our categorical conventions:

\begin{enumerate}[label=\arabic*.]

\item All of our categories are strict as rigid or pivotal categories, and so are functors between them. 
On the representation theoretical side we silently use the usual strictification theorems, e.g. \cite[VII.2]{MaLa-categories} and \cite[Theorem 2.5]{JoSt-braided-tensor}, to ensure that we are in the strict setting.

\item We denote by $\vcirc$ and $\hcirc$ composition and monoidal structure,
respectively.

\item Objects and morphisms are distinguished by font, for example, 
$\obstuff{K}$ denotes an object and $\morstuff{f}$ denotes a morphism.
For example, we let $\obstuff{K}=(k_{1},\dots,k_{r})=k_{1}\hcirc\dots\hcirc k_{r}$ 
be an object of our web categories. We also use {e.g.} $\obstuff{L}$ and similar 
expressions for objects, having the same meaning.

\item The rigid/pivotal structure is denoted by ${}^{\pivo}$, the monoidal unit by $\wtlnit$, and identity morphisms are denoted by $\idmor$, {e.g.} $\idmor_{\obstuff{K}}$ with the notation as in the previous point.

\item\label{N:DiagramsBasics-e} We read diagrams from bottom to top and left to right as summarized by:
\begin{gather*}
(\idmor\hcirc\morstuff{g})\vcirc(\morstuff{f}\hcirc\idmor) 
=
\begin{tikzpicture}[anchorbase,scale=0.22,tinynodes]
\draw[spinach!35,fill=spinach!35] (-5.5,2) rectangle (-2,0.5);
\draw[spinach!35,fill=spinach!35] (-0.5,3) rectangle (3,4.5);
\draw[thick,densely dotted] (-5.5,2.5) node[left,yshift=-2pt]{$\vcirc$} 
to (3.5,2.5) node[right,yshift=-2pt]{$\vcirc$};
\draw[thick,densely dotted] (-1.25,5) 
node[above,yshift=-2pt]{$\hcirc$} to (-1.25,0) node[below]{$\hcirc$};
\draw[usual] (-5,2) to (-5,3.75) node[right]{$\dots$} to (-5,5);
\draw[usual] (-2.5,2) to (-2.5,5);
\draw[usual] (0,0) to (0,1.25) node[right]{$\dots$} to (0,3);
\draw[usual] (2.5,0) to (2.5,3);
\draw[usual] (-5,0) to (-5,0.25) node[right]{$\dots$} to (-5,0.5);
\draw[usual] (-2.5,0) to (-2.5,0.5);
\draw[usual] (2.5,4.5) to (2.5,5);
\draw[usual] (0,4.5) to (0,4.75) node[right]{$\dots$} to (0,5);
\node at (-3.75,1.0) {$\morstuff{f}$};
\node at (1.25,3.5) {$\morstuff{g}$};
\end{tikzpicture}
=
\begin{tikzpicture}[anchorbase,scale=0.22,tinynodes]
\draw[spinach!35,fill=spinach!35] (-5.5,1.75) rectangle (-2,3.25);
\draw[spinach!35,fill=spinach!35] (-0.5,1.75) rectangle (3,3.25);
\draw[usual] (-5,3.25) to (-5,4.125) node[right]{$\dots$} to (-5,5);
\draw[usual] (-2.5,3.25) to (-2.5,5);
\draw[usual] (0,0) to (0,0.625) node[right]{$\dots$} to (0,1.75);
\draw[usual] (2.5,0) to (2.5,1.75);
\draw[usual] (-5,0) to (-5,0.625) node[right]{$\dots$} to (-5,1.75);
\draw[usual] (-2.5,0) to (-2.5,1.75);
\draw[usual] (2.5,3.25) to (2.5,5);
\draw[usual] (0,3.25) to (0,4.125) node[right]{$\dots$} to (0,5);
\node at (-3.75,2.25) {$\morstuff{f}$};
\node at (1.25,2.25) {$\morstuff{g}$};
\end{tikzpicture}
\\
=
\begin{tikzpicture}[anchorbase,scale=0.22,tinynodes]
\draw[spinach!35,fill=spinach!35] (5.5,2) rectangle (2,0.5);
\draw[spinach!35,fill=spinach!35] (0.5,3) rectangle (-3,4.5);
\draw[thick, densely dotted] (5.5,2.5) 
node[right,yshift=-2pt]{$\vcirc$} to (-3.5,2.5) node[left,yshift=-2pt]{$\vcirc$};
\draw[thick, densely dotted] (1.25,5) 
node[above,yshift=-2pt]{$\hcirc$} to (1.25,0) node[below]{$\hcirc$};
\draw[usual] (2.5,2) to (2.5,3.75) node[right]{$\dots$} to (2.5,5);
\draw[usual] (5,2) to (5,5);
\draw[usual] (0,0) to (0,3);
\draw[usual] (-2.5,0) to (-2.5,1.25) node[right]{$\dots$} to (-2.5,3);
\draw[usual] (2.5,0) to (2.5,0.25) node[right]{$\dots$} to (2.5,0.5);
\draw[usual] (5,0) to (5,0.5);
\draw[usual] (-2.5,4.5) to (-2.5,4.75) node[right]{$\dots$} to (-2.5,5);
\draw[usual] (0,4.5) to (0,5);
\node at (3.75,1.0) {$\morstuff{g}$};
\node at (-1.25,3.5) {$\morstuff{f}$};
\end{tikzpicture}
=(\morstuff{f}\hcirc\idmor)\vcirc(\idmor\hcirc\morstuff{g}).
\end{gather*}

\end{enumerate}	
We specify more notation along the way.
\end{notation}

\begin{definition}\label{D:BackgroundWeb}
Let $\webb[\Aa]$ denote the monoidal 
$\Aa$-linear category $\hcirc$-generated by 
objects $k\in\N$ with $0=\wtlnit$ and $k^{\pivo}=k$, and $\vcirc$-$\hcirc$-generated by morphisms
\begin{gather*}
\begin{tikzpicture}[anchorbase,scale=1]
\draw[usual] (0,0)node[below]{$k$} to (0.5,0.5);
\draw[usual] (1,0)node[below]{$l$} to (0.5,0.5);
\draw[usual] (0.5,0.5) to (0.5,1)node[above]{$k{+}l$};
\end{tikzpicture}
\colon k\hcirc l\to k+l
,\quad
\begin{tikzpicture}[anchorbase,scale=1,yscale=-1]
\draw[usual] (0.5,0.5) to (0,0)node[above]{$k$};
\draw[usual] (0.5,0.5) to (1,0)node[above]{$l$};
\draw[usual] (0.5,1)node[below]{$k{+}l$} to (0.5,0.5);
\end{tikzpicture}
\colon k+l\to k\hcirc l
,\quad
\begin{tikzpicture}[anchorbase,scale=1]
\draw[usual] (0,0)node[below]{$k$} to (1,1)node[above]{$k$};
\draw[usual] (1,0)node[below]{$l$} to (0,1)node[above]{$l$};
\end{tikzpicture}
\colon k\hcirc l\to l\hcirc k
,
\end{gather*}
called \emph{merge, split, and crossing}, subject to the relations 
$\vcirc$-$\hcirc$-generated by 
\begin{gather*}
\begin{tikzpicture}[anchorbase,scale=1]
\draw[usual] (0,0)node[below]{${>}\lsym$} to (0,1.5)node[above]{${>}\lsym$};
\end{tikzpicture}
=0
,\quad
\begin{tikzpicture}[anchorbase,scale=1]
\draw[usual] (0,0)node[below]{$k$} to (1,1);
\draw[usual] (1,0)node[below]{$l$} to (0.5,0.5);
\draw[usual] (2,0)node[below]{$\rsym$} to (1,1);
\draw[usual] (1,1) to (1,1.5)node[above]{$k{+}l{+}m$};
\end{tikzpicture}
=
\begin{tikzpicture}[anchorbase,scale=1]
\draw[usual] (0,0)node[below]{$k$} to (1,1);
\draw[usual] (1,0)node[below]{$l$} to (1.5,0.5);
\draw[usual] (2,0)node[below]{$\rsym$} to (1,1);
\draw[usual] (1,1) to (1,1.5)node[above]{$k{+}l{+}m$};
\end{tikzpicture}
,\quad
\begin{tikzpicture}[anchorbase,scale=1]
\draw[usual] (0,0)node[above]{$k$} to (1,-1);
\draw[usual] (1,0)node[above]{$l$} to (0.5,-0.5);
\draw[usual] (2,0)node[above]{$\rsym$} to (1,-1);
\draw[usual] (1,-1) to (1,-1.5)node[below]{$k{+}l{+}m$};
\end{tikzpicture}
=
\begin{tikzpicture}[anchorbase,scale=1]
\draw[usual] (0,0)node[above]{$k$} to (1,-1);
\draw[usual] (1,0)node[above]{$l$} to (1.5,-0.5);
\draw[usual] (2,0)node[above]{$\rsym$} to (1,-1);
\draw[usual] (1,-1) to (1,-1.5)node[below]{$k{+}l{+}m$};
\end{tikzpicture}
,\\
\begin{tikzpicture}[anchorbase,scale=1,rounded corners]
\draw[usual] (0.5,0.35) to (0,0.75)node[left]{$k$} to (0.5,1.15);
\draw[usual] (0.5,0.35) to (1,0.75)node[right]{$l$} to (0.5,1.15);
\draw[usual] (0.5,1.15) to (0.5,1.5)node[above]{$k{+}l$};
\draw[usual] (0.5,0.35) to (0.5,0)node[below]{$k{+}l$};
\end{tikzpicture}
=\binom{k+l}{k}
\cdot
\begin{tikzpicture}[anchorbase,scale=1]
\draw[usual] (0.5,0)node[below]{$k{+}l$} to (0.5,1.5)node[above]{$k{+}l$};
\end{tikzpicture}
,\quad
\begin{tikzpicture}[anchorbase,scale=1]
\draw[usual] (0,0)node[below]{$k$} to (0.5,0.5) to (0.5,1) to (0,1.5)node[above]{$r$};
\draw[usual] (1,0)node[below]{$l$} to (0.5,0.5) to (0.5,1) to (1,1.5)node[above]{$s$};
\end{tikzpicture}
= \sum_{\substack{a,b\ge 0 \\ k-a+b=r}}
(-1)^{ab}
\begin{tikzpicture}[anchorbase,scale=1]
\draw[usual] (1,0.4)node[left,yshift=-0.05cm]{$b$} to (0,1.1);
\draw[usual] (0,0.4)node[right,yshift=-0.1cm]{$a$} to (1,1.1);
\draw[usual] (0,0)node[below]{$k$} to (0,1.5)node[above]{$r$};
\draw[usual] (1,0)node[below]{$l$} to (1,1.5)node[above]{$s$};
\end{tikzpicture}
.
\end{gather*}
These are called, in order, \emph{exterior relation}, \emph{associativity}, \emph{coassociativity}, \emph{digon removal}, and (signed) \emph{Schur relation}.
The morphisms of $\webb[\Aa]$ are called 
\emph{(integral type A exterior) webs}, while 
\emph{diagrams} are webs obtained from the generating webs without 
taking $\Aa$-linear combinations.
\end{definition}

\begin{lemma}\label{L:BackgroundWebRelations}
The following relations hold in $\webb[\Aa]$.
\begin{enumerate}[label=\arabic*.]

\item The \emph{inverse Schur relations}:
\begin{gather*}
\begin{tikzpicture}[anchorbase,scale=1]
\draw[usual] (1,0)node[below]{$l$} to (0,1)node[above]{$l$};
\draw[usual] (0,0)node[below]{$k$} to (1,1)node[above]{$k$};
\end{tikzpicture}
=
(-1)^{kl}
\sum_{b-a=k-l}
(-1)^{k-b}
\begin{tikzpicture}[anchorbase,scale=1]
\draw[usual] (0,0)node[below]{$k$} to (0,1.5)node[above]{$l$};
\draw[usual] (1,0)node[below]{$l$} to (1,1.5)node[above]{$k$};
\draw[usual] (0,0.4) to node[below]{$b$} (1,0.6);
\draw[usual] (1,0.9) to node[above]{$a$} (0,1.1);
\end{tikzpicture}
.
\end{gather*}

\item \emph{Square switch relations}, also called the $E$-$F$-relation of type A:
\begin{gather*}
\begin{tikzpicture}[anchorbase,scale=1]
\draw[usual] (0,0)node[below]{$k$} to (0,1.5)node[above]{$k'$};
\draw[usual] (1,0)node[below]{$l$} to (1,1.5)node[above]{$l'$};
\draw[usual] (0,0.4) to node[below]{$b$} (1,0.6);
\draw[usual] (1,0.9) to node[above]{$a$} (0,1.1);
\end{tikzpicture}
=
\sum_{x\ge 0}
\binom{k'-l}{x}
\begin{tikzpicture}[anchorbase,scale=1]
\draw[usual] (0,0)node[below]{$k$} to (0,1.5)node[above]{$k'$};
\draw[usual] (1,0)node[below]{$l$} to (1,1.5)node[above]{$l'$};
\draw[usual] (0,0.9) to node[above]{$b'$} (1,1.1);
\draw[usual] (1,0.4) to node[below]{$a'$} (0,0.6);
\end{tikzpicture}
\end{gather*}
where $k' = k-b+a$, $l'= k-a+b$, $a'=a-x$ and $b'= b-x$.
\end{enumerate}
\end{lemma}
\begin{proof}
The (non-signed) version of these relations are \cite[4.26-4.28]{BrEnAiEtOs-semisimple-tilting}. The proof that these non-signed relations are given in \cite[Appendix]{BrEnAiEtOs-semisimple-tilting} 
and can easily be adapted to our signed case.
\end{proof}

\begin{lemma}\label{L:BackgroundWebRibbon}
The category $\webb[\Aa]$ together with the crossings is symmetric. We additionally have the following compatibility between crossings and trivalent vertices
\begin{gather*}
\begin{tikzpicture}[anchorbase,scale=1]
\draw[usual] (0,-1)node[below]{$k$} to (1,0) to (0.5,0.5);
\draw[usual] (1,-1)node[below]{$l$} to (0,0) to (0.5,0.5);
\draw[usual] (0.5,0.5) to (0.5,1)node[above]{$k{+}l$};
\end{tikzpicture}
=
(-1)^{kl}
\begin{tikzpicture}[anchorbase,scale=1]
\draw[usual] (0,0)node[below]{$l$} to (0.5,0.5);
\draw[usual] (1,0)node[below]{$k$} to (0.5,0.5);
\draw[usual] (0.5,0.5) to (0.5,1)node[above]{$k{+}l$};
\end{tikzpicture}
.
\end{gather*}
\end{lemma}

\begin{proof}
This can again be proven by adapting proofs from \cite[Appendix]{BrEnAiEtOs-semisimple-tilting} to our signed case. To show the category is symmetric one must check the braid relations \cite[4.32-4.33]{BrEnAiEtOs-semisimple-tilting} and naturality of the braiding, e.g.
\begin{gather*}
\begin{tikzpicture}[anchorbase,scale=1,rounded corners]
\draw[usual] (0,0)node[below]{$k$} to (0.5,0.5);
\draw[usual] (1,0)node[below]{$l$} to (0.5,0.5);
\draw[usual] (0.5,0.5) to (0.5,1.5)node[above]{$k{+}l$};
\draw[usual] (-1,0)node[below]{$j$} to (-1,0.55) to (2,0.95) to (2,1.5)node[above]{$j$};
\end{tikzpicture}
=
\begin{tikzpicture}[anchorbase,scale=1,rounded corners]
\draw[usual] (0,-0.5)node[below]{$k$} to (0,0);
\draw[usual] (1,-0.5)node[below]{$l$} to (1,0);
\draw[usual] (0,0) to (0.5,0.5);
\draw[usual] (1,0) to (0.5,0.5);
\draw[usual] (0.5,0.5) to (0.5,1)node[above]{$k{+}l$};
\draw[usual] (-1,-0.5)node[below]{$j$} to (-1,0.05) to (2,0.45) to (2,1)node[above]{$j$};
\end{tikzpicture}
,
\end{gather*}
\cite[4.31]{BrEnAiEtOs-semisimple-tilting}. The claimed additional compatibility is \cite[4.30]{BrEnAiEtOs-semisimple-tilting}
\end{proof}

We use the following \emph{ladder diagrams}, indicated 
as a composition of generators for the first one, and their shorthand notations $E_{i}^{(a)}$ and $F_{i}^{(a)}$:
\begin{gather}\label{Eq:BackgroundLadders}
\begin{aligned}
E_{i}^{(a)}\idmor_{\obstuff{K}}&:=
\begin{tikzpicture}[anchorbase,scale=1]
\draw[usual] (0,0)node[below]{$k_{i}$} to (0,1)node[above,xshift=-0.1cm]{$k_{i}{+}a$};
\draw[usual] (1,0)node[below]{$k_{i+1}$} to (1,1)node[above,xshift=0.1cm]{$k_{i+1}{-}a$};
\draw[usual] (0,0.6) to node[above]{$a$} (1,0.4);
\end{tikzpicture}
=
\left(
\begin{tikzpicture}[anchorbase,scale=1]
\draw[usual] (0,0)node[below]{$k_{i}$} to (0.5,0.5);
\draw[usual] (1,0)node[below]{$a$} to (0.5,0.5);
\draw[usual] (0.5,0.5) to (0.5,1)node[above]{$k_{i}{+}a$};
\end{tikzpicture}
\hcirc
\begin{tikzpicture}[anchorbase,scale=1]
\draw[usual] (0,0)node[below]{$k_{i+1}{-}a$} to (0,1)node[above]{$k_{i+1}{-}a$};
\end{tikzpicture}
\right)
\vcirc
\left(
\begin{tikzpicture}[anchorbase,scale=1]
\draw[usual] (0,0)node[below]{$k_{i}$} to (0,1)node[above]{$k_{i}$};
\end{tikzpicture}
\hcirc
\begin{tikzpicture}[anchorbase,scale=1,yscale=-1]
\draw[usual] (0,0)node[above]{$a$} to (0.5,0.5);
\draw[usual] (1,0)node[above]{$k_{i+1}{-}a$} to (0.5,0.5);
\draw[usual] (0.5,0.5) to (0.5,1)node[below]{$k_{i+1}$};
\end{tikzpicture}
\right)
,\\
F_{i}^{(a)}\idmor_{\obstuff{K}}&:=
\begin{tikzpicture}[anchorbase,scale=1]
\draw[usual] (0,0)node[below]{$k_{i}$} to (0,1)node[above,xshift=-0.1cm]{$k_{i}{-}a$};
\draw[usual] (1,0)node[below]{$k_{i+1}$} to (1,1)node[above,xshift=0.1cm]{$k_{i+1}{+}a$};
\draw[usual] (0,0.4) to node[above]{$a$} (1,0.6);
\end{tikzpicture}
,
\end{aligned}
\end{gather}
where one acts on the $i$th and $(i+1)$th entry, and the rest is the identity.

\begin{lemma}\label{L:BackgroundWebsSerre}
The \emph{$E$-$F$ relations} and the \emph{type A Serre relations} hold in $\webb[\Aa]$, that is, for all $a\in\N$:
\begin{gather*}
E_{i}F_{i}\idmor_{\wtk}= 
F_{i}E_{i}\idmor_{\wtk}+
(\wtk_{i-1}+\wtk_{i}-\lsym)\cdot\idmor_{\wtk},
\\
E_{i}F_{j}\idmor_{\wtk}= 
F_{j}E_{i}\idmor_{\wtk}\text{ if }i\ne j,
\\
a!\cdot E_{i}^{(a)}\idmor_{\obstuff{K}}=E_{i}^{a}\idmor_{\obstuff{K}}	
,\quad
E_{i}E_{j}\idmor_{\obstuff{K}}=E_{j}E_{i}\idmor_{\obstuff{K}}\text{ if }|i-j|\neq 1
,\\
2\cdot E_{i}E_{j}E_{i}\idmor_{\obstuff{K}}
=E_{i}^{2}E_{j}\idmor_{\obstuff{K}}+E_{j}E_{i}^{2}\idmor_{\obstuff{K}}
\text{ if }|i-j|=1.
\end{gather*}
Similarly for $F$s.
\end{lemma}

\begin{proof}
See \cite[Lemma 2.2.1 and Proposition 5.2.1]{CaKaMo-webs-skew-howe}.
\end{proof}

\begin{notation}\label{N:DiagramsIdempotentCategory}
We consider $\udot[\Aa](\mathfrak{g})$
as an $\Aa$-linear category with the $\one_{\wtk}$ being the objects. 
In this convention it suffices to indicate one $\one_{\wtk}$ 
in a given expression and we will use this below.
The reader unfamiliar with this is referred to \cite[Section 4.1]{CaKaMo-webs-skew-howe}.
\end{notation}

The following statement is known to experts. However, we are not aware of a proof in the literature so we give one. 

\begin{lemma}\label{L:DiagramsIntegralTypeA}
There is a full $\Aa$-linear functor
\begin{gather*}
\ahowe[\Aa]\colon\udot[\Aa](\mathfrak{gl}_{\rsym})\to\webb[\Aa]
\end{gather*}
such that, in the notation of \cref{Eq:BackgroundLadders},
\begin{gather}\label{Eq:DiagramsCKMFunctor}
e^{(a)}_{i}\one_{\wtk}\mapsto
E_{i}^{(a)}\idmor_{\obstuff{K}}=
\begin{tikzpicture}[anchorbase,scale=1]
\draw[usual] (0,0)node[below]{$k_{i}$} to (0,1)node[above,xshift=-0.1cm]{$k_{i}{+}a$};
\draw[usual] (1,0)node[below]{$k_{i+1}$} to (1,1)node[above,xshift=0.1cm]{$k_{i+1}{-}a$};
\draw[usual] (0,0.6) to node[above]{$a$} (1,0.4);
\end{tikzpicture}
,\quad
f^{(a)}_{i}\one_{\wtk}\mapsto
F_{i}^{(a)}\idmor_{\obstuff{K}}=
\begin{tikzpicture}[anchorbase,scale=1]
\draw[usual] (0,0)node[below]{$k_{i}$} to (0,1)node[above,xshift=-0.1cm]{$k_{i}{-}a$};
\draw[usual] (1,0)node[below]{$k_{i+1}$} to (1,1)node[above,xshift=0.1cm]{$k_{i+1}{+}a$};
\draw[usual] (0,0.4) to node[above]{$a$} (1,0.6);
\end{tikzpicture}
.
\end{gather}
The kernel of this functor is the ideal generated by $\one_{\wtk}$ 
with at least one entry $>\lsym$.
\end{lemma}

\begin{proof}
In this proof we also work with $\Q$, the fraction field of $\Aa$.

Consider an injective $\Aa$-module map 
$f_{\Aa}\colon\Aa^{r}\to\Aa^{s}$, 
which can be viewed as a matrix with $s$ 
rows, $r$ columns, and entries in $\Aa$. 
Write $\setstuff{coker}_{\Aa}$ and 
$\setstuff{coker}_{\Q}$ for the cokernel 
of $f_{\Aa}$ and $\Q\hcirc f_{\Aa}\colon\Q^{r}\to\Q^{s}$, 
respectively. 
The torsion $\Aa$-submodule of $\setstuff{coker}_{\Aa}$ is the set
\begin{gather*}
\setstuff{T}=
\{c\in\setstuff{coker}_{\Aa}|\text{there exists $a\in\Aa$ such that $a\cdot c= 0$}\}.
\end{gather*}

\noindent\textit{Claim in \cref{L:DiagramsIntegralTypeA}.}
We have $\setstuff{coker}_{\Aa}/\setstuff{T}\cong
\Aa\cdot\setstuff{coker}_{\Q}$.

\noindent\textit{Proof of Claim in \cref{L:DiagramsIntegralTypeA}.}
There is a homomorphism 
$g\colon\setstuff{coker}_{\Aa}\to
\mathrm{res}^{\Q}_{\Aa}\setstuff{coker}_{\Q}$, 
defined by $v+\setstuff{im}(f_{\Aa})\mapsto 1\hcirc v+
\setstuff{im}(\Q\hcirc f_{\Aa})$. We denote 
the image by $\Aa\cdot\setstuff{coker}_{\Q}$. 

We claim that $\setstuff{T}\subset\setstuff{ker}(g)$. 
To see this, let $t+\setstuff{im}(f_{\Aa})\in\setstuff{T}$, 
so there is non-zero $a\in\Aa$ such that 
$a\cdot t=f_{\Aa}(v)$. Then modulo 
$\setstuff{im}(\Q\hcirc f_{\Aa})$ we have
\begin{gather*}
g(t+\setstuff{im}(f_{\Aa}))
\equiv 
1\hcirc t\equiv a^{-1}\hcirc a\cdot t
\equiv a^{-1}\hcirc f_{\Aa}(v)
\equiv a^{-1}\cdot\Q\hcirc f_{\Aa}(1\hcirc v)\equiv 0.
\end{gather*}
It follows that there is a $\Aa$-linear map $\overline{g}\colon \setstuff{coker}_{\Aa}/\setstuff{T}\to\Aa\cdot\setstuff{coker}_{\Q}$. 

Since $\Aa$ is a PID, a finitely generated 
and torsion free $\Aa$-module is free. Therefore, 
one can easily argue that 
$\setstuff{coker}_{\Aa}/\setstuff{T}$ and 
$\Aa\cdot\setstuff{coker}_{\Q}$ are free 
$\Aa$-modules. Finally, we observe that 
the $\Q$ span of $\Aa\cdot\setstuff{coker}_{\Q}$ 
is equal to $\setstuff{coker}_{\Q}$, 
so $\rank_{\Aa}\Aa\cdot\setstuff{coker}_{\Q}=\dim_{\Q}\setstuff{coker}_{\Q}=s-r$. Also, $\rank_{\Aa}\setstuff{coker}_{\Aa}/T$ is 
less than or equal to $s-r$. Thus, surjectivity 
of $g$ implies injectivity of $\overline{g}$ 
and we conclude that $\overline{g}$ is 
an isomorphism 
$\setstuff{coker}_{\Aa}/\setstuff{T}\cong
\Aa\cdot\setstuff{coker}_{\Q}$.\qed(Claim)

In \cite[Proposition 5.21]{CaKaMo-webs-skew-howe}, it is shown that the assignments in \cref{Eq:DiagramsCKMFunctor} determine a functor
\begin{gather*}
\ahowe[\Q]\colon\udot[\Q](\glm)\to\webb[\Q].
\end{gather*}
Restricting this functor to $\udot[\Aa]\subset\udot[\Q]$ we 
obtain a functor
\begin{gather*}
\ahowe[\Q]|_{\udot[\Aa]}\colon\udot[\Aa](\glm)\to\webb[\Q].
\end{gather*}
\noindent Let $\Aa\cdot\webb[\Q]$
denote the category obtained by the $\Aa$-span of coefficient-free diagrams in $\webb[\Q]$. The functor $\ahowe[\Q]|_{\udot[\Aa]}$ has image in $\Aa\cdot\webb[\Q]$. 

It follows from \cite[Theorem 2.58]{El-ladders-clasps} that the homomorphism spaces in $\webb[\Aa]$ are free and finitely generated $\Aa$-modules. In particular, the homomorphism spaces are torsion free.

It is immediate from the generators and relations description of $\webb[\Aa]$ that there is a functor
\begin{gather*}
\functorstuff{I}\colon\webb[\Aa]\rightarrow\Aa\cdot\webb[\Q]
\end{gather*}
that sends diagrams to their diagrammatic counterparts in $\Aa\cdot\webb[\Q]$, and therefore is full. `Claim in \cref{L:DiagramsIntegralTypeA}' implies that $\functorstuff{I}$ 
is an equivalence.

The composition
\begin{gather*}
\functorstuff{I}^{-1}\vcirc\ahowe[\Q]|_{\udot[\Aa]}
\colon\udot[\Aa](\glm)\to\Aa\cdot\webb[\Q]\to\webb[\Aa]
\end{gather*}
is the desired functor $\ahowe[\Aa]$.

Finally, in (a $\glm$ version of) \cite[Theorem 5.3.1]{CaKaMo-webs-skew-howe}, the authors argue that $\ahowe[\Q]$ is full by showing any web can be rewritten as a composition of $E^{(a)}$ and $F^{(a)}$ webs. This claim is still true over $\Aa$, and therefore $\ahowe[\Aa]$ is full.
\end{proof}

\cref{L:DiagramsIntegralTypeA} implies that the 
higher type A Serre relations hold in $\webb[\Aa]$. In fact, as far as we are aware, the following is not entirely explicit in \cite{CaKaMo-webs-skew-howe} and related literature:

\begin{lemma}\label{L:BackgroundWebsHigherSerre}
\emph{Higher $E$-$F$ relations} and \emph{higher type A Serre relations} hold in the category $\webb[\Aa]$, that is, for all $a,b\in\N$:
\begin{gather*}
E_{i}^{(a)}F_{i}^{(b)}\idmor_{\wtk}= 
\sum_{x\in\N}\sbinom{\wtk_{i}-\wtk_{i+1}+a-b}{x}
\cdot F_{i}^{(b-x)}E_{i}^{(a-x)}\idmor_{\wtk},
\\
E_{i}^{(a)}F_{j}^{(b)}\idmor_{\wtk}= 
F_{j}^{(a)}E_{i}^{(b)}\idmor_{\wtk}\text{ if }i\ne j,
\\
\sbinom{a+b}{a}\cdot E_{i}^{(a)}E_{i}^{(b)}\idmor_{\obstuff{K}}=E_{i}^{(a+b)}\idmor_{\obstuff{K}}
,\quad
E_{i}^{(a)}E_{j}^{(b)}\idmor_{\obstuff{K}}=E_{j}^{(b)}E_{i}^{(a)}\idmor_{\obstuff{K}}\text{ if }|i-j|\neq 1
,\\
\sum_{p+r=n}E_{i}^{(p)}E_{j}^{(m)}E_{i}^{(r)}\idmor_{\obstuff{K}}
=0\text{ if }|i-j|=1\text{ and }n>m.
\end{gather*}
Similarly for $F$s.
\end{lemma}

\begin{proof}
We know that the relation in 
\cref{Eq:BackgroundSerre} holds in 
$\udot[\Aa](\glm)$ by \cref{P:BackgroundGenRels}. 
The claim then follows from applying the 
functor $\ahowe[\Aa]$ defined in
\cref{L:DiagramsIntegralTypeA} to these relations.
\end{proof}

\begin{example}\label{E:BackgroundWebsHigherSerre}
The relation
$E_{i}E_{i+1}E_{i}=E_{i}^{(2)}E_{i+1}+E_{i+1}E_{i}^{(2)}$, in terms of 
diagrams, is
\begin{gather*}
\begin{tikzpicture}[anchorbase,scale=1]
\draw[usual] (0,0)node[below]{$k_{i}$} to (0,2) node[above,white]{$2$};
\draw[usual] (1,0)node[below]{$k_{i+1}$} to (1,2);
\draw[usual] (2,0)node[below]{$k_{i+2}$} to (2,2);
\draw[usual] (0,0.7) to (1,0.5);
\draw[usual] (0,1.7) to (1,1.5);
\draw[usual] (1,1.0) to (2,0.8);
\end{tikzpicture}	
=
\begin{tikzpicture}[anchorbase,scale=1]
\draw[usual] (0,0)node[below]{$k_{i}$} to (0,2) node[above,white]{$2$};
\draw[usual] (1,0)node[below]{$k_{i+1}$} to (1,2);
\draw[usual] (2,0)node[below]{$k_{i+2}$} to (2,2);
\draw[usual] (0,1.6) to node[above]{$2$} (1,1.4);
\draw[usual] (1,0.6) to (2,0.4);
\end{tikzpicture}
+
\begin{tikzpicture}[anchorbase,scale=1]
\draw[usual] (0,0)node[below]{$k_{i}$} to (0,2) node[above,white]{$2$};
\draw[usual] (1,0)node[below]{$k_{i+1}$} to (1,2);
\draw[usual] (2,0)node[below]{$k_{i+2}$} to (2,2);
\draw[usual] (0,0.6) to node[above]{$2$} (1,0.4);
\draw[usual] (1,1.6) to (2,1.4);
\end{tikzpicture}
\end{gather*}
and $E_{i}^{(a)}E_{j}^{(b)}=E_{j}^{(b)}E_{i}^{(a)}$, for $|i-j|> 1$ follows from far commutativity illustrated in \cref{N:DiagramsBasics}.
\end{example}

%%%%%%%%%%%%%%%%

\section{Orthogonal diagrammatics}\label{S:Diagrams}

%%%%%%%%%%%%%%%%

Throughout, our choices of the ground ring for orthogonal groups are:
\begin{enumerate}[label=(\roman*)]

\item The ring $\A=\Z[\frac{1}{2},\sqrt{-1}]$ for integral results (the appearance of $\frac{1}{2}$ 
is probably not surprising, and see \cref{R:SquareRootMinusOne} for why we need $\sqrt{-1}$). 
Working \emph{integrally} means working over $\A$. 

\item The field $\F$ over $\A$ whenever we need to avoid torsion. We assume that $\F$ is infinite (to avoid using group schemes).

\item The field $\K=\Q(\sqrt{-1})$, the fraction field of $\A$, whenever we need characteristic zero.
Note that $\K$ is a special case of a field $\F$.

\end{enumerate}

Note that we can always scalar extend 
$\A$-linear categories and functors
from $\A$ to $\F$ (and hence, to $\K$). 
Our notation to indicate this will be as in the following example: an $\A$-linear functor $\ifunctor[\A]\colon\webb[\A]\to\web[\A]$ scalar extends to an $\F$-linear functor $\ifunctor[\F]\colon\webb[\F]\to\web[\F]$.

%%%%%%%%%%%%%%%%

\subsection{Integral orthogonal webs}\label{SS:DiagramsWeb}

%%%%%%%%%%%%%%%%

The category of \emph{integral orthogonal (exterior) webs} $\web$ is defined as:

\begin{definition}\label{D:DiagramsWeb}
Let $\web$ denote the pivotal $\A$-linear category $\hcirc$-generated by objects $k\in\N$ with $0=\wtlnit$ and $k^{\pivo}=k$, and $\vcirc$-$\hcirc$-generated by morphisms
\begin{gather*}
\begin{tikzpicture}[anchorbase,scale=1]
\draw[usual] (0,0)node[below]{$k$} to (0.5,0.5);
\draw[usual] (1,0)node[below]{$l$} to (0.5,0.5);
\draw[usual] (0.5,0.5) to (0.5,1)node[above]{$k{+}l$};
\end{tikzpicture}
\colon k\hcirc l\to k+l
,\quad
\begin{tikzpicture}[anchorbase,scale=1,yscale=-1]
\draw[usual] (0.5,0.5) to (0,0)node[above]{$k$};
\draw[usual] (0.5,0.5) to (1,0)node[above]{$l$};
\draw[usual] (0.5,1)node[below]{$k{+}l$} to (0.5,0.5);
\end{tikzpicture}
\colon k+l\to k\hcirc l
,\quad
\begin{tikzpicture}[anchorbase,scale=1]
\draw[usual] (0,0)node[below]{$k$} to (1,1)node[above]{$k$};
\draw[usual] (1,0)node[below]{$l$} to (0,1)node[above]{$l$};
\end{tikzpicture}
\colon k\hcirc l\to l\hcirc k
,
\end{gather*}
called \emph{merge, split, and crossing}, subject to the relations 
$\vcirc$-$\hcirc$-generated by the morphisms depicted below. To state the relations we display the pivotal structure using 
\emph{caps and cups}:
\begin{gather*}
\begin{tikzpicture}[anchorbase,scale=1]
\draw[usual] (0,0)node[below]{$k$} to[out=90,in=180] (0.5,0.5)node[above,yshift=0.4cm]{$\wtlnit$} to[out=0,in=90] (1,0)node[below]{$k$};
\end{tikzpicture}
\colon k\hcirc k\to\wtlnit
,\quad
\begin{tikzpicture}[anchorbase,scale=1]
\draw[usual] (0,0)node[above]{$k$} to[out=270,in=180] (0.5,-0.5)node[below,yshift=-0.4cm]{$\wtlnit$} to[out=0,in=270] (1,0)node[above]{$k$};
\end{tikzpicture}
\colon\wtlnit\to k\hcirc k
.
\end{gather*}
The relations are:
\begin{enumerate}[label=\arabic*.]

\item The (exterior) type A web relations, \emph{exterior relation}, \emph{associativity}, \emph{coassociativity}, \emph{digon removal}, and (signed) \emph{Schur relation} as depicted in \cref{D:BackgroundWeb}.

\end{enumerate}

To state the remaining relations, we recall that we can use the caps and cups to rotate merges and splits, e.g.:
\begin{gather*}
\begin{tikzpicture}[anchorbase,scale=1]
\draw[usual] (0,0)node[below]{$k$} to (0.5,0.5);
\draw[usual] (1,0)node[below]{$k{+}l$} to (0.5,0.5);
\draw[usual] (0.5,0.5) to (0.5,1)node[above]{$l$};
\end{tikzpicture}
=
\begin{tikzpicture}[anchorbase,scale=1]
\draw[usual] (0,0)node[below]{$k$} to[out=90,in=180] (0.5,0.5) to[out=0,in=90] (1,0)node[below]{$k$};
\draw[usual] (2,0)node[below]{$l$} to (2,1)node[above]{$l$};
\end{tikzpicture}
\vcirc
\begin{tikzpicture}[anchorbase,scale=1,yscale=-1]
\draw[usual] (0.5,0.5) to (0,0)node[above]{$k$};
\draw[usual] (0.5,0.5) to (1,0)node[above]{$l$};
\draw[usual] (0.5,1)node[below]{$k{+}l$} to (0.5,0.5);
\draw[usual] (-1,0)node[above]{$k$} to (-1,1)node[below]{$k$};
\end{tikzpicture}
.
\end{gather*}

\begin{enumerate}[resume]

\item \emph{Circle removal and lollipop relations} depicted as
\begin{gather*}
\begin{tikzpicture}[anchorbase,scale=1]
\draw[usual] (0,0)node[left]{$k$} to[out=270,in=180] (0.5,-0.5) to[out=0,in=270] (1,0);
\draw[usual] (0,0) to[out=90,in=180] (0.5,0.5) to[out=0,in=90] (1,0);
\end{tikzpicture}	
=\binom{\lsym}{k}
,\quad
\begin{tikzpicture}[anchorbase,scale=1,rounded corners]
\draw[usual] (0.5,0.35) to (0,0.75)node[left]{$k{-}a$} to (0.5,1.15);
\draw[usual] (0.5,0.35) to (1,0.75)node[right]{$a$} to (0.5,1.15);
\draw[usual] (0.5,1.15) to (0.5,1.5)node[above]{$l$};
\draw[usual] (0.5,0.35) to (0.5,0)node[below]{$k$};
\end{tikzpicture}
=0
\text{ for $k>l$}.
\end{gather*}	

\item \emph{Higher even orthogonal $E$-$F$ relations} depicted as
\begin{gather*}
\begin{tikzpicture}[anchorbase,scale=1]
\draw[usual] (0,0)node[below]{$k$} to (0,1.5)node[above,xshift=-0.25cm]{$k{+}a{-}b$};
\draw[usual] (1,0)node[below]{$l$} to (1,1.5)node[above,xshift=0.25cm]{$l{+}a{-}b$};
\draw[usual] (0,0.5) to[out=315,in=225] node[below]{$a$} (1,0.5);
\draw[usual] (0,1) to[out=45,in=135] node[above]{$b$} (1,1);
\end{tikzpicture}
=
\sum_{x\in\N}\binom{k+l-N+a-b}{x}
\cdot
\begin{tikzpicture}[anchorbase,scale=1]
\draw[usual] (0,0)node[below]{$k$} to (0,1.5)node[above,xshift=-0.25cm]{$k{+}a{-}b$};
\draw[usual] (1,0)node[below]{$l$} to (1,1.5)node[above,xshift=0.25cm]{$l{+}a{-}b$};
\draw[usual] (0,1.05) to[out=315,in=225] node[above]{$a{-}x$} (1,1.05);
\draw[usual] (0,0.45) to[out=45,in=135] node[below]{$b{-}x$} (1,0.45);
\end{tikzpicture}
.
\end{gather*}

\end{enumerate}
The morphisms of $\web$ are called 
\emph{(integral orthogonal exterior) webs}, while 
\emph{diagrams} are webs obtained from the generating webs without 
taking $\A$-linear combinations. We also use the notation for objects as in \cref{N:DiagramsBasics}.
\end{definition}

\begin{example}\label{E:Lollipop}
For $l=0$ and $k$ even the lollipop relation becomes
\begin{gather*}
\begin{tikzpicture}[anchorbase,scale=1]
\draw[usual] (0.5,-1)node[below]{$k$} to (0.5,-0.5);
\draw[usual] (0.5,-0.5) to (0,0) to[out=90,in=180] (0.5,0.5)node[above]{$k/2$} to[out=0,in=90] (1,0) to (0.5,-0.5);
\end{tikzpicture}
=0.
\end{gather*}
The name lollipop relation originates in this picture.
The representation theoretical interpretation of this relation and of the digon removal is a nonsemisimple version of \emph{Schur's lemma} which holds over $\A$. Indeed, 
the reader can compare \cref{L:SimpleTilting} with the usual well-known 
results in the theory, e.g. \cite[Corollary 7.4]{AnPoWe-representation-qalgebras} or \cite[Remark 2.29]{AnTu-tilting}.
\end{example}

Recall $\webb$ as in \cref{D:BackgroundWeb}. The defining relations of $\webb$ are a subset of 
those in $\web$.

\begin{definition}\label{D:Functor-AtoO}
There is a monoidal $\A$-linear functor
\begin{gather*}
\ifunctor\colon\webb\to\web
\end{gather*}
sending $k$ to $k$ and diagrams in $\webb$ to their counterparts in $\web$.
\end{definition}

In particular, all relations listed in \cref{SS:BackgroundRecapGLWebs}
hold in $\web$ as well. We will use this (mostly) silently 
throughout.

\begin{proposition}\label{P:DiagramsFromAToO}
The functor $\ifunctor[\F]$ is faithful.
\end{proposition}
\begin{proof}
Proven in 
\cref{SS:FullyFaithful}.
\end{proof}

\begin{remark}
We will not use Proposition \ref{P:DiagramsFromAToO} until we prove it. For now, we are establishing results needed to prove it.
\end{remark}

Recall the notion of a \emph{ribbon category} from, for example, \cite[Section 8.10]{EtGeNiOs-tensor-categories}.
By definition, any ribbon category is braided monoidal and pivotal.

\begin{lemma}\label{L:DiagramsWebRibbon}
The category $\web$ together with the crossings is symmetric. Moreover, $\web$ together with this symmetric and pivotal structure
is ribbon. We additionally have the \emph{Reidemeister I relations}:
\begin{gather*}
\begin{tikzpicture}[anchorbase,scale=1]
\draw[usual] (1,0) to[out=270,in=0] (0.5,-0.35) to[out=180,in=270] (0,1)node[above]{$k$};
\draw[usual] (0,-1)node[below]{$k$} to[out=90,in=180] (0.5,0.35) to[out=0,in=90] (1,0);
\end{tikzpicture}
=
\begin{tikzpicture}[anchorbase,scale=1]
\draw[usual] (0,-1)node[below]{$k$} to (0,1)node[above]{$k$};
\end{tikzpicture}
=
\begin{tikzpicture}[anchorbase,scale=1,xscale=-1]
\draw[usual] (1,0) to[out=270,in=0] (0.5,-0.35) to[out=180,in=270] (0,1)node[above]{$k$};
\draw[usual] (0,-1)node[below]{$k$} to[out=90,in=180] (0.5,0.35) to[out=0,in=90] (1,0);
\end{tikzpicture}
.
\end{gather*}
(Note that being symmetric also implies the \emph{Reidemeister II and III relations}.)
\end{lemma}

\begin{proof}
Because of \cref{L:BackgroundWebRibbon} and the functor from \cref{D:Functor-AtoO}, we just need to derive relations involving the compatibility of cups and caps with the braiding, e.g. the Reidemeister I relation which can be proven by using 
\begin{gather*}
\left(
\begin{tikzpicture}[anchorbase,scale=1]
\draw[usual] (1,0) to[out=270,in=0] (0.5,-0.35) to[out=180,in=270] (0,1)node[above]{$k$};
\draw[usual] (0,-1)node[below]{$k$} to[out=90,in=180] (0.5,0.35) to[out=0,in=90] (1,0);
\end{tikzpicture}
=
\begin{tikzpicture}[anchorbase,scale=1]
\draw[usual] (0,-1)node[below]{$k$} to (0,1)node[above]{$k$};
\end{tikzpicture}
\right)
\Leftrightarrow
\left(
\begin{tikzpicture}[anchorbase,scale=1]
\draw[usual] (1,-1)node[below]{$k$} to (0,0) to[out=90,in=180] (0.5,0.5) to[out=0,in=90] (1,0) to (0,-1)node[below]{$k$};
\end{tikzpicture}
=
\begin{tikzpicture}[anchorbase,scale=1]
\draw[usual] (0,0)node[below]{$k$} to[out=90,in=180] (0.5,0.5) to[out=0,in=90] (1,0)node[below]{$k$};
\end{tikzpicture}
\right)
,
\end{gather*}
and then by applying \cref{L:BackgroundWebRelations}. Note that the naturality of the braiding with respect to the cups and caps follows from being able to rotate diagrams and the Reidemeister II relation.
\end{proof}

Using the Reidemeister I relation from \cref{L:DiagramsWebRibbon}, we can derive the following relation.

\begin{lemma}\label{L:DiagramsWebRelations}
The \emph{sideways digon relations} hold in $\web$, that is:
\begin{gather*}
\begin{tikzpicture}[anchorbase,scale=1]
\draw[usual] (0,0)node[below]{$k$} to (0.5,0.5) to (0.5,1) to (0,1.5)node[above]{$k$};
\draw[usual] (0.5,0.5) to (0.5,1)  
to[out=90,in=90] (1.25,1) to (1.25,0.75)node[right]{$l$} to (1.25,0.5) to[out=270,in=270] (0.5,0.5);
\end{tikzpicture}
=
\binom{\lsym-k}{l}\cdot
\begin{tikzpicture}[anchorbase,scale=1]
\draw[usual] (0,0)node[below]{$k$} to (0,1.5)node[above]{$k$};
\end{tikzpicture}
.
\end{gather*}
\end{lemma}

\begin{proof}
Similar argument to \cite[Lemma 3.5]{SaTu-bcd-webs}.
\end{proof}

For the sake of completeness, we will now make the connection to the dequantized ($\mathbf{z}=q^{\lsym}$ and $q=-1$ in \cite[Section 3]{SaTu-bcd-webs}) of \emph{Sartori's orthogonal (exterior) webs} from \cite[Section 3]{SaTu-bcd-webs}.

\begin{remark}\label{R:DiagramsWebRelations}
For us \cref{L:DiagramsWebRelations} is a consequence 
of the defining relations. In contrast,
\cite[Definition 3.2]{SaTu-bcd-webs} needs to impose this relation 
due to the lack of pivotality in their setting (the dequantization of 
\cite[Definition 3.2]{SaTu-bcd-webs} is however pivotal).
\end{remark}

\begin{proposition}\label{P:DiagramsWeb}
As a symmetric ribbon $\K$-linear category
$\web[\K]$ is equivalent to 
Sartori's orthogonal web category.
\end{proposition}

\begin{proof}[Proof of \cref{P:DiagramsWeb}]
We write $\sweb[\K]$ to denote the category in \cite[Section 3]{SaTu-bcd-webs}, but treated as a $\K$-linear category via the specialization $\mathbf{z}=q^{\lsym}$ and $q=-1$.

Then \cref{L:DiagramsWebRibbon} and \cref{L:DiagramsWebRelations} imply that
there is a symmetric ribbon $\K$-linear functor $\functorstuff{F}\colon\sweb[\K]\to\web[\K]$ given 
by sending the generators of $\sweb[\K]$ to the diagrams with the same name in $\web[\K]$.

An inverse functor can be defined using \emph{explosion} 
(note that we work over $\K$ in this proof), i.e. we define morphisms 
in $\sweb[\K]$ by
\begin{gather*}
\begin{tikzpicture}[anchorbase,scale=1]
\draw[usual] (0,0)node[below]{$2k$} to[out=90,in=180] (0.5,0.5)node[above]{\phantom{k}} to[out=0,in=90] (1,0)node[below]{$2k$};
\end{tikzpicture}
=
\frac{1}{k!}\cdot
\begin{tikzpicture}[anchorbase,scale=1,yscale=-1]
\draw[usual] (0.5,0.5) to (0,0) to[out=270,in=180] (1.25,-0.75)node[below,yshift=0.08cm]{$\vdots$} node[above]{$1$} to[out=0,in=270] (2.5,0);
\draw[usual] (0.5,0.5) to (1,0) to[out=270,in=180] (1.25,-0.25)node[below]{$1$} to[out=0,in=270] (1.5,0) to (2,0.5);
\draw[usual] (0.5,1)node[below]{$2k$} to (0.5,0.5);
\draw[usual] (2,0.5) to (2.5,0);
\draw[usual] (2,1)node[below]{$2k$} to (2,0.5);
\end{tikzpicture}
\colon 2k\hcirc 2k\to\wtlnit,
\end{gather*}
and similarly for thick cups.
We can then define a backwards functor $\functorstuff{F}^{\prime}\colon\web[\K]\to\sweb[\K]$ by sending 
diagrams to their counterparts in the other category. If $\functorstuff{F}^{\prime}$ is a well-defined symmetric ribbon functor,
then it is immediate that $\functorstuff{F}$ and $\functorstuff{F}^{\prime}$
are inverse functors and we are done.

To see that $\functorstuff{F}^{\prime}$ is well-defined we first observe that 
the lollipop relation is a defining relation in $\sweb[\K]$, so it holds in the image of $\functorstuff{F}^{\prime}$. That the circle removal also 
holds is a well-known calculation similar to \cite[Example 1.5]{RoTu-symmetric-howe}. Finally, that $\functorstuff{F}^{\prime}$ is a symmetric ribbon functor is immediate.
\end{proof}

\begin{remark}\label{R:DiagramsBadNotation}
By \cref{P:DiagramsWeb}, $\web$ is philosophically the $\A$-linear analog of the category from \cite[Section 3]{SaTu-bcd-webs}, but
we changed notation since the notation in \cite[Section 3]{SaTu-bcd-webs} gives the impression that $\web$ does not depend on $\lsym$, but it does.
\end{remark}

%%%%%%%%%%%%%%%%

\subsection{Diagrammatic orthogonal Howe duality}\label{SS:DiagramOHowe}

%%%%%%%%%%%%%%%%

Recall $\udot(\mathfrak{g})$ as in \cref{SS:BackgroundIdemForm}. In this section we will prove the analog of Lemma \ref{L:DiagramsIntegralTypeA}, replacing $\udot(\mathfrak{gl}_{\rsym})\to\webb$ with $\udot(\gso)\to\web$. 

To simplify notation, we continue to use $E_{i}^{(a)}$ and $F_{i}^{(a)}$ to denote diagrams as in \cref{Eq:BackgroundLadders}, 
viewed as diagrams in $\web$. In the same spirit we write
\begin{gather}\label{Eq:DiagramsLadders}
\begin{aligned}
e_{i}^{(a)}\idmor_{\obstuff{K}}&:=
\begin{tikzpicture}[anchorbase,scale=1]
\draw[usual] (0,0)node[below]{$k_{i}$} to (0,1)node[above,xshift=-0.1cm]{$k_{i}{+}a$};
\draw[usual] (1,0)node[below]{$k_{i+1}$} to (1,1)node[above,xshift=0.1cm]{$k_{i+1}{+}a$};
\draw[usual] (0,0.5) to[out=315,in=225] node[below]{$a$} (1,0.5);
\end{tikzpicture}
=
\left(
\begin{tikzpicture}[anchorbase,scale=1]
\draw[usual] (0,0)node[below]{$k_{i}$} to (0.5,0.5);
\draw[usual] (1,0)node[below]{$a$} to (0.5,0.5);
\draw[usual] (0.5,0.5) to (0.5,1)node[above]{$k_{i}{+}a$};
\end{tikzpicture}
\hcirc
\begin{tikzpicture}[anchorbase,scale=1]
\draw[usual] (0,0)node[below]{$a$} to (0.5,0.5);
\draw[usual] (1,0)node[below]{$k_{i+1}$} to (0.5,0.5);
\draw[usual] (0.5,0.5) to (0.5,1)node[above]{$k_{i+1}{+}a$};
\end{tikzpicture}
\right)
\vcirc
\left(
\begin{tikzpicture}[anchorbase,scale=1]
\draw[usual] (0,0)node[below]{$k_{i}$} to (0,1)node[above]{$k_{i}$};
\end{tikzpicture}
\hcirc
\begin{tikzpicture}[anchorbase,scale=1]
\draw[usual] (0,0)node[above]{$a$} to[out=270,in=180] (0.5,-0.5)node[below,yshift=-0.4cm]{$\wtlnit$} to[out=0,in=270] (1,0)node[above]{$a$};
\end{tikzpicture}
\hcirc
\begin{tikzpicture}[anchorbase,scale=1]
\draw[usual] (0,0)node[below]{$k_{i+1}$} to (0,1)node[above]{$k_{i+1}$};
\end{tikzpicture}
\right)
,\\
f_{i}^{(a)}\idmor_{\obstuff{K}}&:=
\begin{tikzpicture}[anchorbase,scale=1]
\draw[usual] (0,0)node[below]{$k_{i}$} to (0,1)node[above,xshift=-0.1cm]{$k_{i}{-}a$};
\draw[usual] (1,0)node[below]{$k_{i+1}$} to (1,1)node[above,xshift=0.1cm]{$k_{i+1}{-}a$};
\draw[usual] (0,0.5) to[out=45,in=135] node[below]{$a$} (1,0.5);
\end{tikzpicture}
.
\end{aligned}
\end{gather}
Here, for clarity, we illustrated how $e_{i}^{(a)}$ is obtained 
from the generating morphisms.

The relations in the following Lemma hold in $\web$. Combined with the relations we proved among the $E_{i}$ and $F_{i}$, these relations are sufficient to prove that there is a functor $\udot[\K](\gso)\to\web[\K]$.

\begin{lemma}\label{L:DiagramsWebsSerre}
We have the following in $\web$.
\begin{enumerate}[label=\arabic*.]

\item The \emph{even orthogonal $E$-$F$ relations}, that is:
\begin{gather*}
e_{i}f_{i}\idmor_{\wtk}= 
f_{i}e_{i}\idmor_{\wtk}+
(\wtk_{i-1}+\wtk_{i}-\lsym)\cdot\idmor_{\wtk},
\\
E_{i}f_{m}\idmor_{\wtk}= 
f_{m}E_{i}\idmor_{\wtk} \quad \text{and} \quad F_{i}e_{m}\idmor_{\wtk}= 
e_{m}F_{i}\idmor_{\wtk}.
\end{gather*}

\item The \emph{even orthogonal Serre relations} hold, that is:
\begin{gather*}
e_{i}E_{j}\idmor_{\obstuff{K}}=E_{j}e_{i}\idmor_{\obstuff{K}}\text{ if }|i-j|\neq 1,
\\
2\cdot e_{i}E_{j}e_{i}\idmor_{\obstuff{K}}
=e_{i}^{2}E_{j}\idmor_{\obstuff{K}}+E_{j}e_{i}^{2}\idmor_{\obstuff{K}}
\text{ if }|i-j|=1,
\\
2\cdot E_{i}e_{j}E_{i}\idmor_{\obstuff{K}}
=E_{i}^{2}e_{j}\idmor_{\obstuff{K}}+E_{j}e_{i}^{2}\idmor_{\obstuff{K}}
\text{ if }|i-j|=1
.
\end{gather*}
Similarly for $f$s and $F$s.

\end{enumerate}
\end{lemma}

\begin{proof}
\textit{(a).} The first relation follows by definition. The others are easy to prove using e.g. associativity.

\textit{(b).} We leave it as an easy exercise. For hints, compare to \cite[Lemmas 3.9 to 3.14]{SaTu-bcd-webs}.
\end{proof}

Using the associativity and bigon relations, it is easy to establish the \emph{orthogonal divided power relations} in $\web$, that is for all $a,b\in\N$:
\[
a!\cdot e_{i}^{(a)}\idmor_{\obstuff{K}}=e_{i}^{a}\idmor_{\obstuff{K}} \quad \text{and} \quad e_{i}^{(a)}e_{i}^{(b)}\idmor_{\obstuff{K}}=\sbinom{a+b}{a}\cdot e_{i}^{(a+b)}\idmor_{\obstuff{K}},
\]
and similarly for $f$s. In order to prove there is a functor $\udot(\gso)\to\web$, it remains to show the following relations hold.

\begin{lemma}\label{L:DiagramsWebsHigherSerre}
We have the following in $\web$.	

\begin{enumerate}[label=\arabic*.]

\item The \emph{higher even orthogonal $E$-$F$ relations}, that is, for all $a,b\in\N$:
\begin{gather*}
e_{i}^{(a)}f_{i}^{(b)}\idmor_{\wtk}= 
\sum_{x\in\N}\sbinom{\wtk_{i-1}+\wtk_{i}-N+a-b}{x}
\cdot f_{i}^{(b-x)}e_{i}^{(a-x)}\idmor_{\wtk}.
\end{gather*}

\item The \emph{higher even orthogonal Serre relations} hold
\begin{gather*}
e_{i}^{(a)}E_{j}^{(b)}\idmor_{\obstuff{K}}=E_{j}^{(b)}e_{i}^{(a)}\idmor_{\obstuff{K}}\text{ if }|i-j|\neq 1,
\\
\sum_{p+r=n}e_{i}^{(p)}E_{j}^{(m)}e_{i}^{(r)}\idmor_{\obstuff{K}}
=0
\text{ if }|i-j|=1\text{ and }n>m,
\\
\sum_{p+r=n}E_{i}^{(p)}e_{j}^{(m)}E_{i}^{(r)}\idmor_{\obstuff{K}}
=0\text{ if }|i-j|=1\text{ and }n>m.
\end{gather*}
Similarly for $f$s and $F$s.

\end{enumerate}
\end{lemma}

\begin{proof}
\textit{(a).} By definition.

\textit{(b).} We use Sartori's trick as 
in \cite[Sections 3 and 4]{SaTu-bcd-webs} to prove all relations
inductively by using sequences of type A relations.
The first relation follows easily from far commutativity in $\web$. It is explained in the proofs of \cite[Lemma 3.13, Lemma 3.14]{SaTu-bcd-webs} that for ladder web diagrams, the remaining two even orthogonal Serre relations are a consequence of the type A Serre relations. Since the type A web relations hold in $\web$, it follows from \cref{L:BackgroundWebsHigherSerre} that the higher type A Serre relations hold in $\web$ as well. Thus, the higher even orthogonal Serre relations follows from a completely analogous calculation to the proofs of \cite[Lemma 3.13, Lemma 3.14]{SaTu-bcd-webs}.
\end{proof}

Let $\hideal$ be the $\vcirc$-$\hcirc$-ideal in 
$\udot(\gso)$ generated by 
$\one_{\wtk}$ with at least one entry $>\lsym$.

\begin{proposition}(Diagrammatic Howe duality.)\label{P:DiagramsHowe}
There is a full $\A$-linear functor
\begin{gather*}
\dhowe\colon\udot(\gso)\to\web
\end{gather*}
such that \cref{Eq:DiagramsCKMFunctor} holds for $i\neq\rsym$ and also, in the notation of \cref{Eq:DiagramsLadders},
\begin{gather*}\label{Eq:DiagramsSTFunctor}
e^{(a)}_{\rsym}\one_{\wtk}\mapsto
e_{\rsym}^{(a)}\idmor_{\obstuff{K}}=
\begin{tikzpicture}[anchorbase,scale=1]
\draw[usual] (0,0)node[below]{$k_{\rsym-1}$} to (0,1)node[above,xshift=-0.15cm]{$k_{\rsym-1}{+}a$};
\draw[usual] (1,0)node[below]{$k_{\rsym}$} to (1,1)node[above,xshift=0.15cm]{$k_{\rsym}{+}a$};
\draw[usual] (0,0.5) to[out=315,in=225] node[below]{$a$} (1,0.5);
\end{tikzpicture}
,\quad
f^{(a)}_{\rsym}\one_{\wtk}\mapsto
f_{\rsym}^{(a)}\idmor_{\obstuff{K}}=
\begin{tikzpicture}[anchorbase,scale=1]
\draw[usual] (0,0)node[below]{$k_{\rsym-1}$} to (0,1)node[above,xshift=-0.15cm]{$k_{\rsym-1}{-}a$};
\draw[usual] (1,0)node[below]{$k_{\rsym}$} to (1,1)node[above,xshift=0.15cm]{$k_{\rsym}{-}a$};
\draw[usual] (0,0.5) to[out=45,in=135] node[below]{$a$} (1,0.5);
\end{tikzpicture}
.
\end{gather*}
The kernel of $\dhowe$ contains 
$\hideal$, and the kernel of $\dhowe[\F]$ is spanned by $\hideal$.
\end{proposition}

\begin{proof}[Proof of \cref{P:DiagramsHowe} excluding the identification of the kernel]
To show existence, it follows from \cref{L:BackgroundDivPowerTwo} that it suffices to check that the $\udot(\gso)$ relations are satisfied by $E_{i}^{(a)}$, $F_{i}^{(a)}$, $e_m^{(i)}$, and $f_m^{(i)}$ in $\web$. The most interesting case is the higher Serre relations, which hold thanks to \cref{L:BackgroundWebsHigherSerre} (interpreted in $\web$) and \cref{L:DiagramsWebsHigherSerre}.

To see fullness we can copy \cite[Proof of Theorem 5.3.1]{CaKaMo-webs-skew-howe} 
as follows. We rewrite
\begin{gather*}
\begin{tikzpicture}[anchorbase,scale=1]
\draw[usual] (0,0)node[below]{$v$} to (0,1)node[above]{$v{-}a$};
\draw[usual] (1,0)node[below]{$w$} to (1,1)node[above]{$w{-}a$};
\draw[usual] (0,0.5) to[out=45,in=135] node[below]{$a$} (1,0.5);
\end{tikzpicture}
=
\begin{tikzpicture}[anchorbase,scale=1]
\draw[usual] (-1,0)node[below]{$v{-}a$} to (-1,1)node[above]{$v{-}a$};
\draw[usual] (0,0)node[below]{$a$} to[out=90,in=180] (0.5,0.5)node[above]{\phantom{k}} to[out=0,in=90] (1,0)node[below]{$a$};
\draw[usual] (2,0)node[below]{$w{-}a$} to (2,1)node[above]{$w{-}a$};
\end{tikzpicture}
\vcirc
\left(
\begin{tikzpicture}[anchorbase,scale=1,yscale=-1]
\draw[usual] (0.5,0.5) to (0,0)node[above]{$v{-}a$};
\draw[usual] (0.5,0.5) to (1,0)node[above]{$a$};
\draw[usual] (0.5,1)node[below]{$v$} to (0.5,0.5);
\end{tikzpicture}
\hcirc
\begin{tikzpicture}[anchorbase,scale=1,yscale=-1]
\draw[usual] (0.5,0.5) to (0,0)node[above]{$a$};
\draw[usual] (0.5,0.5) to (1,0)node[above]{$w{-}a$};
\draw[usual] (0.5,1)node[below]{$w$} to (0.5,0.5);
\end{tikzpicture}
\right)
,
\end{gather*}
and similarly for the $e_{m}^{(a)}\one_{\wtk}$s. Now apply the strategy in 
\cite[Proof of Theorem 5.3.1]{CaKaMo-webs-skew-howe}.

The statement about containing the $\vcirc$-$\hcirc$-ideal generated by $\one_{\wtk}$ with at least one entry $>\lsym$. 
follows by construction.
\end{proof}

The remaining statement 
in \cref{P:DiagramsHowe} concerning identification of the kernel will be proven at the end of \cref{SS:FullyFaithful} below.

\subsection{Redundancy of some relations}\label{SS:Redundancy}

The relations in \cref{D:DiagramsWeb}.(c) may all be redundant. We were able to prove the following two lemmas, showing that some of them are redundant.

\begin{lemma}\label{ex-orthoog-e-f}
For $a=b=1$, the higher even orthogonal $E$-$F$ relations become:
\begin{gather}\label{ef-1label}
\begin{tikzpicture}[anchorbase,scale=1]
\draw[usual] (0,0)node[below]{$k$} to (0,1.5)node[above]{$k$};
\draw[usual] (1,0)node[below]{$l$} to (1,1.5)node[above]{$l$};
\draw[usual] (0,0.5) to[out=315,in=225] node[below]{$1$} (1,0.5);
\draw[usual] (0,1) to[out=45,in=135] node[above]{$1$} (1,1);
\end{tikzpicture}
=
(k+l-N)\cdot
\begin{tikzpicture}[anchorbase,scale=1]
\draw[usual] (0,0)node[below]{$k$} to (0,1.5)node[above]{$k$};
\draw[usual] (1,0)node[below]{$l$} to (1,1.5)node[above]{$l$};
\draw[usual] (0,1.05) to[out=315,in=225] node[above]{$1$} (1,1.05);
\draw[usual] (0,0.45) to[out=45,in=135] node[below]{$1$} (1,0.45);
\end{tikzpicture}
+
\begin{tikzpicture}[anchorbase,scale=1]
\draw[usual] (0,0)node[below]{$k$} to (0,1.5)node[above]{$k$};
\draw[usual] (1,0)node[below]{$l$} to (1,1.5)node[above]{$l$};
\end{tikzpicture}
.
\end{gather}
This relation is a consequence of the other (not higher even orthogonal $E$-$F$) relations.
\end{lemma}

\begin{proof}
To see why, observe that the diagram on the left-hand side of \cref{ef-1label} contains two merge-split subdiagrams. Applying the (signed) Schur relation to each merge-split results in a sum of four diagrams. Applying the Reidemeister II relation to one summand results in the $ef$ diagram on the right hand side of \cref{ef-1label}. Applying circle removal to another summand results in $-N$ times the identity diagram. Applying Reidemeister I and then sideways digon relations to the remaining two summands, results in $k+l$ times the identity diagram.
\end{proof}

\begin{lemma}\label{L:Redundancy}
If one only assumes the higher even orthogonal $E$-$F$ relations for the values $k,l\leq\min\{a,b\}$, then 
they follow in general.
\end{lemma}

\begin{proof}
Below we will suppress coefficients and labels to highlight the main steps. We again use a version of Sartori's trick.

The case $a=b=0$ is trivial. To see how this 
can be verified when $a=b=1$, see \cite[Lemma 3.9]{SaTu-bcd-webs}. A similar argument works when $\min\{a,b\}=1$. For $\min\{a,b\}>1$
we write
\begin{gather*}
\begin{tikzpicture}[anchorbase,scale=1]
\draw[usual] (0,0) to (0,1.5);
\draw[usual] (1,0) to (1,1.5);
\draw[usual] (0,0.5) to[out=315,in=225] (1,0.5);
\draw[usual] (0,1) to[out=45,in=135] (1,1);
\end{tikzpicture}
=
\begin{tikzpicture}[anchorbase,scale=1]
\draw[usual] (0,-0.5) to (0,0) to (0.5,0.5) to (1,0);
\draw[usual] (2,0) to (2.5,0.5) to (3,0) to (3,-0.5);
\draw[usual] (0,2.5) to (0,2) to (0.5,1.5) to (1,2);
\draw[usual] (2,2) to (2.5,1.5) to (3,2) to (3,2.5);
\draw[usual] (0.5,0.5) to (0.5,1.5);
\draw[usual] (2.5,0.5) to (2.5,1.5);
\draw[usual] (1,0) to[out=315,in=225] (2,0);
\draw[usual] (1,2) to[out=45,in=135] (2,2);
\end{tikzpicture}
=
\begin{tikzpicture}[anchorbase,scale=1]
\draw[usual] (0,-0.5) to (0,0);
\draw[usual] (0,2) to (0,2.5);
\draw[usual] (3,-0.5) to (3,0);
\draw[usual] (3,2) to (3,2.5);
\draw[usual] (0,0) to (0,2);
\draw[usual] (3,0) to (3,2);
\draw[usual] (1,0) to(1,0.25) to (0,0.75);
\draw[usual] (2,0) to (2,0.25) to (3,0.75);
\draw[usual] (1,2) to (1,1.75) to (0,1.25);
\draw[usual] (2,2) to(2,1.75) to (3,1.25);
\draw[usual] (1,0) to[out=315,in=225] (2,0);
\draw[usual] (1,2) to[out=45,in=135] (2,2);
\draw[orchid,densely dashed,very thick] (-0.25,0) rectangle (1.25,2);
\draw[orchid,densely dashed,very thick] (1.75,0) rectangle (3.25,2);
\end{tikzpicture}
=\sum\sum\text{coeff}\cdot
\begin{tikzpicture}[anchorbase,scale=1]
\draw[usual] (0,-0.5) to (0,0);
\draw[usual] (0,2) to (0,2.5);
\draw[usual] (3,-0.5) to (3,0);
\draw[usual] (3,2) to (3,2.5);
\draw[usual] (0,0) to (0,2);
\draw[usual] (1,0) to (1,2);
\draw[usual] (2,0) to (2,2);
\draw[usual] (3,0) to (3,2);
\draw[usual] (0,0.25) to (1,0.75);
\draw[usual] (3,0.25) to (2,0.75);
\draw[usual] (0,1.75) to (1,1.25);
\draw[usual] (3,1.75) to (2,1.25);
\draw[usual] (1,0) to[out=315,in=225] (2,0);
\draw[usual] (1,2) to[out=45,in=135] (2,2);
\draw[spinach,densely dashed,very thick] (0.75,-0.5) rectangle (2.25,2.5);
\end{tikzpicture}
,
\end{gather*}
where the third step uses \cref{L:BackgroundWebRelations}.(b). Now, the marked piece in the right-hand diagram is a higher even orthogonal $E$-$F$ relation diagram where $k,l\leq\min\{a,b\}$.
\end{proof}

\begin{remark}
A similar argument as in the proof of \cref{ex-orthoog-e-f} above should work when $\min\{a,b\}=1$. That there is such an argument suggests the higher even orthogonal $E$-$F$ relations may follow from the other orthogonal web relations. For example, we convinced ourselves that the case $a=b=2$ also follows from the other relations, and the argument seems to work in general. However, the calculation was rather complicated, and we were unable to definitively prove that the relations in \cref{D:DiagramsWeb}.(c) are redundant.
\end{remark}

%%%%%%%%%%%%%%%%

\section{Howe's action integrally}\label{S:RepsHowe}

%%%%%%%%%%%%%%%%

Let us denote by $\ext(\placeholder)$ the \emph{exterior algebra}. 

\begin{notation}\label{N:RepsNoWedge}
To avoid clutter we will omit the $\wedge$ and write $xy$ instead of $x\wedge y$.
\end{notation}

We study the free $\A$-module $\ext(V\hcirc\A^{\rsym})$ for $\rsym\geq 0$,
where $V$ is a free $\A$-module of rank $\rank_{\A}V=\lsym$.
Thus, we have two crucial (fixed) numbers $\lsym,\rsym\in\Z_{\geq 1}$:
\begin{enumerate}[label=(\roman*)]

\item $\lsym$ is the rank of the left space in $V\hcirc\A^{\rsym}$, or equivalently as we will see later, the \emph{maximal thickness} 
(of horizontal cuts) of strands in $\web$;

\item $\rsym$ is the rank of the right space in $V\hcirc\A^{\rsym}$, or equivalently as we will see later, the \emph{total thickness} of strands in $\web$.
\end{enumerate}
Here, and throughout, $\ogroup=\ogroup[V]$.

\begin{notation}\label{N:RepsDynkinO}
The notation that we will use 
for $\ogroup$ is adapted to the diagram
\scalebox{0.7}{$\dynkin B{}$}, if $\lsym$ is odd, or \scalebox{0.7}{$\dynkin D{}$}, if $\lsym$ is even,
where we have $\lsym$ nodes.
\end{notation}

%%%%%%%%%%%%%%%%%

\subsection{An $\uplain(\gso)$-action on $\ext(\A^{\rsym})$}

%%%%%%%%%%%%%%%%%

We need some notation:

\begin{notation}\label{N:RepsBox}
We fix the following.
\begin{enumerate}[label=\arabic*.]

\item Write $[i,j]=\{i,i+1,\dots,j\}$ for $i\leq j\in\Z$, and	
$\NBoxm=[1,\lsym]\times[1,\rsym]$. In this notation, $\lsym$ 
indexes rows and $\rsym$ indexes columns in illustrations.

\item Actions are always left actions.

\end{enumerate}
As before, we specify more notation as we go.
\end{notation}

Any basis of $\ext(V\hcirc\A^{\rsym})$ 
is indexed by subsets of 
$[1,\lsym]\times[1,\rsym]$.

Suppose $\lsym=1$, and 
consider the $\A$-module 
$\ext(V\hcirc\A^{\rsym})\cong\ext(\A^{\rsym})$. 
Write $\{x_{1},\dots,x_{\rsym}\}$ for the standard 
basis of $\A^{\rsym}$. Given a subset $S\subset[1,\rsym]$, 
such that $S=\{s_{1},\dots,s_{k}\}$ and 
$s_{1}<\dots<s_{k}$, we write $x_{S}=x_{s_{1}}\dots x_{s_{k}}$, 
and get the well-known lemma:

\begin{lemma}\label{L:RepsABasisExt}
The set $\{x_{S}|S\subset[1,\rsym]\}$ is 
an $\A$-basis of $\ext(\A^{\rsym})$.\qed
\end{lemma} 

As usual, there are \emph{differential operators} 
\begin{gather}\label{Eq:RepsClifford}
\begin{gathered}
x_{i}\colon
\ext(\A^{\rsym})\to\ext(\A^{\rsym})
,\quad 
x_{S}\mapsto x_{i}\acts x_{S}=
\begin{cases} 
(-1)^{|S\cap[1,i-1]|}\cdot x_{S\cup\{i\}} & \text{if }i\notin S,
\\
0 & \text{if }i\in S,
\end{cases}
\\
\partial_{i}\colon\ext(\A^{\rsym})\to\ext(\A^{\rsym})
,\quad 
x_{S}\mapsto\partial_{i}\acts x_{S}=
\begin{cases} 
(-1)^{|S\cap[1,i-1]|}\cdot x_{S\setminus\{i\}} & \text{if }i\in S,
\\
0 & \text{if }i\notin S.
\end{cases}
\end{gathered}
\end{gather}
We have the \emph{Leibniz rule} $x_{i}\circ\partial_{j}+\partial_{j}\circ 
x_{i}=\delta_{i,j}$, as one directly verifies.

\begin{remark}\label{R:RepsClifford}
The operators in \cref{Eq:RepsClifford} generate an 
action of the \emph{Clifford algebra} 
\begin{gather*}
\setstuff{Cl}\big(\A^{\rsym}{\oplus}(\A^{\rsym})^{\ast},
(\placeholder,\placeholder)\big),
\quad\text{where }(x_{i},x_{j}^{\ast})=\delta_{ij},
\end{gather*}
on the spin representation.
\end{remark}

Let $a,b\in \{1, \dots, m\}$. Consider the operators on 
$\ext(\A^{\rsym})$ given by
\begin{gather*}
e_{a,b}^{-}=\partial_{a}\circ\partial_{b}
,\quad
e_{a,b}^{o}=x_{a}\circ\partial_{b}-\delta_{a,b}\tfrac{1}{2}
,\quad
e_{a,b}^{+}=x_{a}\circ x_{b}
.
\end{gather*}
Note that $e_{a,b}^{\pm}=-e_{b,a}^{\pm}$ and $e_{a,a}^{\pm}=0$.

\begin{definition}\label{D:RepsHoweChevalleyGens}
We write
\begin{gather*}
e_{i}=e_{i,i+1}^{o}
,\quad 
f_{i}=e_{i+1,i}^{o}
,\quad
h_{i}=e_{i,i}^{o}-e_{i+1,i+1}^{o}
\quad
\text{for }i\in[1,\rsym-1],
\\
e_{\rsym}=e_{\rsym-1,\rsym}^{+}
,\quad 
f_{\rsym}=e_{\rsym,\rsym-1}^{-}
,\quad
h_{\rsym}=e_{\rsym-1,\rsym-1}^{o}+e_{\rsym,\rsym}^{o}. 
\end{gather*}	
These are called \emph{Howe operators for $\gso$} (with $\gso=\gso(\K)$).
\end{definition}

\begin{remark}\label{R:RepsHoweChevalleyGens}
The name in \cref{D:RepsHoweChevalleyGens} comes from Howe's discussion 
on the skew-sym\-metric FFT of invariant theory in \cite[Section 4]{Ho-perspectives-invariant-theory}.
\end{remark}

\begin{notation}\label{N:RepsDynkinSO}
Our notation for $\gso$ is adapted to the diagram
\scalebox{0.7}{$\dynkin D{}$} with $\rsym$ nodes.

In particular, for later use we specify some root conventions.
Let $\Z\Phi\subset\bigoplus_{i=1}^{\rsym}\Z\frac{\epsilon_{i}}{2}$ 
be the \emph{root lattice} for $\gso$. 
We choose the \emph{simple roots} 
by $\alpha_{i}=\epsilon_{i}-\epsilon_{i+1}$, for $i\in[1,\rsym-1]$, 
and $\alpha_{\rsym}=\epsilon_{\rsym-1}+\epsilon_{\rsym}$.
Let $X(\gso)\subset\bigoplus_{i=1}^{\rsym}\Z\frac{\epsilon_{i}}{2}$ 
be the set of \emph{integral weights} 
for $\gso$. Then $X(\gso)_{+}=\bigoplus_{i=1}^{\rsym}\Z\varpi_{i}$, where
\begin{gather*}
\varpi_{1}=\epsilon_{1},\;
\varpi_{2}=\epsilon_{1}+\epsilon_{2},\;
\dots,\;
\varpi_{\rsym-2}=\epsilon_{1}+\dots+\epsilon_{\rsym-2},\;
\begin{aligned}
\varpi_{\rsym-1}&=\tfrac{1}{2}
(\epsilon_{1}+\epsilon_{2}+\dots+\epsilon_{\rsym-1}-\epsilon_{m})
,\\
\varpi_{\rsym}&=\tfrac{1}{2}
(\epsilon_{1}+\epsilon_{2}+\dots+\epsilon_{\rsym-1}+\epsilon_{m}),
\end{aligned}
\end{gather*}
with $i\in[1,\rsym-2]$.
Let moreover $W=\langle s_{1},\dots, s_{\rsym}\rangle$ be the \emph{Weyl group} of $\gso$ generated by the simple reflections corresponding to our choice of simple roots.
\end{notation}

\begin{lemma}\label{L:RepsSOAction}
The operators 
$e_{i}$, $f_{i}$ and $h_{i}$ act on $\ext(\A^{\rsym})$ by:
\begin{gather*}
e_{i}\acts x_{S} = 
\begin{cases} 
x_{(S\setminus\{i+1\})\cup\{i\}} & \text{if }i+1\in S,i\notin S,
\\
0 & \text{else},
\end{cases}
\quad
f_{i}\acts x_{S}=
\begin{cases} 
x_{(S\setminus\{i\})\cup\{i+1\}} & \text{if }i\in S,i{+}1\notin S,
\\
0 & \text{else},
\end{cases}
\\
h_{i}\acts x_{S}= 
\begin{cases}
x_{S} & \text{if }S\cap\{i,i+1\}=\{i\},
\\
-x_{S} & \text{if }S\cap\{i,i+1\}=\{i+1\},
\\
0 & \text{if }S\cap\{i,i+1\}=\emptyset\text{ or }\{i,i+1\},
\end{cases}
\\
e_{\rsym}\acts x_{S}=
\begin{cases} 
x_{S\cup\{m-1,m\}} & \text{if }m-1,m\notin S,
\\
0 & \text{else},
\end{cases}
\quad
f_{\rsym}\acts x_{S}=
\begin{cases} 
x_{S\setminus\{m-1,m\}} & \text{if }m-1,m\in S,
\\
0 & \text{else},
\end{cases}
\\
h_{\rsym}\acts x_{S} = 
\begin{cases}
x_{S} & \text{if }S\cap\{m, m+1\}=\{m, m+1\},
\\
-x_{S} & \text{if }S\cap\{m,m+1\}=\emptyset,
\\
0 & \text{if }S\cap\{m,m+1\}=\{m\}\text{ or }\{m+1\}.
\end{cases}
\end{gather*}
In the first two displays we have $i\in[1,\rsym-1]$.
\end{lemma}

\begin{proof}
Following for example \cite{BoTu-symplectic-howe}, a nice way to think 
about the action is as follows. First, we imagine a row with $\rsym$ 
nodes, where each node is empty or filled with one \emph{dot}. These correspond to 
$S\subset[1,\rsym]$ by, for example,
\begin{gather*}
\ytableausetup{centertableaux,boxsize=0.62cm}
S=\{1,5,6=\rsym\}
\leftrightsquigarrow
\begin{ytableau}
*(magenta!50)\bullet & \phantom{a} & \phantom{a} 
& \phantom{a} & *(magenta!50)\bullet & *(magenta!50)\bullet
\end{ytableau}
\,.
\end{gather*}
In this notation, the operator $x_{i}$ adds a dot in the $i$th box, or acts as zero if there is a dot in the box, up to a factor of $-1$ for each dot to the left of the $i$th box. The operator $\partial_{i}$ removes a dot from the $i$th box, or acts by zero if the box is empty, up to a factor of $-1$ for each dot to the left of the $i$th box. Using this notation the lemma is easy to verify, and we only give an example.
For $\rsym=6$ and $S=\{1,3,4\}$ we have
\begin{gather*}
e_{\rsym}=(x_{\rsym-1}\vcirc
x_{\rsym})\acts
x_{1}x_{3}x_{4}
=x_{1}x_{3}x_{4}x_{5}x_{6}
\\
\leftrightsquigarrow
e_{\rsym}\acts
\begin{ytableau}
*(magenta!50)\bullet & \phantom{a} & *(magenta!50)\bullet 
&*(magenta!50)\bullet & \phantom{a} & \phantom{a}
\end{ytableau}
=
\begin{ytableau}
*(magenta!50)\bullet & \phantom{a} & *(magenta!50)\bullet 
&*(magenta!50)\bullet & *(magenta!50)\bullet & *(magenta!50)\bullet
\end{ytableau}
\,.
\end{gather*}
All other necessary calculations are similar.
\end{proof}

The graphical notation in the above proof is called a \emph{dot diagram}.
We will use these throughout and identify $S$ with such diagrams.

\begin{example}\label{E:RepsDotDiagrams}
Here is an explicit example:
\begin{gather}\label{Eq:RepsDotActionF}
f_{1}\acts x_{\{1,5,6\}}=x_{\{2,5,6\}}
\leftrightsquigarrow
f_{1}\acts
\begin{ytableau}
*(magenta!50)\bullet & \phantom{a} & \phantom{a} 
& \phantom{a} & *(magenta!50)\bullet & *(magenta!50)\bullet
\end{ytableau}
=
\begin{ytableau}
\phantom{a} & *(magenta!50)\bullet & \phantom{a} 
& \phantom{a} & *(magenta!50)\bullet & *(magenta!50)\bullet
\end{ytableau}
\,,
\end{gather}
which is a graphical version of the formulas in \cref{L:RepsSOAction}.
\end{example}

The point is that the action in this notation is visually easier to remember. Explicitly, the action from \cref{L:RepsSOAction} in this dot diagram notation is as follows.
\begin{enumerate}[label=(\roman*)]

\item For $i\in[1,\rsym-1]$ the operators 
$e_{i}$ and $f_{i}$ move dots rightwards or leftwards, if possible, 
and annihilate the diagram otherwise. See \cref{Eq:RepsDotActionF}.

\item The operators $e_{m}$ and $f_{m}$ add or remove dots in the final two 
columns if possible and annihilate the diagram otherwise. For example:
\begin{gather*}
f_{6=\rsym}\acts x_{\{1,5,6\}}=x_{\{1\}}
\leftrightsquigarrow
f_{6}\acts
\begin{ytableau}
*(magenta!50)\bullet & \phantom{a} & \phantom{a} 
& \phantom{a} & *(magenta!50)\bullet & *(magenta!50)\bullet
\end{ytableau}
=
\begin{ytableau}
*(magenta!50)\bullet & \phantom{a} & \phantom{a} 
& \phantom{a} & \phantom{a} & \phantom{a}
\end{ytableau}
\,.
\end{gather*}

\item The $h_{k}$ operators essentially only add signs in case they 
find dots in certain spots. 

\end{enumerate}

More generally, from \cref{SS:TwoBases} below we will use dot diagrams for $S\subset\NBoxm$, where we have 
$\lsym$ rows. For example,
\begin{gather*}
\begin{tikzpicture}[anchorbase]
\draw[thick,->] (0,0) to (0,-0.5)node[left]{$\lsym$} to (0,-1);
\draw[thick,->] (0,0) to (0.5,0)node[above]{$\rsym$} to (1,0);
\end{tikzpicture}
\quad
\begin{ytableau}
*(magenta!50)\bullet & *(magenta!50)\bullet & \phantom{a} 
& \phantom{a} & \phantom{a} & \phantom{a}
\\
*(magenta!50)\bullet & \phantom{a} & \phantom{a} 
& \phantom{a} & \phantom{a} & *(magenta!50)\bullet
\\
*(magenta!50)\bullet & \phantom{a} & \phantom{a} 
& \phantom{a} & *(magenta!50)\bullet & \phantom{a}
\end{ytableau}
\leftrightsquigarrow
S=\{(1,1),(2,1),(3,1),(1,2),(3,5),(2,6)\}.
\end{gather*}

One can check that $W$ acts by even-signed 
permutations on $\oplus_{i=1}^{\rsym}\Z\frac{\epsilon_{i}}{2}$, preserving $X(\gso)$ set-wise, where the simple 
reflection generators $s_{1},\dots,s_{\rsym-1}$ 
permute the $\epsilon_{i}$ and $s_{\rsym}$ 
scales the $\epsilon_{\rsym-1}$ and 
$\epsilon_{\rsym}$ coordinates by $-1$. The longest element $w_{0}\in W$ acts by
\begin{gather*}
w_{0} (a_{1},\dots,a_{\rsym-1}, a_{\rsym}) = 
\begin{cases} 
(-a_{1},\dots,-a_{\rsym-1}, -a_{\rsym}) & \text{if $\rsym$ is even} 
\\
(-a_{1},\dots, -a_{\rsym-1}, a_{\rsym}) & \text{if $\rsym$ is odd}. 
\end{cases}
\end{gather*}
The action of $W$ by even-signed 
permutations on
$\{\frac{1}{2}(\pm\epsilon_{1},\pm\epsilon_{2},\dots,\pm\epsilon_{\rsym})\}$
has two orbits: the vectors with an odd 
number of signs and the vectors with an 
even number of signs. These orbits are the weight spaces of
the $\uplain(\gso)$-representations
$\coweyl[\varpi_{\rsym-1}]$ and 
$\coweyl[\varpi_{\rsym}]$, respectively. The (dual) Weyl representations 
$\weyl[\placeholder]$, $\coweyl[\placeholder]$ are recalled in \cref{S:BackgroundOrthogonal} below.

\begin{remark}\label{R:is-minuscule}
The $\uplain(\gso)$-representations $\coweyl[\varpi_{\rsym-1}]$ and $\coweyl[\varpi_{\rsym}]$ are \emph{minuscule}. It follows that $\coweyl[\varpi_{\rsym-1}]\cong\weyl[\varpi_{\rsym-1}]$ and $\coweyl[\varpi_{\rsym}]\cong \weyl[\varpi_{\rsym}]$.
\end{remark} 

Using the convention that $\varpi_{0}=0$, we see that we have
$\wt_{\gso}(x_{i})=-\varpi_{i-1}+\varpi_{i}-\varpi_{\rsym}$, 
$\wt_{\gso}(x_{\rsym-1})=-\varpi_{m-3}+\varpi_{m-2}$ 
and $\wt_{\gso}(x_{\rsym})=-\varpi_{\rsym-1}$.
In terms of the $\epsilon_{i}$ basis we have
$\wt_{\gso}(x_{i})=\frac{1}{2}(-1,-1,\dots,-1,1,-1\dots,-1,-1)$
with the $1$ in the $i$th entry.
More generally, if $S\subset[1,\rsym]$, 
then $\wt_{\gso}(x_{S})=\sum_{i\in S}\wt_{\gso}(x_{i})$, and in the $\epsilon_{i}$ basis we have $1$ in the $i$th entry for all $i\in S$.
In terms of dot diagrams this reads:
\begin{gather*}
\wt_{\gso}(x_{S})=\tfrac{1}{2}
\big(d_{1}(S),d_{2}(S),\dots,d_{\rsym}(S)\big),
\end{gather*}
where $d_{i}(S)$ is $1$, if the $i$th node for the dot diagram for $S$ 
contains a dot, and $-1$ otherwise.

\begin{example}\label{E:RepsIntegralSOAction}
We get
\begin{gather*}
x_{S}=
\begin{ytableau}
*(magenta!50)\bullet & \phantom{a} & *(magenta!50)\bullet 
&*(magenta!50)\bullet & \phantom{a} & \phantom{a}
\end{ytableau}
\rightsquigarrow
\wt_{\gso}(x_{S})=\tfrac{1}{2}(1,-1,1,1,-1,-1),
\end{gather*}
where $\rsym=6$ and $S=\{1,3,4\}$.
\end{example}

Recall the $\uplain(\gso)$-action on 
$\ext(\A^{\rsym})$ from \cref{L:RepsSOAction}.
Over $\C$ the following is classical, see for example \cite[Section 4.3]{Ho-perspectives-invariant-theory} and the discussion leading to that section.

\begin{proposition}\label{P:RepsIntegralSOAction}
The operators $e_{i}$, $f_{i}$ and $h_{i}$, for $i\in[1,\rsym]$, 
give rise to an $\uplain(\gso)$-action on $\ext(\A^{\rsym})$	
such that the $i$th $\mathfrak{sl}_{2}$-triple $(e_{i},f_{i},h_{i})$
corresponds to the simple root $\alpha_{i}$ and $h_{i}=\alpha_{i}^{\vee}$.
Moreover, 
\begin{gather*}
\ext(\A^{\rsym})\cong\coweyl[\varpi_{\rsym-1}]
\oplus\coweyl[\varpi_{\rsym}]
\end{gather*}
as $\uplain(\gso)$-representations.
\end{proposition}

\begin{proof}
As recalled in \cref{S:BackgroundOrthogonal} below, 
the $\uplain(\gso)$-representation $\coweyl[\lambda]$ specializes 
to the simple $\uplain[\K](\gso)$-representation $\ksimple[\lambda]$.

Using \cref{L:RepsSOAction}, one can easily verify that $\K\hcirc_{\A}\ext(\A^{\rsym})\cong\ksimple[\varpi_{\rsym-1}]
\oplus\ksimple[\varpi_{\rsym}]$, as a $\uplain[\K](\gso)$-representation, such that the operators $e_{i}$, $f_{i}$, $h_{i}$ correspond to $\alpha_{i}$. Moreover, the operators $e_{i}$, 
$f_{i}$, $h_{i}$ preserve the lattice $\ext(\A^{\rsym})$ and the higher divided powers, $e_{i}^{(d)}$ and $f_{i}^{(d)}$, act as zero. Thus, $\uplain(\gso)$ acts on $\ext(\A^{\rsym})$. Analyzing \cref{L:RepsSOAction}, it is easy to see that $\ext(\A^{\rsym})$ is generated over $\uplain(\gso)$ by the highest weight vectors $x_{1}x_{2}\dots x_{\rsym-1}$ and $x_{1}x_{2}\dots x_{\rsym-1}x_{\rsym}$. Indeed, recall that $f_{i}$ acts by moving the $i$th dot to the $i+1$st box (if the $i+1$st box is empty), e.g.
\begin{gather*}
f_{3}\acts
\begin{ytableau}
*(magenta!50)\bullet & *(magenta!50)\bullet & *(magenta!50)\bullet 
& \phantom{a} & *(magenta!50)\bullet & \phantom{a}
\end{ytableau}
=
\begin{ytableau}
*(magenta!50)\bullet & *(magenta!50)\bullet & \phantom{a} 
& *(magenta!50)\bullet & *(magenta!50)\bullet & \phantom{a}
\end{ytableau}
\,,
\end{gather*}
and $f_m$ acts by removing the last two dots (if the last two boxes contain dots), e.g.
\begin{gather*}
f_{6}\acts
\begin{ytableau}
*(magenta!50)\bullet & *(magenta!50)\bullet & *(magenta!50)\bullet 
& \phantom{a} & *(magenta!50)\bullet & *(magenta!50)\bullet
\end{ytableau}
=
\begin{ytableau}
*(magenta!50)\bullet & *(magenta!50)\bullet & *(magenta!50)\bullet 
& \phantom{a} & \phantom{a} & \phantom{a}
\end{ytableau}
\,.
\end{gather*}
Using this one can see that all dot configurations 
can be generated from the vectors $x_{1}x_{2}\dots x_{\rsym-1}$ and $x_{1}x_{2}\dots x_{\rsym-1}x_{\rsym}$.
\end{proof}

%%%%%%%%%%%%%%%

\subsection{Two bases for $\ext(V\hcirc\A^{\rsym})$}\label{SS:TwoBases}

%%%%%%%%%%%%%%%

We now consider the general case 
with arbitrary (but fixed) $\lsym\in\Z_{\geq 1}$, and 
a free $\A$-module $V$ of dimension $\lsym$. 

\begin{notation}\label{N:RepsDistActsTensor}
Throughout, we let the group $\ogroup$ act diagonally on tensor products, as usual.
Moreover, we consider $\uplain(\gso)$ as a Hopf algebra by extending
$\Delta(x)=1\hcirc x+x\hcirc 1$ for all $x\in\gso$.
\end{notation} 

With the structure from \cref{N:RepsDistActsTensor} it follows from \cref{P:RepsIntegralSOAction} that $\uplain(\gso)$ also 
acts on $\ext(\A^{\rsym})^{\hcirc\lsym}\cong\ext(V\hcirc\A^{\rsym})$. Similarly, $\ogroup$ acts on $\ext(V)^{\hcirc\rsym}\cong\ext(V\hcirc\A^{\rsym})$. We will use these two actions below.

\begin{lemma}\label{L:has-Weyl-filt-for-som}
The $\uplain(\gso)$-representation $\ext(V\hcirc\A^{\rsym})$ has a filtration by Weyl representations, and a filtration by dual Weyl representations. 
\end{lemma}

\begin{proof}
Directly from \cref{P:RepsIntegralSOAction}, \cref{R:is-minuscule} and \cref{R:Kaneda}.
\end{proof}

The following classical fact, which can be compared with \cite[Section 4.3]{Ho-perspectives-invariant-theory} again, is fundamental to what follows. 
We prove this result in \cref{SS:FunctorDefiningMaps} by identifying the action of $e_{i}$ and $f_{i}$ as compositions of tensor products of $\ogroup$ intertwiners.

\begin{proposition}(Howe's actions integrally.)\label{P:actions-commute}
The actions of $\ogroup$ and $\uplain(\gso)$ on $\ext(V\hcirc\A^{\rsym})$ commute. Thus, $\ext(V\hcirc\A^{\rsym})$ is an 
$\ogroup$-$\uplain(\gso)^{op}$-birepresentation.
\end{proposition}

\begin{remark}
We have two left actions, and thus we need an ${}^{op}$ in \cref{P:actions-commute}.
\end{remark}

We now specify explicit isomorphisms $\ext(V)^{\hcirc\rsym}\cong\ext(V\hcirc\A^{\rsym})$ and $\ext(\A^{\rsym})^{\hcirc\lsym}\cong\ext(V\hcirc\A^{\rsym})$, see \cref{L:RepsFixIso} below.

Write $V=\bigoplus_{i=1}^{\lsym}\A\cdot v_{i}$. Given, $(i,j)\in\NBoxm$, write
\begin{gather*}
w_{i,j}=v_{i}\hcirc x_{j}\in\ext(V\hcirc\A^{\rsym}).
\end{gather*}

\begin{definition}\label{D:vertical-horizontal-reading}
Let $S\subset\NBoxm$. We will consider two orderings of $\NBoxm$, \emph{horizontal or row reading} and \emph{vertical or column reading}:
\begin{gather*}
(1,1)\ohor(1,2)\dots \ohor(1,m)\ohor(2,1)\ohor(2,2)\ohor\dots\ohor(N,m-1)\ohor(N,m),
\\
(1,1)\overt(2,1)\dots \overt(N,1)\overt(1,2)\overt(2,2)\overt\dots\overt(N-1,m)\overt(N,m).
\end{gather*}
Suppose that $S=\{(i_{1},j_{1})\ohor\dots\ohor(i_{k},j_{k})\}$
and, at the same time, 
$S=\{(i_{1}^{\prime},j_{1}^{\prime})\overt\dots\overt(i_{k}^{\prime},j_{k}^{\prime})\}$. Define
\begin{gather}\label{Eq:RepsTensorBasis}
w^{h}_{S}=w_{i_{1},j_{1}}\dots w_{i_{k}, j_{k}}
,\quad
w^{v}_{S}=w_{i_{1}^{\prime},j_{1}^{\prime}}\dots
w_{i_{k}^{\prime},j_{k}^{\prime}},
\end{gather}
as the given elements of $\ext(V\hcirc\A^{\rsym})$.
\end{definition}

Essentially by definition we get:

\begin{lemma}\label{L:RepsABasisExtThree}
The sets $\{w^{h}_{S}|S\subset\NBoxm\}$,
$\{w^{v}_{S}|S\subset\NBoxm\}$ are 
$\A$-bases of $\ext(V\hcirc\A^{\rsym})$.\qed
\end{lemma} 

There are projections $\pi_{h}\colon\NBoxm\to[1,\rsym]$ and $\pi_{v}\colon\NBoxm\to[1,\lsym]$. Write $S_{i}=\pi_{h}^{-1}(\{i\})$ to denote the subset of $S$ which projects to $i$ under $\pi_{h}$ and $_{j}S=\pi_{v}^{-1}(\{j\})$ to denote the subset of $S$ which projects to $j$ under $\pi_{v}$. Note that $\coprod_{i=1}^{\lsym}{{}_{i}S} =S=\coprod_{j=1}^{\rsym}S_{j}$. 

\begin{example}\label{E:RepsBoxCombinatorics}
In dot diagram notation, now with an $\lsym$-by-$\rsym$ 
rectangle, the elements 
in \cref{Eq:RepsTensorBasis} are simply given by reading along columns or rows, e.g.:
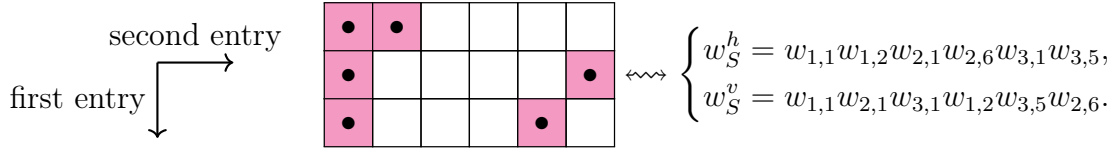
\begin{gather*}
\begin{tikzpicture}[anchorbase]
\draw[thick,->] (0,0) to (0,-0.5)node[left]{first entry} to (0,-1);
\draw[thick,->] (0,0) to (0.5,0)node[above]{second entry} to (1,0);
\end{tikzpicture}
\quad
\begin{ytableau}
*(magenta!50)\bullet & *(magenta!50)\bullet & \phantom{a} 
& \phantom{a} & \phantom{a} & \phantom{a}
\\
*(magenta!50)\bullet & \phantom{a} & \phantom{a} 
& \phantom{a} & \phantom{a} & *(magenta!50)\bullet
\\
*(magenta!50)\bullet & \phantom{a} & \phantom{a} 
& \phantom{a} & *(magenta!50)\bullet & \phantom{a}
\end{ytableau}
\leftrightsquigarrow
\left\{
\begin{gathered}
w^{h}_{S}=w_{1,1}w_{1,2}w_{2,1}w_{2,6}w_{3,1}w_{3,5},
\\
w^{v}_{S}=w_{1,1}w_{2,1}w_{3,1}w_{1,2}w_{3,5}w_{2,6}.
\end{gathered}
\right.
\end{gather*}
Here $(\lsym,\rsym)=(3,6)$. Moreover, we have
\begin{gather*}
S_{1}=\{(1,1),(2,1),(3,1)\}
,\quad
S_{2}=\{(1,2)\}
,\quad
S_{5}=\{(3,5)\}
,\quad
S_{6}=\{(2,6)\}
,
\\
{}_{1}S=\{(1,1),(1,2)\}
,\quad
{}_{2}S=\{(2,1),(2,6)\}
,\quad
{}_{3}S=\{(3,1),(3,5)\}
,
\end{gather*}
and all other of these sets are empty. In other words, the sets
${}_{i}S$ and $S_{j}$ are row and column reading projections.
\end{example}

The following lemma is easy and omitted:

\begin{lemma}\label{L:RepsABasisExtTwo}
The set $\{v_{T}|T\subset[1,\lsym]\}$ is 
an $\A$-basis of $\ext(V)$.\qed
\end{lemma} 

\begin{lemma}\label{L:RepsFixIso}
Horizontal and vertical reading give isomorphisms 
of free $\A$-modules by
\begin{align*}
\hreading&\colon
\ext(V\hcirc\A^{\rsym})\to\ext(\A^{\rsym})^{\hcirc\lsym}
,\quad 
w^{h}_{S}\mapsto 
x_{\pi_{h}({}_{1}S)}\hcirc\dots\hcirc x_{\pi_{h}({}_{\lsym}S)},
\\
\vreading&\colon
\ext(V\hcirc\A^{\rsym})\to\ext(V)^{\hcirc m}
,\quad 
w^{v}_{S}\mapsto 
v_{\pi_{v}(S_{1})}\hcirc\dots\hcirc v_{\pi_{v}(S_{\rsym})}.
\end{align*}
\end{lemma}

\begin{proof}
We combine \cref{L:RepsABasisExt} and \cref{L:RepsABasisExtTwo}.
\end{proof}

In identifications we fix the isomorphisms in \cref{L:RepsFixIso}.

\begin{example}
Let $N=4$, $m=4$, and $S=\{(2,1),(3,1),(1,2),(3,2),(4,2),(2,4)\}$. Then:
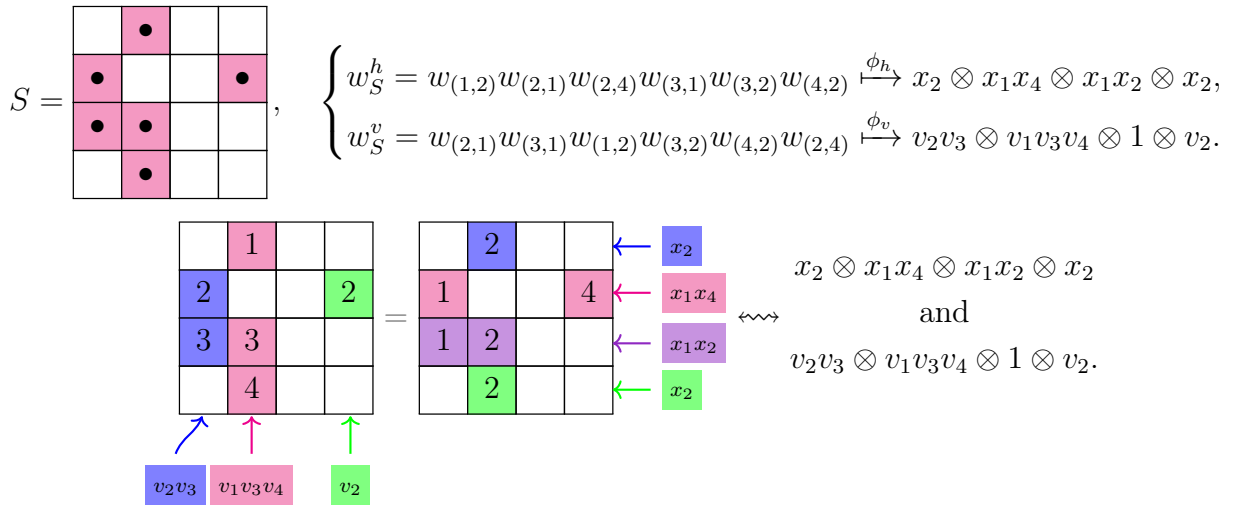
\begin{gather*}
S=
\begin{ytableau}
\phantom{a} & *(magenta!50)\bullet & \phantom{a} 
& \phantom{a}
\\
*(magenta!50)\bullet & \phantom{a} & \phantom{a} 
& *(magenta!50)\bullet
\\
*(magenta!50)\bullet & *(magenta!50)\bullet & \phantom{a} 
& \phantom{a}
\\
\phantom{a} & *(magenta!50)\bullet & \phantom{a} 
& \phantom{a}
\end{ytableau}
\,,\quad
\left\{
\begin{gathered}
w_{S}^{h}=w_{(1,2)}w_{(2,1)}w_{(2,4)}w_{(3,1)}w_{(3,2)}w_{(4,2)} \xmapsto{\hreading}x_{2}\hcirc x_{1}x_{4}\hcirc x_{1}x_{2}\hcirc x_{2},
\\
w_{S}^{v}=w_{(2,1)}w_{(3,1)}w_{(1,2)}w_{(3,2)}w_{(4,2)}w_{(2,4)}
\xmapsto{\vreading} v_{2}v_{3}\hcirc v_{1}v_{3}v_{4}\hcirc 1\hcirc v_{2}.
\end{gathered}
\right.
\\
\begin{tikzpicture}[anchorbase]
\node at (0,0) {
\begin{ytableau}
\phantom{a} & *(magenta!50)1 & \phantom{a} 
& \phantom{a}
\\
*(blue!50)2 & \phantom{a} & \phantom{a} 
& *(green!50)2
\\
*(blue!50)3 & *(magenta!50)3 & \phantom{a} 
& \phantom{a}
\\
\phantom{a} & *(magenta!50)4 & \phantom{a} 
& \phantom{a}
\end{ytableau}
=
\begin{ytableau}
\phantom{a} & *(blue!50)2 & \phantom{a} 
& \phantom{a}
\\
*(magenta!50)1 & \phantom{a} & \phantom{a} 
& *(magenta!50)4
\\
*(orchid!50)1 & *(orchid!50)2 & \phantom{a} 
& \phantom{a}
\\
\phantom{a} & *(green!50)2 & \phantom{a} 
& \phantom{a}
\end{ytableau}};
\draw[blue,thick,->] (-2.9,-1.8)node[below,black]{\colorbox{blue!50}{\tiny\mystrut$v_{2}v_{3}$}} to[out=70,in=250] (-2.55,-1.3);
\draw[magenta,thick,->] (-1.9,-1.8)node[below,black]{\colorbox{magenta!50}{\tiny\mystrut$v_{1}v_{3}v_{4}$}} to (-1.9,-1.3);
\draw[green,thick,->] (-0.6,-1.8)node[below,black]{\colorbox{green!50}{\tiny\mystrut$v_{2}$}} to (-0.6,-1.3);
\draw[blue,thick,->] (3.35,0.95)node[right,black]{\colorbox{blue!50}{\tiny\mystrut$x_{2}$}} to (2.85,0.95);
\draw[magenta,thick,->] (3.35,0.335)node[right,black]{\colorbox{magenta!50}{\tiny\mystrut$x_{1}x_{4}$}} to (2.85,0.335);
\draw[orchid,thick,->] (3.35,-0.335)node[right,black]{\colorbox{orchid!50}{\tiny\mystrut$x_{1}x_{2}$}} to (2.85,-0.335);
\draw[green,thick,->] (3.35,-0.95)node[right,black]{\colorbox{green!50}{\tiny\mystrut$x_{2}$}} to (2.85,-0.95);
\end{tikzpicture}
\raisebox{0.6cm}{$
\leftrightsquigarrow
\begin{gathered}
x_{2}\hcirc x_{1}x_{4}\hcirc x_{1}x_{2}\hcirc x_{2}
\\
\text{and}
\\
v_{2}v_{3}\hcirc v_{1}v_{3}v_{4}\hcirc 1\hcirc v_{2}.
\end{gathered}$}
\end{gather*}
We have also illustrated how to get the expression for $w_{S}^{v}$ and $w_{S}^{h}$ on the bottom and right, respectively.
Note that $w_{S}^{h}=-w_{S}^{v}$.
\end{example}

Let $\sigma_{h}^{v}(S)$ be the permutation 
that sends the ordered set $\{(i_{1},j_{1}),\dots,(i_{k},j_{k})\}$ 
to the ordered set $\{(i^{\prime}_{1},j^{\prime}_{1}),
\dots,(i^{\prime}_{k},j^{\prime}_{k})\}$. Let $\ell$ denote its length.

\begin{lemma}\label{L:RepsSignReading}
We have $w_{S}^{h}=(-1)^{\ell(\sigma_{h}^{v}(S))}\cdot w_{S}^{v}$.
\end{lemma}

\begin{proof}
Directly from the signed commutation rules of the exterior algebra.
\end{proof}

%%%%%%%%%%%%%%%%%%%%%%%%%

\subsection{Two explicit actions on $\ext(V\hcirc\A^{\rsym})$}

%%%%%%%%%%%%%%%%%%%%%%%%%

Each of the two bases $\{w_{S}^{v}\}$ and $\{w_{S}^{h}\}$ for $\ext(V\hcirc\A^{\rsym})$ are adapted to the action by $\ogroup$ and $\uplain(\gso)$, respectively. We describe these actions now. The two actions commute, see \cref{P:actions-commute}. 

%%%%%%%%%%%%%%%%%%%%%%%%%

\subsubsection{The $\ogroup$-action}

%%%%%%%%%%%%%%%%%%%%%%%%%

The group $\ogroup$ acts naturally on 
$\ext(V)^{\hcirc\rsym}$ and on the exterior 
algebra $\ext(V\hcirc\A^{\rsym})$ of 
$V\hcirc\A^{\rsym}\cong V\oplus\dots\oplus V$. 
Fixing this $\ogroup$-action we have:

\begin{lemma}\label{L:RepsOActs}
The map $\vreading$ intertwines the $\ogroup$-action on $\ext(V\hcirc\A^{\rsym})$.
\end{lemma}

\begin{proof}
A direct check.
\end{proof}

\begin{remark}\label{R:SquareRootMinusOne}
There are two common choices of symmetric bilinear form when defining $\ogroup$. These have \emph{Gram matrices} either the diagonal or the antidiagonal. Since $\sqrt{-1}\in\A$ one can prove that these two forms are equivalent.     
\end{remark}

We define another basis for $V$, such that $(\placeholder,\placeholder)$ will have antidiagonal Gram 
matrix with respect to the new basis. 

\begin{definition}\label{D:abu-basis}
For $i\in\{1,\dots,n\}$, define
\begin{gather*}
a_{i}=v_{i}-\sqrt{-1}\cdot v_{\lsym-i+1}, 
\quad 
b_{i}=\frac{v_{i}+\sqrt{-1}\cdot v_{\lsym-i+1}}{2}, \quad \text{and}
\\
u=v_{\llsym+1}\text{ if $\lsym=2\llsym+1$}.
\end{gather*}
(Note that we define $u$ only if $\lsym$ is odd.)
\end{definition}

For $i\in\{1,\dots,n\}$ we can express the $v_{i}$ basis in terms of this new basis as
\begin{gather*}
v_{i}=\frac{a_{i}+2\cdot b_{i}}{2}
,\quad 
v_{\lsym-i+1}=\sqrt{-1}\cdot\frac{a_{i}-2\cdot b_{i}}{2}, \quad \text{and}
\\
v_{\llsym+1}=u,\text{ if $\lsym=2\llsym+1$}. 
\end{gather*}

\begin{lemma}\label{L:Pairing}
The pairing $(v_{i},v_{j})=\delta_{ij}$ gives
\begin{gather*}
(a_{i},b_{i})=1=(b_{i},a_{i})=(u,u),
\end{gather*}
while all other pairings of basis vectors $a_{i}$, $b_{j}$ and $u$ vanish. 
\end{lemma}

\begin{proof}
A routine calculation.
\end{proof}

\begin{definition}\label{D:split-torus}
Let the \emph{split torus} $T\subset\sogroup$ be the diagonal matrices in $\sogroup\subset\ogroup$ with respect to the basis in \cref{D:abu-basis}.
\end{definition}

The group $T$ is the subgroup generated by operators $\alpha^{\vee}_{\epsilon_{i}}(t)\in\sogroup$, for $i\in\{1,\dots,n\}$, determined by the action
$\alpha^{\vee}_{\epsilon_{i}}(t)\acts a_{j}=t^{(\epsilon_{i},\epsilon_{j})}\cdot a_{j}$, $\alpha^{\vee}_{\epsilon_{i}}(t)\acts b_{j}=t^{(\epsilon_{i},-\epsilon_{j})}\cdot b_{j}$
and $\alpha^{\vee}_{\epsilon_{i}}(t)\acts u=u$.

\begin{remark}\label{R:v-not-wt-basis}
The action of $T$ in the $v_{i}$ basis is now easily computed, for example:
\begin{gather*}
\alpha^{\vee}_{\epsilon_{1}}(t)\acts v_{1}=
\frac{t+t^{-1}}{2}\cdot v_{1}
-\sqrt{-1}\cdot\frac{t-t^{-1}}{2}\cdot v_{\lsym}. 
\end{gather*}
This action of $T$ in the $v_{i}$ basis is not as easy to work with for our purposes, which is the reason we introduce the new basis. 
\end{remark}

\begin{definition}\label{D:OThingy}
Let $\sigma\in\ogroup$ be the 
following element determined by its action.
\begin{enumerate}[label=\arabic*.]

\item If $\lsym=2\llsym+1$, then $\sigma\in\ogroup$ acts on $V$ by 
\begin{gather*}
\sigma\acts a_{i}=a_{i}
,\quad
\sigma\acts b_{i}=b_{i}
,\quad
\sigma\acts u=-u.
\end{gather*}

\item If $\lsym=2\llsym$, then $\sigma\in\ogroup$ acts on $V$ by 
\begin{gather*}
\sigma \acts a_{i}=a_{i}
,\quad 
\sigma \acts b_{i}=b_{i}
,\quad
\sigma\acts a_{\lsym}=\tfrac{1}{2}\cdot b_{\lsym}
,\quad
\sigma\acts b_{\lsym}=2\cdot a_{\lsym}.
\end{gather*}

\end{enumerate}
Here $i\in\{1,\dots,n-1\}$.
\end{definition}

Note that conjugation by $\sigma$ preserves $\sogroup$ and induces an automorphism of $\lso$, which agrees with the automorphism, also denoted by $\sigma$, 
defined in \cref{D:U(oN)-defn}

%%%%%%%%%%%%%%%%%%%%%%%%%

\subsubsection{The $\uplain(\gso)$-action}

%%%%%%%%%%%%%%%%%%%%%%%%%

\begin{definition}\label{D:Howes-action-via-phi-hor}
Using the isomorphism $\phi_{h}$ and the action of $\uplain(\gso)$ on $\ext(\A^{\rsym})^{\hcirc N}$ from \cref{N:RepsDistActsTensor}, we can define an action of $\uplain(\gso)$ on $\ext(V\hcirc\A^{\rsym})$ by
\begin{gather*}
u\acts w=\phi_{h}^{-1}\big(u\acts\phi_{h}(w)\big),
\end{gather*}
for all $u\in\uplain(\gso)$ and $w\in\ext(V\hcirc\A^{\rsym})$.
\end{definition}

\begin{notation}\label{N:Howes-action-via-phi-hor}
For $S\subset\NBoxm$, write $d_{j}(S)=\sum_{i=1}^{\lsym}d_{j}({}_{i}S)$ 
for $j\in[1,\rsym]$. That is $d_{j}(S)$ is equal to the number of boxes in the $j$th column with a dot minus the number of boxes in the $j$th column without a dot.
\end{notation}

\begin{example}\label{E:Howes-action-via-phi-hor}
For the dot diagram
\begin{gather*}
S=
\begin{ytableau}
\phantom{a} & *(magenta!50)\bullet & \phantom{a} 
& \phantom{a}
\\
*(magenta!50)\bullet & \phantom{a} & \phantom{a} 
& *(magenta!50)\bullet
\\
*(magenta!50)\bullet & *(magenta!50)\bullet & \phantom{a} 
& \phantom{a}
\\
\phantom{a} & *(magenta!50)\bullet & \phantom{a} 
& \phantom{a}
\end{ytableau}
\rightsquigarrow
d_{1}=2-2 = 0,d_{2}=3-1=2,d_{3}=0-4=-4,d_{4}=1-3=-2,
\end{gather*}
where $d_{i}=d_{i}(S)$.
\end{example}

\begin{lemma}\label{L:RepsMoreOActions}
The following defines a
$\uplain(\gso)$-action on $\ext(V\hcirc\A^{\rsym})$:
\begin{gather*}
e_{j}\acts w_{S}^{h}=
\sum_{\substack{1\leq i\leq\lsym\\(i,j+1)\in S\\(i,j)\notin S}} 
w^{h}_{(S\setminus\{(i,j+1)\})\cup\{(i,j)\}}
,\quad
f_{j}\acts w_{S}^{h}=
\sum_{\substack{1\leq i\leq\lsym\\(i,j)\in S\\(i,j+1)\notin S}} w^{h}_{(S\setminus\{(i,j)\})\cup\{(i,j+1)\}}
,
\\
e_{\rsym}\acts w_{S}^{h}=
\sum_{\substack{1\leq i\leq\lsym\\(i,m-1)\notin S\\(i,m)\notin S}}w^{h}_{S\cup\{(i,m-1),(i,m)\}}
,\quad
f_{\rsym}\acts w_{S}^{h}=
\sum_{\substack{1\leq i\leq\lsym\\(i,m-1)\in S\\ (i,m)\in S}}w^{h}_{S\setminus\{(i,m-1),(i,m)\}}
,
\\
h_{k}\acts w_{S}^{h}
=\alpha_{k}^{\vee}\Big(\tfrac{1}{2}\big(d_{1}(S), \ldots,d_{\rsym}(S)\big)\Big)\cdot w_{S}^{h},
\end{gather*}
for $j\in[1,\rsym-1]$ and $k\in[1,\rsym]$.
\end{lemma}

\begin{proof}
Using the definition of $\hreading$, this follows 
from the description of the $\uplain(\gso)$-action on $\ext(\A^{\rsym})^{\hcirc\lsym}$ in 
\cref{N:RepsDistActsTensor}. In fact, $\hreading$ is designed so that this lemma is true.

More explicitly, in dot diagrams the action is row-wise.
For example,
\begin{gather*}
f_{1}\acts
\begin{ytableau}
*(magenta!50)\bullet & *(magenta!50)\bullet & \phantom{a} 
& \phantom{a} & \phantom{a} & \phantom{a}
\\
*(magenta!50)\bullet & \phantom{a} & \phantom{a} 
& \phantom{a} & \phantom{a} & *(magenta!50)\bullet
\\
*(magenta!50)\bullet & \phantom{a} & \phantom{a} 
& \phantom{a} & *(magenta!50)\bullet & \phantom{a}
\end{ytableau}
=
\begin{ytableau}
*(magenta!50)\bullet & *(magenta!50)\bullet & \phantom{a} 
& \phantom{a} & \phantom{a} & \phantom{a}
\\
\phantom{a} & *(magenta!50)\bullet &\phantom{a} 
& \phantom{a} & \phantom{a} & *(magenta!50)\bullet
\\
*(magenta!50)\bullet & \phantom{a} & \phantom{a} 
& \phantom{a} & *(magenta!50)\bullet & \phantom{a}
\end{ytableau}
+
\begin{ytableau}
*(magenta!50)\bullet & *(magenta!50)\bullet & \phantom{a} 
& \phantom{a} & \phantom{a} & \phantom{a}
\\
*(magenta!50)\bullet & \phantom{a} & \phantom{a} 
& \phantom{a} & \phantom{a} & *(magenta!50)\bullet
\\
\phantom{a} & *(magenta!50)\bullet & \phantom{a} 
& \phantom{a} & *(magenta!50)\bullet & \phantom{a}
\end{ytableau}
,
\end{gather*}
which is the rule for rows as exemplified in 
\cref{Eq:RepsDotActionF} plus the coproduct.
\end{proof}

%%%%%%%%%%%%%%%%%%%%%%%%%

\subsection{Semisimple Howe duality}

%%%%%%%%%%%%%%%%%%%%%%%%%

For the following \emph{semisimple Howe duality} we recall that 
$\ParNm$ denotes the set of $\rsym$-restricted dominant $\ogroup$-weights, cf. \cref{S:BackgroundWeights}, which 
index the simple $\ogroup$-representations (as well as Weyl and 
indecomposable tilting representations). Moreover, below we will 
give a combinatorial map $\dagger$, again defined later in \cref{S:BackgroundWeights}, that takes a dominant $\ogroup$-weight and produces a dominant $\gso$-weight.
The following is our version of \cite[Section 4.3.5]{Ho-perspectives-invariant-theory}:

\begin{proposition}(Semisimple Howe duality.)\label{P:bimod-decomp-over-Q}
As an $\ogroup$-$\uplain[\K](\gso)^{op}$-birepresentation we have:
\begin{gather*}
\K\hcirc_{\A}\ext(V\hcirc\A^{\rsym})\cong 
\bigoplus_{\lambda\in\ParNm}
\ksimple[\lambda]\boxtimes\ksimple[\lambda^{\dagger}].
\end{gather*}
\end{proposition}

\begin{proof}
The key ingredient is the result \cref{P:actions-commute}, which makes 
the question about a $\ogroup$-$\uplain[\K](\gso)^{op}$-birepresentation decomposition of $\ext(V\hcirc\A^{\rsym})$ well-defined, 
and character calculations in the spirit of Howe's original construction. 

In more details, let $\lambda\in\ParNm$. By \cref{L:z-lambda-is-joint-high-low-wt-vector}, \cref{R:z-generates-L(lambda)-L(w0(wtb))}, and \cref{D:lambda-dagger}, it follows from complete reducibility that $\bigoplus_{\lambda\in\ParNm}\ksimple[\lambda]\boxtimes\ksimple[\lambda^{\dagger}]$ is isomorphic to a direct summand $\K\hcirc_{\A}\ext(V\hcirc\A^{\rsym})$. The claim then follows from the character calculation in \cite[Section 4.3.5]{Ho-perspectives-invariant-theory} or
\cite[Proposition 3.2]{AdRy-tilting-howe-positive-char}. 
\end{proof}

Let $\DCompN$ denote the set of $\lsym$-restricted dominant $\gso$-weights, specified in \cref{S:BackgroundWeights}.

\begin{proposition}\label{P:dim-End-over-O}
We have
\begin{gather*}
\dim_{\F}\End_{O(V)}\big(\F\hcirc_{\A}\ext(V\hcirc\A^{\rsym})\big)
\!=\!\sum_{\lambda\in\ParNm}\big(\dim_{\K}\ksimple[\lambda^{\dagger}]\big)^{2}
\!=\!\sum_{\wtk\in\DCompN}\big(\dim_{\K}\ksimple[\wtk]\big)^{2}.
\end{gather*}
\end{proposition}

\begin{proof}
Using that $\F\hcirc_{\A}\ext(V\hcirc\A^{\rsym})$ is a tilting representations for both actions, as follows from \cref{L:has-Weyl-filt-for-som}, and \cref{L:O-delta/nabla-Ext-vanishing}, a standard argument, similar to \cite[Proposition 2.3]{AnStTu-semisimple-tilting}, yields
\begin{gather*}
\dim_{\F}\End_{O(V)}\big(\F\hcirc_{\A}\ext(V\hcirc\A^{\rsym})\big)=
\dim_{\K}\End_{O(V)}\big(\K\hcirc_{\A}\ext(V\hcirc\A^{\rsym})\big)
\end{gather*}
Therefore, \cref{P:bimod-decomp-over-Q} implies
\begin{gather*}
\dim_{\F}\End_{O(V)}
\big(\F\hcirc_{\A}\ext(V\hcirc\A^{\rsym})\big)
=\sum_{\lambda\in\ParNm}\big(\dim_{\K}\ksimple[\lambda^{\dagger}]\big)^{2},
\end{gather*}
and the final claim then follows from \cref{P:order-reversing-bijection} proven later on.
\end{proof}

\cref{P:dim-End-over-O} gives an effective way to compute 
the dimension of the endomorphism space over $\F$ since $\dim_{\K}\ksimple[\lambda^{\dagger}]$ 
can be computed using \emph{Weyl's dimension formula}.

%%%%%%%%%%%%%%%%

\section{The diagrammatic presentation}\label{S:fullyfaithfulfunctor}

%%%%%%%%%%%%%%%%

We now prove our first main theorem: the equivalence of symmetric ribbon $\F$-linear categories between
$\webF$ and $\FundF$, see \cref{T:mainthm-integral-webs}.
Upon additive idempotent completion, $\webF$ is thus a \emph{diagrammatic 
version} of $\TiltF$ (recall tilting representations from 
\cref{SS:RepsTiltingO}).
Our main tool is a 
version of Howe's $\ogroup$-$\gso$ duality in positive characteristic,  cf. \cref{SS:FunctorIntegralHowe}.

%%%%%%%%%%%%%%%%

\subsection{From webs to reps}\label{SS:FunctorDefiningMaps}

%%%%%%%%%%%%%%%%

We begin by defining some $\ogroup$-equivariant maps which 
correspond to our generating webs.

\begin{notation}\label{Nota:Lmabda-tuple}
We let $\ext[i]$ denote $\ext[i](V)$ and write
$\ext[(i_{1},\dots,i_{\rsym})]=\ext[i_{1}](V)\hcirc
\dots\hcirc\ext[i_{\rsym}](V)$. 
\end{notation}

Let $T=\{t_{1}<\dots<t_{i}\} $ and $U=\{u_{i+1}<\dots<u_{i+j}\}$ be subsets of $[1,\lsym]$ such that 
$U\cap T=\emptyset$ and 
write $T\cup U=\{s_{1}<\dots< s_{i+j}\}$. Consider the 
permutation $\sigma_{T,U}\in\sgroup[i+j]$ in the symmetric group $\sgroup[i+j]$ of $[1,i+j]$ which 
sends $\{t_{1}<\dots,t_{i},u_{i+1}<\dots<u_{i+j}\}$ 
to $\{s_{1}<\dots<s_{i+j}\}$. We let 
$\ell(T,U)=\ell(\sigma_{T,U})$ (here the length $\ell$ is in terms of number of simple transpositions). 

\begin{definition}\label{D:FunctorMaps}
We define \emph{merge, split and crossing maps} to be
\begin{align*}
\Pa_{k,l}^{k+l}&\colon\ext[k]\hcirc\ext[l]\to\ext[{k+l}],
v_{T}\hcirc v_{U}\mapsto(-1)^{\ell(T,U)}v_{T\cup U},
\\
\Ya_{k+l}^{k,l}&\colon\ext[{k+l}]\to\ext[k]\hcirc\ext[l],
v_{S}\mapsto
\sum_{\substack{S=T\coprod U \\|T|=k,|U|=l}}
(-1)^{\ell(T,U)}v_{T}\hcirc v_{U},
\\
\XX_{k,l}^{l,k}&\colon\ext[k]\hcirc\ext[l]\to\ext[l]\hcirc\ext[k],
v_{T}\hcirc v_{S}\mapsto v_{S}\hcirc v_{T},
\end{align*}
and \emph{cap and cup maps}
\begin{align*}
\ca_{k}&\colon\ext[k]\hcirc\ext[k]\to\wtlnit,
v_{T}\hcirc v_{U}\mapsto(-1)^{\ell(T,U)}v_{T\cup U},
\\
\cu_{k}&\colon\wtlnit\to\ext[k]\hcirc\ext[k],
1\mapsto (-1)^{\binom{i}{2}}\sum_{\substack{S\subset [1,\lsym] \\ |S|=k}} v_{S}\hcirc v_{S}.
\end{align*}
(The notation is hopefully suggestive.)
\end{definition}

Recall that $\ogroup$ acts on the spaces in \cref{D:FunctorMaps}.

\begin{lemma}\label{L:FunctorMaps}
The maps in \cref{D:FunctorMaps} are 
$\ogroup$-equivariant.
\end{lemma}

\begin{proof}
A direct calculation.
\end{proof}

If $\wtk\in \Pi_{\rsym}^{\le \lsym}$, then $\wtk=(k_1, \dots, k_m)$ where $k_{i}\in \{0, 1, \dots, N\}$. We write $\Lambda^{\wtk}:=\Lambda^{k_1}\otimes \dots \otimes \Lambda^{k_m}$ as in \cref{Nota:Lmabda-tuple}.

\begin{definition}
Define $\FundF$ as the category with objects $\ext[\wtk]$ for all $\rsym\in\N$ and all 
$\wtk\in\Pi_{\rsym}^{\le\lsym}$, and morphisms all $\F$-linear maps which commute with $\ogroup$.
\end{definition}

In other words, $\FundF$ is the category of \emph{fundamental $\ogroup$-representations}, i.e. the full subcategory 
of all finite dimensional $\ogroup$-representations monoidally generated by the Weyl representations for the fundamental $\ogroup$-weights.

\begin{lemma}
The category $\FundF$ is a symmetric ribbon 
$\F$-linear category with:
\begin{enumerate}[label=(\roman*)]

\item with monoidal structure given by the usual tensor product of $\ogroup$-representations; 

\item with symmetry given by the tensor flip;

\item with pivotal structure given by
$X^{\pivo}=\Hom_{\F}(X,\F)$ and 
$X\rightarrow X^{\pivo\pivo}$ defined 
by $x\mapsto(f\mapsto f(x))$.

\end{enumerate}
\end{lemma}

\begin{proof}
$\FundF$ is a full subcategory 
of all finite dimensional $\ogroup$-representations, and inherits all structures from the parent category.
\end{proof}

\begin{proposition}(The (diagrammatic) presentation functor.)\label{P:def-of-webs-to-fund}
There is a symmetric ribbon $\F$-linear functor 
\begin{gather*}
\pfunctor\colon\webF\to\FundF,
\\
k\mapsto\ext[k],
\\
\begin{tikzpicture}[anchorbase,scale=1]
\draw[usual] (0,0)node[below]{$k$} to (0.5,0.5);
\draw[usual] (1,0)node[below]{$l$} to (0.5,0.5);
\draw[usual] (0.5,0.5) to (0.5,1)node[above]{$k{+}l$};
\end{tikzpicture}
\mapsto\text{\upshape$\Pa_{k,l}^{k+l}$}
,\;
\begin{tikzpicture}[anchorbase,scale=1,yscale=-1]
\draw[usual] (0.5,0.5) to (0,0)node[above]{$k$};
\draw[usual] (0.5,0.5) to (1,0)node[above]{$l$};
\draw[usual] (0.5,1)node[below]{$k{+}l$} to (0.5,0.5);
\end{tikzpicture}
\mapsto\Ya^{k,l}_{k+l}
,\;
\begin{tikzpicture}[anchorbase,scale=1]
\draw[usual] (0,0)node[below]{$k$} to (1,1)node[above]{$k$};
\draw[usual] (1,0)node[below]{$l$} to (0,1)node[above]{$l$};
\end{tikzpicture}
\mapsto\XX_{k,l}^{l,k}
,\;
\begin{tikzpicture}[anchorbase,scale=1]
\draw[usual] (0,0)node[below]{$k$} to[out=90,in=180] (0.5,0.5)node[above,yshift=0.4cm]{$\wtlnit$} to[out=0,in=90] (1,0)node[below]{$k$};
\end{tikzpicture}
\mapsto\text{\upshape$\ca_{k}$}
,\;
\begin{tikzpicture}[anchorbase,scale=1]
\draw[usual] (0,0)node[above]{$k$} to[out=270,in=180] (0.5,-0.5)node[below,yshift=-0.4cm]{$\wtlnit$} to[out=0,in=270] (1,0)node[above]{$k$};
\end{tikzpicture}
\mapsto\cu_{k}
.
\end{gather*}
\end{proposition}

\begin{proof}
If $\pfunctor$ is well-defined, 
the other claims follow easily.
Hence, it suffices to check that the defining relations of $\webF$, see \cref{D:DiagramsWeb}, are satisfied in $\FundF$. That the exterior type A web relations hold follows from \cite[Proposition 5.2.1]{CaKaMo-webs-skew-howe} via \cref{P:DiagramsFromAToO}. Circle removal relations hold since $\ext[k]$ has categorical dimension $\binom{N}{k}$ and lollipop relations hold by Schur's lemma and the fact that $\ext[k]$ are simple and pairwise nonisomorphic, cf. \cref{L:SimpleTilting}.

It remains to prove that the higher even orthogonal $E$-$F$ relations hold in $\FundF$. As explained in \cref{L:SimpleTilting}, the objects in $\FundF$ are tensor products of the objects
\[
\Lambda^i\F^N\cong \F\otimes_{\A}\Lambda^i \cong \F\otimes\Delta_{\A}(1^i)\cong \F\otimes_{\A}\nabla_{\A}(1^i), \text{for} \ i=1, \dots, N.
\]
By \cref{L:o-tensor-weyl-is-weyl} and \cref{L:O-delta/nabla-Ext-vanishing} it follows that the dimension of homomorphism spaces in $\FundF$ are independent of $\F$. Thus, the corresponding homomorphism spaces are free, finitely generated, and torsion free over $\A$, and it follows that we can check relations over $\K$. Working over $\K$, we can use $X^{(a)} = a!X^a$ to replace divided powers with Chevalley generators, so it suffices to prove the orthogonal $E$-$F$ relation in \cref{ef-1label}. But, as explained in \cref{ex-orthoog-e-f} this relation is known to be a consequence of the other orthogonal web relations and therefore holds in $\catstuff{Fund}_{\K}(\ogroup)$.
\end{proof}

\begin{definition}\label{D:e/f-web}
We define \emph{ladder operators} as follows. First:
\begin{align*}
E_{i,j}^{i+1,j-1}&=(\Pa_{i,1}^{i+1}\hcirc
\id_{\ext[{j-1}]})\circ(\id_{\ext[i]}
\hcirc\Ya_{j}^{1,j-1})
\colon
\ext[i]\hcirc\ext[j]\to\ext[{i+1}]\hcirc\ext[{j-1}],
\\
F_{i,j}^{i-1,j+1}&=(\id_{\ext[{i-1}]}
\hcirc\Pa_{1,j}^{j+1})\circ
(\Ya_{i}^{i-1, 1}\hcirc\id_{\ext[j]})\colon
\ext[i]\hcirc\ext[j]\rightarrow
\ext[{i-1}]\hcirc\ext[{j+1}],
\\
e_{i,i}^{i+1,i+1}&=(\Pa_{i,1}^{i+1}\hcirc
\Pa_{1,i}^{i+1})\circ(\id_{\ext[i]}\hcirc
\cu_{1}\hcirc\id_{\ext[i]})\colon
\ext[i]\hcirc\ext[i]\to\ext[{i+1}]\hcirc\ext[{i+1}],
\\
f_{i,i}^{i-1,i-1}&=(\id_{\ext[{i-1}]}\hcirc\ca_{1}
\hcirc\id_{\ext[{i-1}]})\circ(\Ya_{i}^{i-1, 1}
\hcirc\Ya_{i}^{1,i-1})\colon
\ext[i]\hcirc\ext[i]\to\ext[{i-1}]\hcirc\ext[{i-1}]. 
\end{align*}
Moreover, let $\wtk\in\Pi_{\rsym}$. We define operators on 
$\ext(V\hcirc\A^{\rsym})$ by
\begin{align*}
\textbf{E}_j\one_{\wtk}\acts w_{S}^{v}&=\id_{\Lambda^{(k_{1}, \dots,k_{j-1})}}\hcirc E_{k_{j},k_{j+1}}^{k_{j}+1,k_{j+1}-1}\hcirc \id_{\Lambda^{(k_{j+2},\dots,k_{m})}}\big(\phi_v(w_{S}^{v})\big),
\\
\textbf{F}_j\one_{\wtk}\acts w_{S}^{v}&=\id_{\Lambda^{(k_{1}, \dots,k_{j-1})}}\hcirc F_{k_{j},k_{j+1}}^{k_{j}-1,k_{j+1}+1}\hcirc \id_{\Lambda^{(k_{j+2},\dots,k_{m})}}\big(\phi_v(w_{S}^{v})\big),
\\
\textbf{e}_{\rsym}\one_{\wtk}\acts w_{S}^{v}&=\id_{\Lambda^{(k_{1}, \dots,k_{m-2})}}\hcirc e_{k_{m-1},k_{m}}^{k_{m-1}+1,k_{m}+1} \big(\phi_v(w_{S}^{v})\big),
\\
\textbf{f}_{\rsym}\one_{\wtk}\acts w_{S}^{v}&=\id_{\Lambda^{(k_{1}, \dots,k_{m-2})}}\hcirc e_{k_{m-1},k_{m}}^{k_{m-1}-1,k_{m}-1} \big(\phi_v(w_{S}^{v})\big),
\end{align*}
if $\wt w_{S}^{v}=\wtk$, and $w_{S}^{v}\mapsto 0$ otherwise.
(The notation is again hopefully suggestive.)
\end{definition}

\begin{lemma}\label{L:explicit-idempotent-Howe-operators}
Let $\wtk\in\Pi_{\rsym}^{\le\lsym}$ and let $\wt x_{S}^{h}=\wtk$. 
We have the following explicit description of 
the action of the ladder operators:
\begin{align*}
\textbf{E}_j\one_{\wtk}\acts w_{S}^{v}&= \sum_{\substack{1\le i  \le N\\ (i,j)\notin S, (i,j+1)\in S}} (-1)^{\ell(S_{j}, \{i\}) + \ell(\{i\}, S_{j+1}\setminus \{i+1\})} w^{v}_{S\cup \{(i,j)\}\setminus \{(i,j+1)\}},
\\
\textbf{F}_j\one_{\wtk}\acts w_{S}^{v}&=\sum_{\substack{1\le i  \le N\\ (i,j)\in S, (i,j+1)\notin S}} (-1)^{\ell(S_{j}\setminus\{i\}, \{i\}) + \ell(\{i\}, S_{j+1})} w^{v}_{S\setminus \{(i,j)\}\cup \{(i,j+1)\}},
\\
\textbf{e}_{i}\one_{\wtk}\acts w_{S}^{v}&=\sum_{\substack{1\le i\le N \\ (i,j)\notin S, (i,j+1)\notin S}} (-1)^{\ell(S_{j}, \{i\}) + \ell(\{i\}, S_{j+1})} w^{v}_{S\cup\{(i,j), (i,j+1)\}},
\\
\textbf{f}_{i}\one_{\wtk}\acts w_{S}^{v}&= \sum_{\substack{1\le i\le N \\ (i, j)\in S, (i, j+1)\in S}} (-1)^{\ell(S_{j}\setminus\{i\}, \{i\}) + \ell(\{i\}, S_{j+1}\setminus\{i\})} w^{v}_{S\setminus\{(i,j), (i, j+1)\}}.
\end{align*}
\end{lemma}

\begin{proof}
This follows from the Definition of $\vreading$ and the formulas
\begin{align*}
E_{i,j}^{i+1,j-1}(v_{S}\hcirc v_{T})&=
\sum_{\substack{t\in T \\ t\notin S}}(-1)^{\ell(S,\{t\})
+
\ell(\{t\},T\setminus\{t\})}v_{S\cup\{t\}}\hcirc 
v_{T\setminus\{t\}},
\\
F_{i,j}^{i-1,j+1}(v_{S}\hcirc v_{T})&= 
\sum_{\substack{s\in S \\ s\notin T}}
(-1)^{\ell(S\setminus\{s\},\{s\})+\ell(\{s\},T)}
v_{S\setminus\{s\}}\hcirc v_{T\cup\{s\}},
\\
e_{i,i}^{i+1,i+1}(v_{S}\hcirc v_{T})&=
\sum_{\substack{x\notin S\\ x\notin T}}
(-1)^{\ell(S,\{x\})+\ell(\{x\},T)}
v_{S\cup\{x\}}\hcirc v_{T\cup\{x\}},
\\
f_{i,i}^{i-1,i-1}(v_{S}\hcirc v_{S})&= 
\sum_{\substack{y\in S\\ y\in T}}(-1)^{\ell(S\setminus\{y\},\{y\})
+\ell(\{y\},T\setminus\{y\})}v_{S\setminus\{y\}}\hcirc
v_{T\setminus\{y\}}.
\end{align*}
That these formulas hold is a direct calculation.
\end{proof}

\begin{proposition}(Howe's actions diagrammatically.)\label{P:two-actions-agree}
The ladder operators above define an $\uplain(\gso)$-action that
agrees with the $\uplain(\gso)$-action in \cref{D:Howes-action-via-phi-hor}.
\end{proposition}

\begin{proof}
Using \cref{L:RepsMoreOActions} and \cref{L:explicit-idempotent-Howe-operators}, it is easy to see the actions agree up to a sign. Since the action in \cref{D:Howes-action-via-phi-hor} has no signs, it suffices to show that the signs in \cref{L:RepsMoreOActions} cancel the signs coming from \cref{L:RepsSignReading}. 

For example, to verify the two $e_{j}$ actions agree, we have to check that
\begin{gather*}
(-1)^{\ell(\sigma_{h}^{v}(S))}
=(-1)^{\ell(\sigma_{h}^{v}(S\cup\{(i,j)\}\setminus\{(i,j+1)\}))}  (-1)^{\ell(S_{j},\{i\})+
\ell(\{i\}, S_{j+1}\setminus\{i+1\})}.
\end{gather*}
We leave the verification of the remaining cases to the reader.
\end{proof}

We finally arrive at a proof, using webs, that the actions of $\ogroup$ and $\uplain(\gso)$ on the space $\ext(V\hcirc\A^{\rsym})$ commute.

\begin{proof}[Proof of \cref{P:actions-commute}]
From \cref{P:two-actions-agree}, we see the action of the generators of the algebra $\uplain(\gso)$ on $\ext(V\hcirc\A^{\rsym})$ is given in \cref{L:explicit-idempotent-Howe-operators} which are compositions of tensor products of the $\ogroup$ intertwiners in \cref{D:FunctorMaps}, see \cref{D:e/f-web}
\end{proof}

%%%%%%%%%%%%%%%%

\subsection{Howe duality}\label{SS:FunctorIntegralHowe}

%%%%%%%%%%%%%%%%

It follows from \cref{L:check-saturated}, see also \cref{R:check-saturated}, that the set of $\gso$-weights appearing in $\ext(V\hcirc\A^{\rsym})$, denoted $\CompN$ above, is a 
\emph{saturated set}, see \cref{D:BackgroundSaturated}. Recall also
that $\DCompN$ is a saturated set of dominant $\gso$-weights, see \cref{D:BackgroundSaturatedTwo}.

\begin{notation}
Write $\splain(\gso)$ to denote the \emph{Schur algebra} quotient of $\udot(\gso)$ determined by the saturated set $\DCompN$.
The reader unfamiliar with these is referred to the background in \cref{SS:SchurAlgebra} below.
\end{notation}

The algebra $\splain(\gso)$, is the quotient of $\udot(\gso)$ by the two sided ideal generated by $\one_{\wtk}$ such that $\wtk\notin \CompN$. Moreover, this algebra contains orthogonal idempotents $\one_{\wtk}$, for $\wtk\in\CompN$, such that 
$\sum_{\wtk\in\CompN}\one_{\wtk}=1$.

One can think of $\splain(\gso)$ as an algebraic version of webs with $m$ boundary points on the top and bottom of the diagram.

\begin{notation}
Similar to \cref{N:DiagramsIdempotentCategory}, for fixed $\rsym$, we view $\splain(\gso)$ as a category with objects $\one_{\wtk}$, for $\wtk\in\CompN$, and morphisms $\Hom_{\splain(\gso)}(\one_{\wtk}, \one_{\wtl})=\one_{\wtl}\splain(\gso)\one_{\wtk}$. The composition of morphisms is multiplication in the algebra.  
\end{notation}

\begin{lemma}(A Schur functor.)
The action of $\uplain(\gso)$ on $\ext(V\hcirc\A^{\rsym})$ induces a functor
\begin{gather*}
\sfunctor[\A]\colon\splain(\gso)\to\Fund
\end{gather*}
such that $\one_{\wtk}\mapsto\ext[\wtk]$, 
for all $\wtk\in\CompN$.
\end{lemma}

\begin{proof}
\cref{L:wtk-space-is-Lambda-k} implies that 
\begin{gather*}
\one_{\wtk}\acts\ext(V\hcirc\A^{\rsym})= 
\begin{cases}
\ext[\wtk] & \text{if $\wtk\in\CompN$},
\\
0 & \text{if $\wtk\notin\CompN$}.
\end{cases}
\end{gather*}
Therefore, the action of $\uplain(\gso)$ induces a functor $\splain(\gso)\to\Fund$.
\end{proof}

We have already discussed how 
$\ext(V\hcirc\A^{\rsym})\cong\ext(\A^{\rsym})^{\hcirc\lsym}$. We also have: 

\begin{lemma}\label{L:full-tilting-som}
The space $\ext(V\hcirc\F^{\rsym})$ is a full tilting representation for $\splain[\F](\gso)$. 
\end{lemma}

\begin{proof}
Because of \cref{L:has-Weyl-filt-for-som}, it suffices to show that if $\wtk\in\DCompN$, then the indecomposable 
tilting $\uplain[\F](\gso)$-representation $\ftilting[\wtk]$ (see \cref{S:BackgroundOrthogonal}) is a summand of $\ext(V\hcirc\F^{\rsym})$, which is \cref{L:T(wtk)-is-summand}. 
\end{proof}

\begin{proposition}[One-sided double commutant]\label{P:O-so-Howe-duality}
The functor $\sfunctor$ is fully faithful.
\end{proposition}

\begin{proof}
Using \cref{L:full-tilting-som}, \cref{L:AR-lemma}, and \cref{P:actions-commute} we find that $\splain[\F](\gso)$ injects into $\End_{\ogroup}\big(\ext(V\hcirc\F^{\rsym})\big)$. It then follows from \cref{E:schuralgebras-weyl-character}, \cref{P:dim-End-over-O}, and the equality $\rank_{\A}\weyl[\wtk]=
\dim_{\K}\ksimple[\wtk]=\dim_{\F}\fweyl[\wtk]$, that 
\begin{gather*}
\dim_{\F}\splain[\F](\gso)=
\sum_{\wtk\in\DCompN}\big(\rank_{\A}\weyl[\wtk]\big)^{2}
=\dim_{\F}\End_{O(V)}\big(\ext(V\hcirc\F^{\rsym})\big).
\end{gather*}
Injectivity then implies that $\splain[\F](\gso)\rightarrow \End_{O(V)}\big(\ext(V\hcirc\F^{\rsym})\big)$ is surjective.
\end{proof}

The other side of the double commutant theorem, \cref{P:O-so-Howe-duality} for the orthogonal action, 
can be done similarly. However, we do not need it in this work.

%%%%%%%%%%%%%%%%

\subsection{fully faithful}\label{SS:FullyFaithful}

%%%%%%%%%%%%%%%%

The full functor $\dhowe$ in \cref{P:DiagramsHowe} factors through the Schur quotient, inducing a full functor
\begin{gather*}
\ddhowe\colon
\splain(\gso)\twoheadrightarrow\web,
\one_{\wtk}\mapsto\wtk_{1}\hcirc\dots\hcirc\wtk_{\rsym},
\end{gather*}
and both \cref{Eq:DiagramsCKMFunctor} and \cref{Eq:DiagramsSTFunctor} hold.
Thus, we have the commuting diagram
\begin{gather*}
\begin{tikzcd}[ampersand replacement=\&,column sep=1.5cm,row sep=1cm]
\udot(\gso)\ar[r,thick,->>,"\dhowe"]
\ar[dr,thick,->>,"\text{quo.}"'] \& \web
\\
\& \splain(\gso)\ar[u,thick,->>,"\ddhowe"']
\end{tikzcd}
.
\end{gather*}

\begin{proposition}\label{P:alpha-equals-epsilon-varphi}
We have the following commuting diagram:
\begin{gather*}
\begin{tikzcd}[ampersand replacement=\&,column sep=1.5cm,row sep=1cm]
\udot(\gso)\ar[r,thick,->>,"\dhowe"]
\ar[dr,thick,->>,"\text{quo.}"'] \& \web\ar[r,thick,->,"\pfunctor"] \& \Fund
\\
\& \splain(\gso)\ar[u,thick,->>,"\ddhowe"']
\ar[ur,line width=0.025cm,hook,two heads,"\sfunctor"']\&
\end{tikzcd}
.
\end{gather*}
\begin{gather}\label{Eq:action-equals-web}
\sfunctor=\pfunctor\circ\ddhowe.
\end{gather}
\end{proposition}
\begin{proof}
Note that both sides of \cref{Eq:action-equals-web} send $\one_{\wtk}$ to $\ext[\wtk]$. The desired claim then follows from \cref{D:e/f-web}, \cref{P:two-actions-agree}, and \cref{P:def-of-webs-to-fund}.
\end{proof}

Denote by $\Kar[\placeholder]$ the additive idempotent completion.
For the following theorem recall the category of tilting $\ogroup$-representations $\TiltF$ as in \cref{SS:RepsTiltingO}.

\begin{theorem}(The orthogonal web calculus.)\label{T:mainthm-integral-webs}
The presentation functor 
\begin{gather*}
\pfunctor\colon\webF\rightarrow\FundF
\end{gather*}
is an equivalence of symmetric ribbon $\F$-linear categories.
Moreover,
\begin{gather*}
\Kar[\webF]\cong\TiltF
\end{gather*}
as symmetric ribbon $\F$-linear categories.
\end{theorem}

\begin{proof}
That the functor preserves the structures as in the statement follows immediately 
from the definitions. The functor is also essentially surjective by 
the construction of $\FundF$, so it remains to argue that $\pfunctor$ is fully faithful.

Suppose $W$ is a morphism in $\webF$ which $\pfunctor$ maps to zero. Since $\ddhowe$ is full by \cref{P:DiagramsHowe}, there is a morphism $u$ in $\dot{S}_{\F}(\gso)^{\leq\lsym}$ so that 
$W=\ddhowe(u)$. Then \cref{Eq:action-equals-web} implies $\sfunctor(u)=\pfunctor(W)=0$. \cref{P:O-so-Howe-duality} says that $u=0$, so $W=\ddhowe(u)=0$. Thus, $\pfunctor$ is faithful. A similar easy argument using \cref{P:O-so-Howe-duality} and \cref{Eq:action-equals-web} shows that $\pfunctor$ is full.

The final claim follows then from the equivalence and \cref{P:Kar-Fund-is-Tilt}.
\end{proof}

We actually proved our main theorem:

\begin{proof}[Proof of \cref{T:Main}]
\cref{T:mainthm-integral-webs} is a more 
refined formulation of \cref{T:Main}.
\end{proof}

\begin{proof}[Proofs of \cref{P:DiagramsFromAToO} and \cref{P:DiagramsHowe}, wrap-up]
We use the result \cref{T:mainthm-integral-webs} to fill in some of the remaining statements from \cref{S:Diagrams}.

\textit{The case of \cref{P:DiagramsFromAToO}.} We note that we have the commuting diagram (we marked the functor we are interested in)
\begin{gather*}
\begin{tikzcd}[ampersand replacement=\&,column sep=1.5cm,row sep=1cm]
\webb[\F]\ar[r,thick,->,"\ppfunctor","\cong"'] 
\ar[d,thick,->,dotted,"\ifunctor"']
\& \FundFF\ar[d,thick,->,hook]
\\
\web[\F]\ar[r,thick,->,"\pfunctor","\cong"'] \& \FundF
\end{tikzcd}
.
\end{gather*}
Here $\FundFF$ is the analog of $\FundF$ but for the general linear group, and $\ppfunctor$ is the presentation functor 
from \cite{CaKaMo-webs-skew-howe}. We note:
\begin{enumerate}[label=$\bullet$]

\item Commutativity of this diagram follows by careful comparison of 
the definitions.

\item The top and bottom functors are equivalences by 
\cite[Theorem 3.3.1]{CaKaMo-webs-skew-howe} and 
\cref{T:mainthm-integral-webs}, respectively.

\item The right functor is 
faithful since $\ogroup[\lsym]$ is a subgroup of $\glgroup$.

\end{enumerate}
These together 
implies that the left functor is faithful.

\textit{The case of \cref{P:DiagramsHowe}.}
Directly from \cref{P:O-so-Howe-duality}
and \cref{T:mainthm-integral-webs}.
\end{proof}

%%%%%%%%%%%%%%%%

\section{Semisimplification}\label{S:Semisimple}

%%%%%%%%%%%%%%%%

We first rapidly recall the notion of \emph{semisimplification} of a rigid symmetric monoidal category such that the endomorphisms of the unit object is spanned by the identity. Then we prove \cref{T:O-ss} by constructing an equivalence, in analogy with \cite[Section 3]{BrEnAiEtOs-semisimple-tilting}.
The basic material below can be found in many sources such 
as \cite{EtGeNiOs-tensor-categories} or \cite{MR4486913}.

%%%%%%%%%%%%%%%%

\subsection{A reminder on semisimplifications}\label{SS:SemisimpleReminder}

%%%%%%%%%%%%%%%%

All categories are assumed to be $\F$-linear and have finite dimensional homomorphism spaces. Let $\catstuff{C}$ be a rigid symmetric monoidal category with monoidal unit $\wtlnit$, such that $\End_{\catstuff{C}}(\wtlnit)=\F\cdot\id_{\wtlnit}$. 

\begin{notation}\label{N:SemisimpleEvaluation}
Write $\ev_{\obstuff{X}}\colon\obstuff{X}^{\pivo}\hcirc\obstuff{X}\rightarrow\wtlnit$, $\coev\colon\wtlnit\rightarrow\obstuff{X}\hcirc\obstuff{X}^{\pivo}$ for the unit and counit of adjunction realizing $\obstuff{X}^{\pivo}$ as the \emph{right dual} of $\obstuff{X}$. Let $P_{\obstuff{X},\obstuff{Y}}\colon\obstuff{X}\hcirc\obstuff{Y}\rightarrow\obstuff{Y}\hcirc\obstuff{X}$ denote the symmetry.  
\end{notation}

\begin{definition}\label{D:SemisimpleTrace}
The \emph{trace} of an endomorphism $f\colon\obstuff{X}
\rightarrow\obstuff{X}$ in $\catstuff{C}$, denoted $\tr_{\catstuff{C}}(\morstuff{f})$, is the element of $\F$ such that
\begin{gather*}
\ev_{\obstuff{X}}\circ P_{\obstuff{X},\obstuff{X}^{\pivo}}\circ 
(\morstuff{f}\hcirc\id)\circ\coev_{\obstuff{X}}=\tr_{\catstuff{C}}(\morstuff{f})
\cdot\id_{\wtlnit}.
\end{gather*}
We define the \emph{categorical dimension} of $\obstuff{X}$ as $\dim_{\catstuff{C}}(\obstuff{X})=\tr_{\catstuff{C}}(\id_{\obstuff{X}})$. 
\end{definition}

\begin{remark}
The category of vector spaces over $\F$, denoted $\catstuff{Vec}_{\F}$, is a rigid symmetric monoi\-dal category such that the endomorphisms of the monoidal identity are spanned by the identity endomorphism. The categorical dimension of a vector space is the usual dimension of the vector space under the map $\mathbb{Z}\rightarrow \F$. More generally, if there is an $\F$-linear symmetric monoidal functor $F:\catstuff{C}\rightarrow \catstuff{Vec}_{\F}$, then $\dim_{\catstuff{C}}(X)$ is equal to the dimension of the vector space $F(X)$ when one views $\dim_{\F}F(X)$ as an element of $\F$. 
\end{remark}

We make the following simplifying assumption: every nilpotent endomorphism in $\catstuff{C}$ has trace zero. 

\begin{lemma}\label{L:SemiAssumption}
If there is a symmetric monoidal functor from $\catstuff{C}$ to an abelian category, then $\catstuff{C}$ satisfies the simplifying assumption.
\end{lemma}

\begin{proof}
See \cite[Remark 2.9]{MR4486913}.
\end{proof}

\begin{remark}\label{R:SemisimpleNeg}
All the categories we consider have the property in \cref{L:SemiAssumption}, and therefore satisfy the simplifying assumption.
\end{remark}

\begin{definition}\label{D:SemisimpleNeg}
The subcategory of \emph{negligible morphisms} in $\catstuff{C}$, denoted $\mathcal{N}_{\catstuff{C}}$, is the category with the same objects as $\catstuff{C}$, and with morphisms
\begin{gather*}
\Hom_{\mathcal{N}_{\catstuff{C}}}(\obstuff{X},\obstuff{Y})
= 
\{\morstuff{f}\colon\obstuff{X}\rightarrow\obstuff{Y}|\text{such that $\tr_{\catstuff{C}}(\morstuff{f}\circ\morstuff{g})=0$, for all $\morstuff{g}\colon\obstuff{Y}\rightarrow\obstuff{X}$}\}
.
\end{gather*}
The \emph{semisimplification} of $\catstuff{C}$ is the quotient category 
\begin{gather*}
\semisimple{\catstuff{C}}=\catstuff{C}/\mathcal{N}_{\catstuff{C}},
\end{gather*}
which is defined as the category with the same objects as $\catstuff{C}$, and with morphisms
\begin{gather*}
\Hom_{\semisimple{\catstuff{C}}}(\obstuff{X},\obstuff{Y})
=\Hom_{\catstuff{C}}(\obstuff{X},\obstuff{Y})/\Hom_{\mathcal{N}_{\catstuff{C}}}(\obstuff{X},\obstuff{Y}).
\end{gather*}
Write $\pi_{\catstuff{C}}\colon\catstuff{C}\rightarrow \semisimple{\catstuff{C}}$ for the quotient functor.
\end{definition}

\begin{lemma}
The category $\semisimple{\catstuff{C}}$ is monoidal.
\end{lemma}

\begin{proof}
It is well-known that $\mathcal{N}_{\catstuff{C}}$ is a $\vcirc$-$\hcirc$-ideal, see e.g. \cite[Theorem 2.6]{MR4486913}, which implies the claim.
\end{proof}

\begin{example}\label{E:SemisimpleNeg}
Key examples of semisimplifications are 
the \emph{Verlinde categories} in the spirit of \cite{AnPa-fusion-lie-algebras}, see also \cite[Section 8.18.2]{EtGeNiOs-tensor-categories} 
for background and references. These are constructed as 
follows. Let $G$ be a simple algebraic group over $\Z$ with Coxeter number $h$. let $\F$ be an infinite field of
characteristic $p\geq h$. Then the Verlinde category $\Ver_{\F}(G)$
for $G$ is the semisimplification of the
category of tilting $G$-representations:
\begin{gather*}
\Ver_{\F}(G):=\semisimple{\Tilt_{\F}(G)}.
\end{gather*}
This construction works more generally, e.g. also 
for $G=\glgroup$ or $G=\ogroup[\lsym]$. The latter category
$\Vero$ plays an important role in this paper.
\end{example}

For the purpose of this paper, a category is \emph{semisimple} if it is abelian and every object is isomorphic to a finite direct sum of simple objects.

\begin{lemma}\label{L:SemisimpleNeg}
We have the following.
\begin{enumerate}[label=\arabic*.]

\item The category $\semisimple{\catstuff{C}}$ is rigid symmetric monoidal and $\End_{\semisimple{\catstuff{C}}}(\one_{\semisimple{\catstuff{C}}})=\F\cdot\id_{\one_{\semisimple{\catstuff{C}}}}$. 

\item The category $\semisimple{\catstuff{C}}$ is semisimple.

\item An object $\semisimple{\obstuff{X}}$ is a simple objects in $\semisimple{\catstuff{C}}$ if and only if $\obstuff{X}$ is an indecomposable object in $\catstuff{C}$ with $\dim_{\catstuff{C}}(\obstuff{X})\neq 0$. Two such simple objects are isomorphic in $\semisimple{\catstuff{C}}$ if and only if the corresponding indecomposable objects are isomorphic in $\catstuff{C}$. 

\item Suppose $\catstuff{D}$ is a rigid symmetric monoidal semisimple category. If $\functorstuff{F}\colon\catstuff{C}\rightarrow\catstuff{D}$ is a full symmetric monoidal functor, then there is a fully faithful symmetric monoidal functor $\semisimple{\functorstuff{F}}\colon\semisimple{\catstuff{C}}\rightarrow\catstuff{D}$ such that $\functorstuff{F}=\semisimple{\functorstuff{F}}\circ\pi_{\catstuff{C}}$, {i.e.} the following diagram commutes:
\begin{gather*}
\begin{tikzcd}[ampersand replacement=\&,column sep=1em,row sep=1em]
\catstuff{C}\ar[rr,thick,->,"\functorstuff{F}"]
\ar[dr,thick,->,"\pi_{\catstuff{C}}"'] \& \& \catstuff{D}
\\
\& \semisimple{\catstuff{C}}\ar[ur,thick,dotted, ->,"\semisimple{\functorstuff{F}}"'] \&
\end{tikzcd}
.
\end{gather*}

\end{enumerate}
\end{lemma}

\begin{proof}
Because of the simplifying assumption, all of this is in 
\cite[Section 2]{MR4486913}.
\end{proof}

%%%%%%%%%%%%%%%%

\subsection{Colored Brauer diagrams}\label{SS:ColorBrauer}

%%%%%%%%%%%%%%%%

By a \emph{two-part partition diagram} with $b$ bottom and $t$ top points we mean a diagram of the following form:
\begin{gather}\label{Eq:DiagramsBrauer}
\begin{tikzpicture}[anchorbase]
\draw[usual] (0.5,0) to[out=90,in=180] (1.25,0.45) to[out=0,in=90] (2,0);
\draw[usual] (1,0) to[out=90,in=180] (1.25,0.25) to[out=0,in=90] (1.5,0);
\draw[usual] (0,1) to[out=270,in=180] (0.75,0.55) to[out=0,in=270] (1.5,1);
\draw[usual] (1,1) to[out=270,in=180] (1.75,0.55) to[out=0,in=270] (2.5,1);
\draw[usual] (0,0) to (0.5,1);
\draw[usual] (2.5,0) to (2,1);
\draw[usual] (3,1) to[out=270,in=180] (3.25,0.55) to[out=0,in=270] (3.5,1);
\end{tikzpicture}
\quad\text{where }b=6,t=8.
\end{gather}
In words, a two-part partition diagram with $b$ bottom and $t$ top points is 
a diagram corresponding to a partition of $\{1,\dots,b\}\cup\{1,\dots,t\}$
where every block has two parts. A \emph{Brauer diagram}
with $\rsym$ bottom and $\lsym$ top points is then any representative 
of diagrams that represent the same partition.
These assemble into the \emph{Brauer category} $\brauertwo_d$ \cite{Br-brauer-algebra-original, LeZh-brauer-invariant-theory}.

This category can be thought of as the symmetric ribbon $\F$-linear category $\hcirc$-generated by a selfdual object $\bullet$ of categorical dimension 
$d$. In other words, the category $\brauertwo_d$ has a universal mapping property \cite[Theorem 2.6]{LeZh-brauer-invariant-theory}.

\begin{proposition}\label{P:DiagramsBrauer}
For $k=1, \dots, N$, there is a symmetric ribbon $\F$-linear functor
\begin{gather*}
\brauertwo_{\binom{N}{k}}\to\web[\F],
\\
\bullet\mapsto k,
\begin{tikzpicture}[anchorbase,scale=1]
\draw[usual] (0,0)node[below]{\phantom{$k$}} to (1,1)node[above]{\phantom{$k$}};
\draw[usual] (1,0)node[below]{\phantom{$k$}} to (0,1)node[above]{\phantom{$k$}};
\end{tikzpicture}
\mapsto
\begin{tikzpicture}[anchorbase,scale=1]
\draw[usual] (0,0)node[below]{$k$} to (1,1)node[above]{$k$};
\draw[usual] (1,0)node[below]{$k$} to (0,1)node[above]{$k$};
\end{tikzpicture}
,
\begin{tikzpicture}[anchorbase,scale=1]
\draw[usual] (0,0)node[below]{\phantom{$k$}} to[out=90,in=180] (0.5,0.5) to[out=0,in=90] (1,0)node[below]{\phantom{$k$}};
\end{tikzpicture}
\mapsto
\begin{tikzpicture}[anchorbase,scale=1]
\draw[usual] (0,0)node[below]{$k$} to[out=90,in=180] (0.5,0.5) to[out=0,in=90] (1,0)node[below]{$k$};
\end{tikzpicture}
,\quad
\begin{tikzpicture}[anchorbase,scale=1]
\draw[usual] (0,0)node[above]{\phantom{$k$}} to[out=270,in=180] (0.5,-0.5) to[out=0,in=270] (1,0)node[above]{\phantom{$k$}};
\end{tikzpicture}
\mapsto
\begin{tikzpicture}[anchorbase,scale=1]
\draw[usual] (0,0)node[above]{$k$} to[out=270,in=180] (0.5,-0.5) to[out=0,in=270] (1,0)node[above]{$k$};
\end{tikzpicture}
,
\end{gather*}
where we draw the defining structures of $\brauertwo_{\binom{N}{k}}$ 
in the standard way.
\end{proposition}

\begin{proof}
Using the Brauer categories universal mapping property, this follows from \cref{L:DiagramsWebRibbon}.
\end{proof}

For $r\in\Z_{\geq 1}$, a \emph{colored Brauer diagram} with colors $\{0,\dots,r-1\}$ is a diagram of the form
\begin{gather*}
\begin{tikzpicture}[anchorbase]
\draw[usual1] (0.5,0) to[out=90,in=180] (1.25,0.45) to[out=0,in=90] (2,0);
\draw[usual1] (1,0) to[out=90,in=180] (1.25,0.25) to[out=0,in=90] (1.5,0);
\draw[usual2] (0,1) to[out=270,in=180] (0.75,0.55) to[out=0,in=270] (1.5,1);
\draw[usual0] (1,1) to[out=270,in=180] (1.75,0.55) to[out=0,in=270] (2.5,1);
\draw[usual0] (0,0) to (0.5,1);
\draw[usual0] (2.5,0) to (2,1);
\draw[usual1] (3,1) to[out=270,in=180] (3.25,0.55) to[out=0,in=270] (3.5,1);
\end{tikzpicture}
\quad\text{where }r=3\text{ and }
\left(\,
\begin{tikzpicture}[anchorbase]
\draw[usual0] (0,0)to (0,1);
\end{tikzpicture}
\,,
\begin{tikzpicture}[anchorbase]
\draw[usual1] (0,0)to (0,1);
\end{tikzpicture}
\,,
\begin{tikzpicture}[anchorbase]
\draw[usual2] (0,0)to (0,1);
\end{tikzpicture}
\,\right)
\leftrightsquigarrow
(0,1,2).
\end{gather*}

\begin{definition}\label{D:ColorBrauer}
Let $(d_{i})=(d_{0},\dots,d_{r-1})\in\Z^{r}$. Define 
the \emph{colored Brauer category} $\brauertwo_{(d_{i})}$ to be the rigid symmetric monoidal category with objects tensor products of the self-dual objects of dimension $d_{i}$ for $i=0, \dots, r-1$, and morphisms $\{0,\dots,r-1\}$-colored Brauer diagrams. 
\end{definition}

Recall that $\Kar[\placeholder]$ denotes the additive idempotent completion.

\begin{lemma}\label{L:Abelian}
The category $\Kar[\brauertwo]_{(d_{i})}$ is an additive idempotent closed symmetric ribbon $\F$-linear category such that the endomorphisms of $\one$ are the $\F$-span of $\id_{\one}$. 
Moreover, $\Kar[\brauertwo]_{(d_{i})}$ admits an 
abelian envelope which is also a symmetric 
ribbon $\F$-linear category.
\end{lemma}

\begin{proof}
The additive idempotent completion of a symmetric 
ribbon $\F$-linear category is a symmetric 
ribbon $\F$-linear category. To see that the endomorphisms of the unit object are all scalar multiples of the identity, note that the only Brauer diagram with empty bottom and top is the empty diagram.
For the final claim we refer to \cite[Theorem A]{MR4277111}.
\end{proof}

Note that \cref{L:SemisimpleNeg} implies that all the semisimplifications that we will see below are semisimple. Recall the notion of \emph{Deligne tensor product} \cite[Section 1.11]{EtGeNiOs-tensor-categories}, which we denote by 
$\boxtimes$. The Deligne tensor product preserves the class of $\F$-linear semisimple rigid symmetric monoidal categories \cite[Section 4.6]{EtGeNiOs-tensor-categories}.

\begin{lemma}\label{L:ColorBrauerEquivalence}
There is an equivalence of symmetric ribbon $\F$-linear categories
\begin{gather*}
\semisimple{\Kar[\brauertwo_{(d_{i})}]}\to
\boxtimes_{i=0}^{r-1}\semisimple{\Kar[\brauertwo_{d_{i}}]}.
\end{gather*}
\end{lemma}

\begin{proof}
We follow \cite[Proof of Lemma 3.3]{BrEnAiEtOs-semisimple-tilting}. The universal mapping property inherent in the definition of $\Kar[\brauertwo_{(d_{i})}]$ gives rise to a functor $\Kar[\brauertwo_{(d_{i})}]\to\boxtimes_{i=0}^{r-1}\semisimple{\Kar[\brauertwo_{d_{i}}]}$. This functor is clearly full, so we can apply \cref{L:SemisimpleNeg} to deduce there is an equivalence $\semisimple{\Kar[\brauertwo_{(d_{i})}]}\to\boxtimes_{i=0}^{r-1}\semisimple{\Kar[\brauertwo_{d_{i}}]}$. 
\end{proof}

After replacing $\Kar[\brauertwo_{(d_{i})}]$ by its abelian envelope 
from \cref{L:Abelian}, one could hope that \cref{L:ColorBrauerEquivalence} 
holds without the semisimplification.

%%%%%%%%%%%%%%%%

\subsection{Bases at infinity}\label{SS:DiagramsAtInfinity}

%%%%%%%%%%%%%%%%

To prove that after semisimplification, the image of the colored Brauer diagrams in the webs $\web[\F]$ generate the category, we will use a spanning set for homomorphism spaces in 
the orthogonal web category which is analogous to the chicken foot diagrams in \cite[Lemma 4.9]{BrEnAiEtOs-semisimple-tilting}.

For completeness, we define two spanning sets for the morphisms in $\web$. 
One of these is built from \emph{many-to-few-to-many (mfm)}
diagrams and the other from \emph{few-to-many-to-few (fmf)} diagrams, with many and few referring to the total number of strings.

\begin{remark}\label{R:DiagramsAtInfinity}
We expect that the spanning sets we are defining are actually bases whenever the total thickness 
of the strands is $\ll N$, but we will not need or prove this.
\end{remark}

A web is a mfm bottom part if its associated partition contains no splits and caps 
and has a minimal number of crossings, its a 
mfm top part if its associated partition contains no merges and cups 
and has a minimal number of crossings.
Dually, we define fmf bottom part and
and a fmf top part by swapping the roles of splits and merges.
Finally, a web is a sandwiched part if it contains only crossings.

\begin{definition}\label{D:DiagramsAtInfinity}
A web is called a \emph{mfm sandwich diagram} if it is of the form
\begin{gather*}
\begin{tikzpicture}[anchorbase,scale=1]
\draw[mor] (0,-0.5) to (0.25,0) to (0.75,0) to (1,-0.5) to (0,-0.5);
\node at (0.5,-0.25){$\rsym$};
\draw[mor] (0,1) to (0.25,0.5) to (0.75,0.5) to (1,1) to (0,1);
\node at (0.5,0.75){$m^{\prime}$};
\draw[mor] (0.25,0) to (0.25,0.5) to (0.75,0.5) to (0.75,0) to (0.25,0);
\node at (0.5,0.25){$f$};
\draw[thick,->] (1.5,0.75)node[right]{mfm top part - splits and cups} to (1,0.75);
\draw[thick,->] (1.5,0.25)node[right]{sandwiched part - crossings} to (1,0.25);
\draw[thick,->] (1.5,-0.25)node[right]{mfm bottom part - merges and caps} to (1,-0.25);
\end{tikzpicture}
.
\end{gather*}
The set of all many-to-few sandwich diagrams from $\obstuff{K}$ to $\obstuff{L}$ by $\sand{K}{L}$.

Similarly, a web is called a \emph{fmf sandwich diagram} if it is of the form
\begin{gather*}
\begin{tikzpicture}[anchorbase,scale=1]
\draw[mor] (0,0) to (0.25,-0.5) to (0.75,-0.5) to (1,0) to (0,0);
\node at (0.5,-0.25){$f$};
\draw[mor] (0,0.5) to (0.25,1) to (0.75,1) to (1,0.5) to (0,0.5);
\node at (0.5,0.75){$f^{\prime}$};
\draw[mor] (0,0) to (0,0.5) to (1,0.5) to (1,0) to (0,0);
\node at (0.5,0.25){$\rsym$};
\draw[thick,->] (1.5,0.75)node[right]{fmf top part - merges and cups} to (1,0.75);
\draw[thick,->] (1.5,0.25)node[right]{sandwiched part - crossings} to (1,0.25);
\draw[thick,->] (1.5,-0.25)node[right]{fmf bottom part - splits and caps} to (1,-0.25);
\end{tikzpicture}
.
\end{gather*}
For the through strands, say the thicknesses of strands add up to $\rsym$ (this number is the same at every generic horizontal cut). Then we also require the crossings in the middle have to be shortest coset representatives of
types $(b_{1},\dots,b_{k})$ and $(t_{1},\dots,t_{l})$ in $\sgroup[\rsym]$, where the $b$ and $t$ are the bottom and top endpoints of the through strands. The set of all few-to-many sandwich diagrams from $\obstuff{K}$ to $\obstuff{L}$ is denoted by $\sandtwo{K}{L}$.
\end{definition}

In \cref{D:DiagramsAtInfinity} we sandwich a symmetric group in between merges, splits, caps and caps.

\begin{example}\label{E:DiagramsAtInfinity}
The diagram
\begin{gather*}
\begin{tikzpicture}[anchorbase,scale=1]
\draw[usual] (0,0)node[below]{$1$} to (0.5,0.5);
\draw[usual] (1,0)node[below]{$1$} to (0.5,0.5);
\draw[usual] (0.5,0.5) to (0.5,1) to (0.5,1.5);
\draw[usual] (-0.5,0)node[below]{$1$} to[out=90,in=180] (0.5,1) to[out=0,in=90] (1.5,0)node[below]{$1$};
\draw[usual] (-1,0)node[below]{$2$} to (-1,1.5);
\draw[usual] (0.5,1.5) to (-1,2);
\draw[usual] (-1,1.5) to (0.5,2);
\draw[usual] (-0.5,2.5)node[above]{$1$} to (-1,2) to (-1.5,2.5)node[above]{$1$};
\draw[usual] (0,2.5)node[above]{$1$} to (0.5,2) to (1,2.5)node[above]{$1$};
\end{tikzpicture}
\in\sand{\obstuff{(2,1,1,1,1)}}{\obstuff{(1,1,1,1)}}
,
\\
\begin{tikzpicture}[anchorbase,scale=1]
\draw[mor] (0,-0.5) to (0.25,0) to (0.75,0) to (1,-0.5) to (0,-0.5);
\node at (0.5,-0.25){$\rsym$};
\end{tikzpicture}
=
\begin{tikzpicture}[anchorbase,scale=1]
\draw[usual] (0,0)node[below]{$1$} to (0.5,0.5);
\draw[usual] (1,0)node[below]{$1$} to (0.5,0.5);
\draw[usual] (0.5,0.5) to (0.5,1) to (0.5,1.5)node[above]{$2$};
\draw[usual] (-0.5,0)node[below]{$1$} to[out=90,in=180] (0.5,1) to[out=0,in=90] (1.5,0)node[below]{$1$};
\draw[usual] (-1,0)node[below]{$2$} to (-1,1.5)node[above]{$2$};
\end{tikzpicture}
,\quad
\begin{tikzpicture}[anchorbase,scale=1]
\draw[mor] (0,1) to (0.25,0.5) to (0.75,0.5) to (1,1) to (0,1);
\node at (0.5,0.75){$m^{\prime}$};
\end{tikzpicture}
=
\begin{tikzpicture}[anchorbase,scale=1,yscale=-1]
\draw[usual] (0.5,0.5) to (0,0)node[above]{$1$};
\draw[usual] (0.5,0.5) to (1,0)node[above]{$1$};
\draw[usual] (0.5,1)node[below]{$2$} to (0.5,0.5);
\draw[usual] (2,0.5) to (1.5,0)node[above]{$1$};
\draw[usual] (2,0.5) to (2.5,0)node[above]{$1$};
\draw[usual] (2,1)node[below]{$2$} to (2,0.5);
\end{tikzpicture}
,\quad
\begin{tikzpicture}[anchorbase,scale=1]
\draw[mor] (0.25,0) to (0.25,0.5) to (0.75,0.5) to (0.75,0) to (0.25,0);
\node at (0.5,0.25){$f$};
\end{tikzpicture}
=
\begin{tikzpicture}[anchorbase,scale=1]
\draw[usual] (0,0)node[below]{$2$} to (1,1)node[above]{$2$};
\draw[usual] (1,0)node[below]{$2$} to (0,1)node[above]{$2$};
\end{tikzpicture}
,
\end{gather*}
is an example of a many-to-few sandwich diagram that we also split into its defining pieces.

Moreover, the diagram
\begin{gather*}
\begin{tikzpicture}[anchorbase,scale=1]
\draw[usual] (-1.5,0) to (-1.5,1.5)node[above]{$1$};
\draw[usual] (-0.5,0) to (-0.5,1.5)node[above]{$1$};
\draw[usual] (0,0) to (0.5,0.5);
\draw[usual] (1,0) to (0.5,0.5);
\draw[usual] (0.5,0.5) to (0.5,1) to (0.5,1.5)node[above]{$2$};
\draw[usual] (0,1.5)node[above]{$2$} to[out=270,in=180] (0.5,1) to[out=0,in=270] (1,1.5) node[above]{$2$};
\draw[usual] (-1.5,-0.5) to (-1.5,0);
\draw[usual] (-0.5,-0.5) to (0,0);
\draw[usual] (0,-0.5) to (-0.5,0);
\draw[usual] (1,-0.5) to (1,0);
\draw[usual] (-1,-1) to (-1.5,-0.5);
\draw[usual] (-1,-1) to (-0.5,-0.5);
\draw[usual] (-1,-1.5)node[below]{$2$} to (-1,-1);
\draw[usual] (0.5,-1) to (0,-0.5);
\draw[usual] (0.5,-1) to (1,-0.5);
\draw[usual] (0.5,-1.5)node[below]{$2$} to (0.5,-1);
\end{tikzpicture}
\in\sandtwo{\obstuff{(2,2)}}{\obstuff{(1,1,2,2,2)}},
\\
\begin{tikzpicture}[anchorbase,scale=1]
\draw[mor] (0,1) to (0.25,0.5) to (0.75,0.5) to (1,1) to (0,1);
\node at (0.5,0.75){$f$};
\end{tikzpicture}
=
\begin{tikzpicture}[anchorbase,scale=1,yscale=-1]
\draw[usual] (0.5,0.5) to (0,0)node[above]{$1$};
\draw[usual] (0.5,0.5) to (1,0)node[above]{$1$};
\draw[usual] (0.5,1)node[below]{$2$} to (0.5,0.5);
\draw[usual] (2,0.5) to (1.5,0)node[above]{$1$};
\draw[usual] (2,0.5) to (2.5,0)node[above]{$1$};
\draw[usual] (2,1)node[below]{$2$} to (2,0.5);
\end{tikzpicture}
,\quad
\begin{tikzpicture}[anchorbase,scale=1]
\draw[mor] (0,-0.5) to (0.25,0) to (0.75,0) to (1,-0.5) to (0,-0.5);
\node at (0.5,-0.25){$f^{\prime}$};
\end{tikzpicture}
=
\begin{tikzpicture}[anchorbase,scale=1]
\draw[usual] (-1.5,0)node[below]{$1$} to (-1.5,1.5)node[above]{$1$};
\draw[usual] (-0.5,0)node[below]{$1$} to (-0.5,1.5)node[above]{$1$};
\draw[usual] (0,0)node[below]{$1$} to (0.5,0.5);
\draw[usual] (1,0)node[below]{$1$} to (0.5,0.5);
\draw[usual] (0.5,0.5) to (0.5,1) to (0.5,1.5)node[above]{$2$};
\draw[usual] (0,1.5)node[above]{$2$} to[out=270,in=180] (0.5,1) to[out=0,in=270] (1,1.5) node[above]{$2$};
\end{tikzpicture}
,\quad
\begin{tikzpicture}[anchorbase,scale=1]
\draw[mor] (0.25,0) to (0.25,0.5) to (0.75,0.5) to (0.75,0) to (0.25,0);
\node at (0.5,0.25){$\rsym$};
\end{tikzpicture}
=
\begin{tikzpicture}[anchorbase,scale=1]
\draw[usual] (-1,0)node[below]{$1$} to (-1,1)node[above]{$1$};
\draw[usual] (0,0)node[below]{$1$} to (1,1)node[above]{$1$};
\draw[usual] (1,0)node[below]{$1$} to (0,1)node[above]{$1$};
\draw[usual] (2,0)node[below]{$1$} to (2,1)node[above]{$1$};
\end{tikzpicture}
,
\end{gather*}
is an example of an element in 
$\sandtwo{\obstuff{(2,2)}}{\obstuff{(1,1,2,2,2)}}$.
\end{example}

\begin{proposition}\label{P:DiagramsAtInfinity}
The sets $\sand{K}{L}$ and $\sandtwo{K}{L}$are $\A$-linear spanning sets of $\Hom_{\Web}(\obstuff{K},\obstuff{L})$.
\end{proposition}

\begin{proof}
By \cref{L:DiagramsWebRibbon}, $\web$, we can push 
all trivalent vertices and Morse points (cups and caps) to wherever we want them to be.
The relations
\begin{gather*}
\begin{tikzpicture}[anchorbase,scale=1]
\draw[usual] (0,0)node[below]{$k$} to (1,1);
\draw[usual] (1,0)node[below]{$l$} to (0.5,0.5);
\draw[usual] (2,0)node[below]{$\rsym$} to (1,1);
\draw[usual] (1,1) to (1,1.5)node[above]{$k{+}l{+}m$};
\end{tikzpicture}
=
\begin{tikzpicture}[anchorbase,scale=1]
\draw[usual] (0,0)node[below]{$k$} to (1,1);
\draw[usual] (1,0)node[below]{$l$} to (1.5,0.5);
\draw[usual] (2,0)node[below]{$\rsym$} to (1,1);
\draw[usual] (1,1) to (1,1.5)node[above]{$k{+}l{+}m$};
\end{tikzpicture}
,\quad
\begin{tikzpicture}[anchorbase,scale=1]
\draw[usual] (0,-1)node[below]{$k$} to (1,0) to (0.5,0.5);
\draw[usual] (1,-1)node[below]{$l$} to (0,0) to (0.5,0.5);
\draw[usual] (0.5,0.5) to (0.5,1)node[above]{$k{+}l$};
\end{tikzpicture}
=(-1)^{kl}
\begin{tikzpicture}[anchorbase,scale=1]
\draw[usual] (0,0)node[below]{$l$} to (0.5,0.5);
\draw[usual] (1,0)node[below]{$k$} to (0.5,0.5);
\draw[usual] (0.5,0.5) to (0.5,1)node[above]{$k{+}l$};
\end{tikzpicture}
,
\end{gather*}
then ensure that the crossings that end up in the sandwiched part 
are given by shortest coset representatives for $\sandtwo{K}{L}$.
\end{proof}

\begin{remark}\label{R:DiagramsAtInfinityTwo}
Our strategy to construct $\sand{K}{L}$ is borrowed from 
semigroup and monoid theory where similar constructions known under the slogan of 
\emph{Green relations or cells}, see for example \cite[Section 4]{Tu-sandwich-cellular}. Explicitly 
and in the spirit of sandwich cellularity, 
\cite{Br-algebra-orthogonal} worked out a semisimple version of \cref{P:DiagramsAtInfinity} 
for the Brauer algebra (which sits inside the category $\webb$ by \cref{P:DiagramsBrauer}).

Constructions similar to $\sandtwo{K}{L}$ are sometimes known as \emph{chicken feet bases} 
and have appeared in several different contexts in the literature, 
see {e.g.} \cite[Definition 5.26]{StWe-quiver-schur}, 
\cite[Proof of Theorem 1.10]{RoTu-symmetric-howe}, 
\cite[Section 4]{BrEnAiEtOs-semisimple-tilting} or 
\cite[Section 3.5]{CaKuMu-webs-type-p}.

Neither of these should be confused with cellular 
or light-ladder-type bases as in \cite{AnStTu-cellular-tilting}, 
\cite{El-ladders-clasps} or \cite{Bo-c2-tilting}.
\end{remark}

%%%%%%%%%%%%%%%%

\subsection{Orthogonal semisimplifications and colored Brauer diagrams}\label{SS:SemisimpleBrauer}

%%%%%%%%%%%%%%%%

Recall that $\F$ is an infinite field over $\A$. As usual 
in modular representation theory, let $p=\mathrm{char}\,\F\in\{3,5,7,\dots\}\cup\{\infty\}$ with $p=\infty$ in case the characteristic of $\F$ is zero.

\begin{lemma}\label{L:SemisimpleBrauer}
If $p>\lsym$, then there is an equivalence of symmetric ribbon $\F$-linear categories
\begin{gather*}
\semisimple{\Kar[\brauertwo_{\lsym}]}\to\semisimple{\TiltF}.
\end{gather*}
\end{lemma}

\begin{proof}
Our hypothesis on $p$ implies that $\ext[k]$ appears as a direct summand of $V^{\otimes k}$ for $k\in[0,\lsym]$. Thus $\Kar[\FundF]$ is equivalent to the additive idempotent completion of the full monoidal subcategory generated by $V$.

The universal property of $\Kar[\brauertwo_{\lsym}]$ gives us a functor $\Kar[\brauertwo_{\lsym}]\to\Rep_{\F}\big(\ogroup\big)$ sending the generating object to $V$. This functor is full, see \cite[Section 7]{CoPr-invariant-theory}, so \cref{L:SemisimpleNeg} implies there is an equivalence between $\semisimple{\Kar[\brauertwo_{\lsym}]}$ and the additive idempotent completion of the full monoidal subcategory generated by $V$. 

The result then follows from \cref{P:Kar-Fund-is-Tilt}.
\end{proof}

\begin{lemma}\label{L:ss-coloredBrauer}
If $p>d_{i}$ for $i\in\{0, \dots, r-1\}$, then
\begin{gather*}
\semisimple{\Kar[\brauertwo_{(d_{i})}]}
\to
\boxtimes_{i=0}^{r-1}
\semisimple{\Kar[\brauertwo_{d_{i}}]}
\to
\boxtimes_{i=0}^{r-1}
\semisimple{\catstuff{Tilt}_{\F}(\ogroup[d_{i}])}
\end{gather*}
are equivalences of symmetric ribbon $\F$-linear categories.
\end{lemma}

\begin{proof}
Combine \cref{L:ColorBrauerEquivalence} and \cref{L:SemisimpleBrauer}.
\end{proof}

Now we drop the assumption that $p>\lsym$ and study the semisimplification for $\TiltF$. 

\begin{definition}
For $k\in\N$ we use $(k)_{p}=(k_{0},k_{1},\dots)$
to denote the \emph{$p$-adic digits} of $k$ as in 
\cref{S:Intro}.

Let $x,y\in\N$. Then define $x\leq_{p} y$ 
if $(x)_{p}$ is less than or equal to $(y)_{p}$ entrywise, meaning $x_{i}\leq y_{i}$ for all $i\in\N$. 
\end{definition}

Because of \cref{P:Kar-Fund-is-Tilt}, $\TiltF$ is the additive idempotent completion of the category $\FundF$, which is monoidally generated by the exterior powers $\ext[k]$ for $k\in[0,\lsym]$. In fact, after passing to the semisimplification something stronger is true.

\begin{remark}
The main player below is \emph{Lucas' theorem}:
\begin{gather*}
\binom{a}{b}\equiv
\prod_{i\in\N}\binom{a_{i}}{b_{i}}.
\end{gather*}
Here we again use $p$-adic digits $(a_{i})$ and 
$(b_{i})$ for $a,b\in\N$.
\end{remark}

\begin{lemma}\label{L:gen-by-Lambda^p's}
The category $\semisimple{\Kar[\FundF]}$ is $\hcirc$-generated 
by $\ext[p^{i}]$ for $i\geq 0$.     
\end{lemma}

\begin{proof}
For $p=\infty$ there is nothing to show, so assume $p<\infty$.
The same argument in the proof of \cite[Lemma 3.4]{BrEnAiEtOs-semisimple-tilting} shows that
\begin{enumerate}[label=\arabic*.]

\item If $k\nleq_{p}\lsym$, then 
Lucas' theorem implies that $p$ divides 
$\dim_{\F}\ext[k]$ and therefore 
$\semisimple{\ext[k]}\cong 0$. 

\item If $k\leq_{p}\lsym$, then $\ext[k]$ is a direct summand of $\bigotimes_{i\in\N}(\ext[p^{i}])^{\otimes k_{i}}$.

\end{enumerate}
We conclude that $\semisimple{\Kar[\FundF]}$ is 
$\otimes$-generated by the claimed exterior powers.
\end{proof}

\cref{L:gen-by-Lambda^p's} tells us that every object in $\semisimple{\Kar[\FundF]}$ is a sum of summands of tensor products of $p^{i}$th exterior powers. The following lemma helps us understand morphisms between tensor products of $p^{i}$th exterior powers.  

\begin{lemma}\label{L:merge-split-negligible}
Let $a,b\in\Z_{>0}$ such that $a+b=p^{i}$. The morphisms $\Ya_{p^{i}}^{a,b}$ and {\upshape$\Pa_{a,b}^{p^{i}}$} are zero in $\semisimple{\Kar[\FundF]}$.
\end{lemma}

\begin{proof}
Suffices to show the merge in split morphisms are in the negligible ideal. To this end, use \cref{P:DiagramsFromAToO}
and then we can use the same argument as for the proof of \cite[Lemma 4.16]{BrEnAiEtOs-semisimple-tilting}.
\end{proof}

\begin{lemma}\label{L:dim-Lambdapi-mod-p}
Using $p$-adic digits, in $\webF$ we have
\begin{gather*}
\begin{tikzpicture}[anchorbase,scale=1]
\draw[usual] (0,0)node[left]{$k$} to[out=270,in=180] (0.5,-0.5) to[out=0,in=270] (1,0);
\draw[usual] (0,0) to[out=90,in=180] (0.5,0.5) to[out=0,in=90] (1,0);
\end{tikzpicture}=\lsym_{i}.
\end{gather*}
\end{lemma}

\begin{proof}
Another important consequence of Lucas' theorem is that
\begin{gather*}
\dim_{\mathbf{Fund}_{\F}\big(\ogroup\big)}\ext[p^{i}]=\lsym_{i},
\end{gather*}
where we use the $p$-adic digits $\lsym_{i}$. In particular, the $p^{i}$ labeled circle in $\webF$ is equal to $\lsym_{i}$ times the empty diagram. 
\end{proof}

It follows from \cref{L:dim-Lambdapi-mod-p} that there is a symmetric ribbon $\F$-linear functor 
\begin{gather*}
\brauertwo_{(N)_{p}}\to\webF
\end{gather*}
which sends the color $i$ generating object -- which has dimension $N_{i}$ in $\brauertwo_{(N)_{p}}$ -- to the generating object in $\webF$ labeled $p^{i}$. Crossings colored $i$ and $j$ are sent to crossings labeled $p^{i}$ and $p^{j}$, while cups and caps colored $i$ are sent to cups and caps labeled $p^{i}$. Composing with the functor 
\begin{gather*}
\webF\to\FundF\to\semisimple{\Kar[\FundF]},
\end{gather*}
then taking additive idempotent completion of $\brauertwo_{(N)_{p}}$ we get the following.

\begin{lemma}\label{L:Br-to-ssFund-ess-surj}
The functor $\Kar[\brauer_{(N)_{p}}]\to\semisimple{\Kar[\FundF]}$ is essentially surjective.
\end{lemma}

\begin{proof}
Immediate from \cref{L:gen-by-Lambda^p's}.
\end{proof}

We now prove our second main theorem:

\begin{proof}[Proof of \cref{T:O-ss}]
Our argument is analogous to the proof of \cite[Theorem 4.17]{BrEnAiEtOs-semisimple-tilting}. It follows from \cref{P:DiagramsAtInfinity} and \cref{L:merge-split-negligible} that the functor
\begin{gather*}
\Kar[\brauertwo_{(N)_{p}}]\twoheadrightarrow\semisimple{\Kar[\FundF]}
\end{gather*}
is full. It follows then from \cref{L:Br-to-ssFund-ess-surj} and \cref{L:SemisimpleNeg} that there is an equivalence
\begin{gather*}
\semisimple{\Kar[\brauertwo_{(N)_{p}}]}\xrightarrow{\cong}\semisimple{\Kar[\FundF]}.
\end{gather*}
Thus, we have a chain of equivalences
\begin{gather*}
\scalebox{0.95}{$\boxtimes_{i=0}^{r-1}
\semisimple{\TiltFi}\xleftarrow{\cref{L:ss-coloredBrauer}} 
\semisimple{\Kar[\brauertwo_{(N)_{p}}]}
\to
\semisimple{\Kar[\FundF]}
\xrightarrow{\cref{P:Kar-Fund-is-Tilt}} 
\semisimple{\TiltF}.$}
\end{gather*}
The proof is complete.
\end{proof}

\begin{remark}
By \cref{T:O-ss},
$\semisimple{\TiltFi}$ has finitely many simple objects 
if and only if all $p$-adic digits are not $2$. To see this note 
that $\Rep_{\F}\big(\ogroup[2]\big)$ contributes 
infinitely many simple objects, while all other cases contribute 
finitely many simple objects.
\end{remark}

%%%%%%%%%%%%%%%%

\section{Background: highest weight categories for orthogonal groups}\label{S:BackgroundOrthogonal}

%%%%%%%%%%%%%%%%

This section summarizes the highest weight theory 
of the orthogonal group, and also of the special orthogonal group.
The former is difficult to find in the literature since it is 
not simply connected, so we decided to give the details although 
the material is well-known to experts.

We will work over $\Aa$ as in \cref{N:GeneralA} whose faction field is $\Q$, and then switch to $\A$ and $\F$ for the orthogonal group.

%%%%%%%%%%%%%%%%

\subsection{Tilting representations in general}\label{SS:RepsTilting}

%%%%%%%%%%%%%%%%

The following can be found in many works, e.g. \cite{Do-tilting-alg-groups} or \cite{Ri-good-filtrations}. Also the appendix of \cite{Do-q-schur} covers lot of material relevant for us, and so does \cite{Ja-reps-algebraic-groups}. See also the additional material to \cite{AnStTu-cellular-tilting}, and the setting in \cite{BrSt-semi-infinite-highest-weight} that we will use from time to time.

\subsubsection{Integral representation theory for semisimple groups}

Let $\mathfrak{g}$ be a semisimple Lie algebra over $\Q$. We have the algebra $\uplain[\Aa]=\uplain[\Aa](\mathfrak{g})$ which is the $\Aa$-subalgebra of $\uplain[\Q]=\uplain[\Q](\mathfrak{g})$ generated by $e_{\alpha}^{(a)}$, 
$f_{\alpha}^{(a)}$, and $\binom{h_{\alpha}}{c}$ 
for all $\alpha\in\Delta$, and $a,b,c\in\N$. Also associated to $\mathfrak{g}$ is a simply connected semisimple group scheme $G_{\Aa}$.
Write $\Rep_{\Aa}=\Rep(G_{\Aa})$ for its category of free finite rank $G_{\Aa}$-representations ($\Rep$ means in general free finite rank representations). Such a representation gives rise to a free $\Aa$-module of finite rank with an action of $\uplain[\Aa]$. This gives rise to a fully faithful monoidal functor $\Rep_{\Aa}\to\Rep_{\Aa}\uplain[\Aa]$.

The \emph{Chevalley involution} $\omega\colon\uplain[\Q]\rightarrow \uplain[\Q]$, which swaps $e_{\alpha}$ and $f_{\alpha}$, and negates $h_{\alpha}$, also preserves $\uplain[\Aa]\subset\uplain[\Q]$. We also denote the restriction of $\omega$ to $\uplain[\Aa]$ by $\omega$. Given a $\uplain[\Aa]$-representation $\rsym$, we obtain another $\uplain[\Aa]$-representation, denoted $M^{\omega}$, by twisting the action of $\uplain[\Aa]$ by $\omega$, i.e. $\rho_{M^{\omega}}=\rho_{\rsym}\circ\omega$. If $\rsym$ is a free finite rank $\Aa$-module, then $M^{\pivo}=\Hom_{\Aa}(M,\Z)^{\omega}$ is too. Moreover, since $\omega^{2}=1$, the natural identification of $\rsym$ with its double dual gives a canonical isomorphism of $\uplain[\Aa]$-representations: $M\cong M^{\pivo\pivo}$. 

Let $\wta\in X^{+}=X^{+}(\mathfrak{g})$. Then, after choosing a Borel subalgebra $B$, $\uplain[\Aa]$ has \emph{induced representations} $\acoweyl[\wta]=\textbf{Ind}_{B}^{G}(-\wta)$, and \emph{Weyl representations} $\aweyl=\acoweyl[\wta]^{\pivo}$. Since $\acoweyl[\wta]\cong\acoweyl[\wta]^{\pivo\pivo}=\aweyl[\wta]^{\pivo}$, we also refer to induced representations as \emph{dual Weyl representations}.

If $\qsimple[\wta]$ is the simple $\uplain[\Q]$-representation with fixed highest weight vector $v_{\wta}^{+}$, of weight $\wta$, then $\aweyl[\wta]\cong\uplain[\Aa]\cdot v_{\wta}\subset\qsimple[\wta]$. In fact, $\aweyl[\wta]$ is a free $\Aa$-module which is a direct sum of its weight spaces, and therefore has a \emph{character}. This character is equal to the character of $\qsimple[\wta]$, which in turn is given by \emph{Weyl's character formula}. From $\acoweyl[\wta]=\aweyl[\wta]^{\pivo}$, we find that $\acoweyl[\wta]$ is also a free $\Aa$-module, which is a direct sum of its weight spaces, and has character given by the Weyl character formula.

For each $\wta\in X^{+}$, there is a unique $\uplain[\Aa]$-representation homomorphism, $\hs_{\wta}\colon\aweyl[\wta]\to\acoweyl[\wta]$ such that $v_{\wta}^{+}\mapsto (v_{\wta}^{+})^{\pivo}$. We may write $\hs$ in place of $\hs_{\wta}$ if the weight is clear form context. Let $\wta,\wtb\in X^{+}$. Then we have \emph{Ext-vanishing}:
\begin{gather*}
\Ext^{i}\big(\aweyl[\wta],\acoweyl[\wtb]\big)=
\begin{cases}
\Aa\cdot\hs & \text{if $i=0$ and $\wta=\wtb$}, 
\\
0 & \text{otherwise}.
\end{cases}
\end{gather*}
Write $\Wcat_{\Aa}$ to denote the full subcategory of objects in $\Rep_{\Aa}$ which admit a filtration by Weyl representations. Similarly, write $\DWcat_{\Aa}$ for the full subcategory with object admitting filtrations by dual Weyl representations. Define
\begin{gather*}
\Tilt_{\Aa}=\Tilt_{\Aa}(G)
=\Wcat_{\Aa}\cap\DWcat_{\Aa}.
\end{gather*}
The objects of this category are \emph{tilting representations}.

\begin{remark}\label{R:Kaneda}
Using Lusztig's work on canonical bases for quantum groups \cite[Part IV]{Lu-intro-quantum-groups}, \cite{Kan-based-modules-filtrations} shows that each of the subcategories $\Wcat_{\Aa}$, $\DWcat_{\Aa}$, and $\Tilt_{\Aa}$, is closed under tensor product. Over a field this results is Paradowski's \cite{Pa-tilting-tensor}.
\end{remark}

\subsubsection{Highest weight theory for reductive groups}

Let $\Ff$ be a field. We can define all the notions from 
above by scalar extension from $\Aa$ to $\Ff$, and we also get a
fully faithful monoidal functor $\Rep_{\Ff}\to\Rep_{\Ff}\uplain[\Ff]$.

All the results of the previous section still hold over $\Ff$. But now, there is more since we can talk about \emph{simple representations}. For each $\wta\in X^{+}$, there is a finite dimensional simple representation $\ffsimple[\wta]$, which is the unique simple quotient representation of $\ffweyl[\wta]$ and the unique simple subrepresentation of $\ffcoweyl[\wta]$. This implies that the map $\hs$ factors nontrivially through $\ffsimple[\wta]$. Moreover, the set $\{\ffsimple[\wta]|\wta\in X^{+}\}$ is a complete and irredundant set of simple objects in $\Rep_{\Ff}$. 

Recall that the usual partial order on weights $X$ is defined by 
$\wtb\leq\wta$ if $\wta-\wtb$ is an $\N$-linear combination of positive roots. If $\wta\in X$, then we write 
\begin{gather*}
X(\leq\wta)=\{\wtb\in X|\wtb\leq\wta\},
\end{gather*}
and $X^{+}(\leq\wta)=X(\leq \wta)\cap X^{+}$. Here $X^{+}$ means dominant (integral) weights. 

For $\pi\subset X^{+}$ and $M$ in $\Rep_{\Ff}$, let $M_{\pi}$ be the largest subrepresentation with all composition factors isomorphic to $\ffsimple[\wtb]$, where $\wtb\in\pi$. Define $\Rep_{\Ff}(\pi)$ to be the full subcategory of $\Rep_{\Ff}$ with objects $M=M_{\pi}$. In the case that $\pi=X^{+}(\leq\wta)$, we simply write $M_{\leq\wta}$ and $\Rep_{\Ff}(\leq\wta)$.

Fix $\wta\in X^{+}$. Since $\ffsimple[\wtb]=\ffsimple[\wtb]_{\leq\wtb}$ and Weyl and dual Weyl representations have the same character as $\qsimple[\wta]$, it is easy to see that $\ffweyl[\wta]=\ffweyl[\wta]_{\leq\wta}$ and $\ffcoweyl[\wta]=\ffcoweyl[\wta]_{\leq\wta}$. One can show that
\begin{gather*}
\Ext^{>0}_{\Rep_{\Ff}(\leq\wta)}\big(\ffweyl[\wta],\ffsimple[\wta]\big)=0 
,\quad
\Ext^{>0}_{\Rep_{\Ff}(\leq\wta)}\big(\ffsimple[\wta],\ffcoweyl[\wta]\big).
\end{gather*}
Since $\ffweyl[\wta]$ has $\ffsimple[\wta]$ as its unique simple quotient, we can say that $\ffweyl[\wta]$ is a projective cover of $\ffsimple[\wta]$ in $\Rep_{\Ff}(\leq\wta)$. Similarly, $\ffcoweyl[\wta]$ is an injective hull of $\ffsimple[\wta]$ in $\Rep_{\Ff}(\leq\wta)$.

Thus, the category $\Rep_{\Ff}$ is a (semi-infinity) \emph{highest weight category}. It follows from that we have the following classification of indecomposable objects in $\Tilt_{\Ff}$:
\begin{enumerate}[label=(\roman*)]

\item For each $\wta\in X^{+}$, there is an \emph{indecomposable tilting representations}, denoted $\fftilting[\wta]$, which has as part of its Weyl filtration a subrepresentation $\ffweyl[\wta]\to \fftilting[\wta]$.

\item If $\wta\neq\wtb$, then $\fftilting[\wta]\not\cong\fftilting[\wtb]$.

\item Every indecomposable object in $\Tilt_{\Ff}$ 
is of the form $\fftilting[\wta]$ for some $\wta\in X^{+}$.

\end{enumerate}

%%%%%%%%%%%%%%%%

\subsubsection{Saturated sets}\label{SS:BackgroundRepBasics}

%%%%%%%%%%%%%%%%

It is not difficult to check that $\Rep_{\Ff}(\leq\wta)$ is also a highest weight category. We want to generalize this property for other subsets $\pi\subset X^{+}$. The key is the following: 

\begin{definition}\label{D:BackgroundSaturated}
A set $S\subset X$ is \emph{saturated} 
if for all $\wta\in X$ and for all 
$\alpha\in\Phi^{+}$, then
\begin{gather*}
\wta-i\alpha\in S 
\quad\text{when } 
\begin{cases} 
0\leq i\leq \alpha^{\vee}(\wta) &
\text{if }\alpha^{\vee}(\wta)\geq 0,
\\ 
\alpha^{\vee}(\wta)\leq i\leq 0 &
\text{if }\alpha^{\vee}(\wta)<0.
\end{cases}
\end{gather*}
Here $\Phi^{+}$ denotes the set of positive roots.
\end{definition}

\begin{lemma}\label{L:BackgroundSaturated}
Saturated sets are invariant under unions, intersections,
and under the action of the Weyl group $W$ associated to $\mathfrak{g}$.
\end{lemma}

\begin{proof}
The first two claims are immediate. Moreover, 
since $W$ is generated by reflections, the claim follows from considering the formula for the reflection perpendicular to $\alpha\in\Phi$, $s_{\alpha}(\wta)=\wta-\alpha^{\vee}(\wta)\cdot\alpha$, while noting that $s_{\alpha}=s_{-\alpha}$.
\end{proof}

The prototypical example of a saturated set of 
weights is the set of weights in a Weyl 
representation. Moreover, from \cite[Theorem 1.9]{St-partial-order-dom-wts} 
we have
\begin{gather*}
\wt V(\wta)=W\cdot X^{+}(\leq\wta).
\end{gather*}
Here $\wt\qsimple[\wta]$ denotes the weights of the simple highest weight $\mathfrak{g}(\Q)$-representation $\qsimple[\wta]$.
This suggests the following definition.

\begin{definition}\label{D:BackgroundSaturatedTwo}
A set of dominant weights $\pi\subset X^{+}$ is \emph{saturated} if for all $\wta\in\pi$, we have $X^{+}(\leq\wta)\subset\pi$.
\end{definition}

\begin{example}
The prototypical example of a set of dominant weights which is saturated is $X^{+}(\leq\wta)$. 
\end{example}

\begin{lemma}\label{L:BackgroundDominant}
If $S\subset X$ is a saturated set of weights, then $S\cap X^{+}$ is a saturated set of dominant weights.
\end{lemma}

\begin{proof}
See \cite[Proof of Lemma 1.8]{St-partial-order-dom-wts}.
\end{proof}

\begin{remark}\label{R:BackgroundSaturated}
The definition of a set of weights being saturated is
classical, see e.g. \cite[Exercises VI.1.23-24 and VI.2.5]{Bo-chapters-4-6}.
The notion of a set of dominant weights being saturated 
came later, see e.g. \cite[Definition A3]{Do-q-schur}.
\end{remark}

\begin{proposition}\label{P:sat-is-high-wt}
If $\pi\subset X^{+}$ is saturated, then $\Rep_{\Ff}(\pi)$ is a highest weight category, with indexing set $\pi$ and partial order induced from $X^{+}$ by $\pi\subset X^{+}$. 
\end{proposition}

\begin{proof}
This is \cite[Proposition A3.4]{Do-q-schur}.
\end{proof}

Finally, we state a Lemma which makes it easy to verify certain sets are saturated. This lemma is comparable to 
\cite[Proposition 1.3.2]{DoGiSu-gen-schur}.

\begin{lemma}\label{L:check-saturated}
Suppose $V$ is any finite dimensional 
$\mathfrak{g}(\Q)$-representation, then $\wt V$ is saturated. 
\end{lemma}

\begin{proof}
Since $V$ is completely reducible, we have
\begin{gather*}
\wt V=\wt\bigoplus_{\wta\in X^{+}}\qsimple[\wta]^{\oplus [V:\qsimple[\wta]]}=
\bigcup_{\wta\in X^{+},[V:\qsimple[\wta]]\neq 0}
\wt\qsimple[\wta],
\end{gather*}
The claim follows from observing that each 
$\wt\qsimple[\wta]$ is saturated and being saturated 
is closed under unions, cf. \cref{L:BackgroundSaturated}.
\end{proof}

\begin{remark}\label{R:check-saturated}
This criterion is particularly useful when we have a finite dimensional representation in $\Rep_{\Ff}$ which comes from a representation over $\uplain[\Aa]$, since then we can extend scalars from $\Aa$ to $\Q$ to verify the weight spaces are saturated.
\end{remark}

\subsubsection{Schur algebras}\label{SS:SchurAlgebra}

\cref{P:sat-is-high-wt} suggests the following definition.

\begin{definition}
The \emph{generalized Schur algebra} 
associated to a saturated set of dominant weights $\pi\subset X^{+}$, denoted $\spi[\Aa](\mathfrak{g})$, or $\spi[\Aa]$ if $\mathfrak{g}$ is understood, is defined as the quotient of $\udot[\Aa]$ by the ideal generated by $\one_{\chi}$ for all 
$\chi\notin W\cdot\pi$.
\end{definition}

The algebra $\spi[\Aa]$ is an associative algebra with unit $\one_{\pi}=\sum_{\chi\in W\cdot\pi}\one_{\chi}$.

\begin{remark}
If $V$ is a $\uplain[\Ff]$-representation, then $V_{\pi}$ is naturally a representation over $\spi[\Ff]$. In fact, using \cite[Proposition A3.2(ii)]{Do-q-schur}, one finds there is an equivalence of additive $\Ff$-linear categories $\Rep_{\Ff}(\pi)\cong\Rep_{\Ff}\spi[\Ff]$.    
\end{remark}

The canonical basis $\mathbb{B}$ for $\dot{U}_{\Aa}$ descends to a canonical basis $\mathbb{B}[\pi]=\coprod_{\wta\in\pi}\mathbb{B}[\wta]$, where $\mathbb{B}[\wta]$ is as defined in \cite[29.1]{Lu-intro-quantum-groups}. This renders $\spi[\Aa]$ a 
\emph{based module}, as a left representation over $\dot{U}_{\Aa}$, and therefore $\spi[\Aa]$ has a filtration by Weyl representations, see \cite[Section 27.1.7]{Lu-intro-quantum-groups}. 

A representation with a Weyl filtration will always embed into a tilting representations, cf. \cite[Lemma 5B.11]{BoTu-symplectic-howe}. In particular, $\spi[\Ff]$ embeds in a tilting representations.

We learned the following key lemma from \cite{AdRy-tilting-howe-positive-char}:

\begin{lemma}\label{L:AR-lemma}
A full tilting representations for $\spi[\Ff]$ is faithful. 
\end{lemma}

\begin{proof}
See \cite[Proposition 5B.13]{BoTu-symplectic-howe}.
\end{proof}

Complete reducibility of finite dimensional representations over $\uplain[\Q]$ implies that
\begin{gather*}
\spi[\Q]\cong\prod_{\wta\in\pi}
\End\big(\qsimple[\wta]\big).
\end{gather*}
Since $\qsimple[\wta]\cong\qweyl[\wta]$, and $\aweyl[\wta]$ has the same formal character as $\qweyl[\wta]$, it follows that $\spi[\Aa]$ has Weyl character 
\begin{gather}\label{E:schuralgebras-weyl-character}
\big(\spi[\Aa]:\aweyl[\wta]\big)=
\begin{cases}
\rank_{\Aa}\aweyl[\wta] & \text{if $\wta\in\pi$},
\\
0 & \text{otherwise}.
\end{cases}
\end{gather}
The analog equality then follows for $\spi[\Ff]$, since $\spi[\Aa]$ is a free $\Aa$-module with basis $\mathbb{B}[\pi]$. 

%\subsubsection{Howe duality for Schur algebras}

%%%%%%%%%%%%%%%%

\subsection{Tilting representations for orthogonal groups}\label{SS:RepsTiltingO}

%%%%%%%%%%%%%%%%

The orthogonal group is disconnected, with identity component the special orthogonal group and component group $\Z/2\Z$. However, the usual theory of tilting representations for connected reductive groups can be modified as follows. First, following \cite{AcHaRi-disconnected-groups} (taking a slightly different perspective in some places), we describe how to think about representations of the orthogonal group as a highest weight category. General theory from \cite{BrSt-semi-infinite-highest-weight} then implies the existence of tilting representations for orthogonal groups. 

%%%%%%%%%%%%%%%%%%%%

\subsubsection{Representations of $\ogroup$}\label{SS:highestweights-for-on}

%%%%%%%%%%%%%%%%%%%%

Following the conventions from before:

\begin{notation}
We will write $\Rep_{\A}\big(\ogroup[\lsym]\big)$ and $\Rep_{\A}\big(\sogroup[\lsym]\big)$ for the respective categories of finite dimensional representations over $\A$. We also use similar notation
that should be easy to guess from the context.
\end{notation}

The following defines an involutive algebra automorphism
as one easily checks:

\begin{definition}\label{D:U(oN)-defn}
We define a map $\sigma\colon\uplain(\mathfrak{so}_{\lsym})\to\uplain(\mathfrak{so}_{\lsym})$ by:
\begin{enumerate}[label=\arabic*.]

\item When $\lsym=2\llsym+1$ we let
\begin{gather*}
e_{i}\mapsto e_{i},
f_{i}\mapsto f_{i}, 
h_{i}\mapsto h_{i}\text{ for $i\in[1,\lsym-1]$},
\\
e_{\lsym}\mapsto -e_{\lsym},
f_{\lsym}\mapsto -f_{\lsym},
h_{\lsym}\mapsto h_{\lsym}.
\end{gather*}

\item When $\lsym=2\llsym$ we let
\begin{gather*}
e_{i}\mapsto e_{i},
f_{i}\mapsto f_{i},
h_{i}\mapsto h_{i}\text{ for 
$i\in[1,\lsym-2]$},
\\
e_{\lsym-1}\mapsto e_{\lsym},
f_{\lsym-1}\mapsto f_{\lsym},
h_{\lsym-1}\mapsto h_{\lsym}.
\end{gather*}

\end{enumerate}
In the even case $\sigma$ 
is the automorphism induced by the type D Dynkin diagram automorphism swapping the fishtail vertices.
\end{definition}

Write $\uplain(\mathfrak{so}_{\lsym})^{\sigma}$ to denote the $\A$-algebra generated by $\uplain(\mathfrak{so}_{\lsym})$ and $\A[\sigma]/(\sigma^{2})$ subject to the relation 
\begin{gather}\label{E:ON-smash-product-reln}
\sigma X\sigma^{-1}=\sigma(X) 
\quad \text{for all $X\in\uplain(\mathfrak{so}_{\lsym})$}.
\end{gather}

\begin{lemma}
As a right $\uplain(\mathfrak{so}_{\lsym})$-representation 
$\uplain(\mathfrak{so}_{\lsym})^{\sigma}$ is freely generated by $1$ and $\sigma$. 
\end{lemma}

\begin{proof}
Boring and omitted.
\end{proof}

We can view a finite dimensional representation of $\sogroup[\lsym]$ as a finite dimensional $\uplain(\mathfrak{so}_{\lsym})$-representation, with a weight decomposition, such that the dominant weights are contained in $X^{+}(\sogroup[\lsym])$ \cite[Sections 7.14--7.17]{Ja-reps-algebraic-groups}. From this perspective, a finite dimensional representation of $\ogroup[\lsym]$ can be viewed as a finite dimensional $\uplain(\mathfrak{so}_{\lsym})^{\sigma}$-representation, with a weight decomposition, such that the dominant weights are contained in $X^{+}(\sogroup[\lsym])$.

The \emph{Chevalley involution} $\omega\colon\uplain(\mathfrak{so}_{\lsym})\to\uplain(\mathfrak{so}_{\lsym})$ swaps $e_{\alpha}$ and $f_{\alpha}$ and negates $h_{\alpha}$. 

\begin{lemma}
The Chevalley involution commutes with $\sigma$ and preserves the relations $\sigma^{2}=1$ and $\sigma X \sigma^{-1}=\sigma(X)$, for $X\in \uplain(\mathfrak{so}_{\lsym})$. Thus, we can extend $\omega$ to an automorphism of $\uplain(\mathfrak{so}_{\lsym})^{\sigma}$.
\end{lemma}

\begin{proof}
Easy and omitted.
\end{proof}

As usual, we can use the Chevalley involution to 
define the \emph{dual $\ogroup[\lsym]$-representation} by $U^{\pivo}=\Hom_{\A}(U,\A)^{\omega}$ where $U\in\Rep_{\A}\big(\ogroup[\lsym]\big)$.
Moreover, suppose that $W$ is a finite dimensional $\sogroup[\lsym]$-representation. View $W$ as a $\uplain(\mathfrak{so}_{\lsym})$-representation. Then we define
\emph{induction and restriction}
\begin{gather*}
\IndN(\placeholder)
\colon\Rep_{\A}\big(\ogroup[\lsym]\big)\leftrightarrows
\Rep_{\A}\big(\sogroup[\lsym]\big)
\colon
\ResN(\placeholder)
\end{gather*}
as follows. Before doing so, note that, given a finite dimensional representation $U$ of $\ogroup[\lsym]$, we obtain a $\uplain(\mathfrak{so}_{\lsym})^{\sigma}$-representation structure on $U$. We then define:
\begin{gather*}
\IndN(W)=\uplain(\mathfrak{so}_{\lsym})^{\sigma}\hcirc_{\uplain(\mathfrak{so}_{\lsym})}W,
\\
\ResN(U)= \functorstuff{R}^{\uplain(\mathfrak{so}_{\lsym})^{\sigma}}_{\uplain(\mathfrak{so}_{\lsym})}(U).
\end{gather*}

\begin{lemma}
Induction and restriction are exact.
\end{lemma}

\begin{proof}
Since $\uplain(\mathfrak{so}_{\lsym})^{\sigma}$ is free as a right $\uplain(\mathfrak{so}_{\lsym})$-representation, $\IndN$ is exact.
A similar argument works for $\ResN$.
\end{proof}

%%%%%%%%%%%%%%%%%%%%

\subsubsection{Dominant weights for $\ogroup$}\label{SS:highestweights-for-on-dom}

%%%%%%%%%%%%%%%%%%%%

We will write
\begin{gather*}
\llsym=
\begin{cases}
\tfrac{\lsym-1}{2} & \text{if $\lsym$ is odd},
\\
\tfrac{\lsym}{2} & \text{if $\lsym$ is even}.
\end{cases}
\end{gather*}
Let $W$ be a finite dimensional $\ogroup[\lsym]$-representation. Then $W$ is naturally an $\sogroup[\lsym]$-represen\-tation, by restriction, and therefore decomposes into a direct sum of weight spaces indexed by $X\big(\sogroup[\lsym]\big)=X(\mathfrak{so}_{\lsym})\cap
\bigoplus_{i=1}^{\lsym}\Z\epsilon_{i}$. If $w\in W$ is an $\sogroup[\lsym]$-weight vector, then we write $\wt_{SO}(w)$ for the corresponding element in $X\big(\sogroup[\lsym]\big)$. The dominant weights of simple $\sogroup$-representations 
are parameterized by the set
\begin{gather*}
X^{+}\big(\sogroup[\lsym]\big)=
\{a_{1}\epsilon_{1}+\dots+a_{\llsym}\epsilon_{\lsym}|a_{i}\in\Z,a_{1} \geq\dots\geq a_{\llsym},
\substack{a_{\llsym}\geq 0\text{ if $\lsym$ is odd} 
\\
a_{\llsym-1}\geq|a_{\llsym}|\text{ if $\lsym$ is even}}\}
\subset X^{+}(\mathfrak{so}_{\lsym}).
\end{gather*}

Weights for $\ogroup[\lsym]$ are pairs of data, the $\sogroup[\lsym]$-weight, and a ``weight'' for $\sigma$:

\begin{definition}\label{D:X+on}
If $U$ is a finite dimensional $\ogroup[\lsym]$-representation and $u\in U$ is a $\sigma$ eigenvector with $\sigma\acts u=\epsilon\cdot u$, then we write $\epsilon_{\sigma}(u)=\epsilon$. If $u\in U$ is not a $\sigma$ eigenvector, then we write $\epsilon_{\sigma}(u)=0$. If $u$ is also a weight vector for $\sogroup$, then we write
\begin{gather*}
\wt_{O}(u)=\big(\wt_{SO}(u),\epsilon_{\sigma}(u)\big),
\end{gather*}
and call it the \emph{$\ogroup[\lsym]$-weight} of $u$.
\end{definition}

The following partial order is taken from \cite[Section 1]{AdRy-tilting-howe-positive-char}.

\begin{definition}\label{D:dominant-orthogonal-wt}
Let $X\big(\ogroup[\lsym]\big)$ be the set of all pairs $(\wta,\epsilon)$ which appear as weights $\wt_{O}(u)$ for $u\in U$, where $U$ ranges over all finite dimensional $\ogroup[\lsym]$-representations. Let $X^{+}\big(\ogroup[\lsym]\big)$ denote the \emph{dominant} weights, that is pairs of the form
\begin{gather*}
(\wta,\pm 1),\text{ for $\wta\in X^{+}\big(\sogroup[\lsym]\big)$, such that $\sigma(\wta)=\wta$},
\\
(\wta,0),\text{ for $\wta\in X^{+}\big(\sogroup[\lsym]\big)$, 
such that $\sigma(\wta)\neq\wta$.}
\end{gather*}
The partial order that we will use is: $(\wta,\epsilon)<(\wtb,\epsilon^{\prime})$ if and only if $\wta<\wtb$ or $\wta<\sigma(\wtb)$.
\end{definition}

Let $\sigma$ denote the generator of the group of automorphisms of the Dynkin diagram for $\mathfrak{so}_{\lsym}$. This is trivial when $\lsym$ is odd and a nontrivial involution when $\lsym$ is even. The Dynkin diagram automorphism induces maps, which we also denote $\sigma$, on all objects which are determined by the $\mathfrak{so}_{\lsym}$ Dynkin diagram. In particular, $\sigma$ acts on $X(\mathfrak{so}_{\lsym})$, preserving the subset $X^{+}(\mathfrak{so}_{\lsym})$. Note that $\sigma$ is the identity when $\lsym$ is odd. When $\lsym$ is even, $\sigma$ acts on $X(\mathfrak{so}_{\lsym})$ by $(a_{1},\dots,a_{\llsym-1},a_{\llsym})\mapsto (a_{1},\dots,a_{\llsym-1},-a_{\llsym})$.

\begin{definition}\label{D:ParN}
Let $\Par$ be the set of all partitions, that is weakly decreasing sequences of elements in $\Z_{\geq 0}$. We identify partitions with their Young diagram. Taking the transpose of the Young diagram determines an involution of $\Par$, denoted by $\lambda\mapsto \lambda^{T}$ (the transpose diagram). Define the \emph{dominant $\ogroup$-weights} to be
\begin{gather*}
\ParN=\{\lambda\in\Par|(\lambda^{T})_{1}+(\lambda^{T})_{2}\leq\lsym\}.
\end{gather*}
Let further $\xp\colon X^{+}\big(\sogroup[\lsym]\big)\to\Par$ be defined by $\sum_{i=1}^{\lsym}a_{i}\epsilon_{i}\mapsto 
(a_{1},\dots,a_{\llsym-1},|a_{\llsym}|)$.
\end{definition} 

\begin{remark}\label{R:xp-SON-sigma-compatibility}
If $\wta\in X^{+}\big(\sogroup[\lsym]\big)$, then $\xp(\wta)=\xp(\sigma(\wta))$.
\end{remark}

\begin{remark}
The image of the map $\xp$ is contained in $\ParN$. We saw in \cref{R:xp-SON-sigma-compatibility} that $\xp$ is not injective when $\lsym$ is even, and $\xp$ is injective when $\lsym$ is odd. Moreover, $\xp$ is not surjective. The image of $\xp$ is the subset of $\lambda\in\ParN$ such that $(\lambda^{T})_{1}\leq n$.    
\end{remark}

The set $\ParN$ is not closed undertaking the 
transpose, but there is another involution on this set.

\begin{definition}\label{D:twist}
Define the \emph{twisting} involution on $\ParN$, denoted $\lambda\mapsto \lambda^{tw}$, by $\lambda^{tw}= (\lsym-(\lambda^{T})_{1}, (\lambda^{T})_{2}, \dots)^{T}$.
\end{definition}

In words: the twist of $\lambda$ has the same Young diagram, except the first column is replaced with $\lsym-(\lambda^{T})_{1}$ boxes.

\begin{remark}
The fixed points of $\lambda\mapsto\lambda^{tw}$ are exactly the $\lambda$ such that $(\lambda^{T})_{1}=(\lambda^{T})_{2}$. In particular, if $\lsym$ is odd, then the twisting involution on $\ParN$ does not have any fixed points. 
\end{remark}

\begin{lemma}
We have the following.
\begin{enumerate}[label=\arabic*.]

\item The map $\tau\colon X^{+}\big(\ogroup[\lsym]\big)\to\ParN$
\begin{gather*}
\tau\colon X^{+}\big(\ogroup[\lsym]\big)\to\ParN,
\\
(\wta, +1)\mapsto
(\wta_{1},\dots,\wta_{\lsym}), 
\quad(\wta,-1)\mapsto(\wta_{1},\dots,\wta_{\lsym})^{tw},
\text{ for $\sigma(\wta)=\wta$},
\\
(\wta, 0)\mapsto (\wta_{1}, \dots, \wta_{n-1}, |\wta_{\lsym}|)
\text{ otherwise},
\end{gather*}
is a bijection.

\item The image of the map 
\begin{gather*}
\xp\colon X^{+}\big(\sogroup[\lsym]\big)\to\ParN
\end{gather*}
is a fundamental domain for the twisting involution acting on $\ParN$, and $\sigma(\wta)\neq\wta$ if and only if $\xp(\wta)^{tw}=\xp(\wta)$. 

\end{enumerate}

\end{lemma}

\begin{proof}
Not difficult and omitted.
\end{proof}

Thus, we get:

\begin{remark}\label{R:three-ways-on-wts}
There are three ways to encode a dominant weight for the orthogonal group in the literature and that we use in this paper:
\begin{enumerate}[label=\arabic*.]

\item $(\wta,\epsilon)\in X^{+}\big(\sogroup[\lsym]\big)\times\{0, \pm 1\}$,

\item $(\xp(\wta),\epsilon)\in\Lambda_{+}\times 
\{0,\pm 1\}$, 
\quad\text{and}

\item $\tau(\wta,\epsilon)\in\ParN$.
\end{enumerate}
Which one is more convenient depends on the context. 
\end{remark}

\begin{notation}
For a finite dimensional $\ogroup$-representation $U$, and $(\wta, \epsilon)\in X\big(\ogroup[\lsym]\big)$, we write
\begin{gather*}
[(\wta, \epsilon)]U=\{u\in U|\wt_{O}(u)=(\wta, \epsilon)\} \end{gather*}
for the $(\wta,\epsilon)$-weight space of $U$. 
For $\lambda\in\ParN$, there is a corresponding 
$(\wta,\epsilon)\in X^{+}\big(\ogroup[\lsym]\big)$, and we will write
\begin{gather*}
[\lambda]U=[(\wta, \epsilon)]U,
\end{gather*}
for any finite dimensional $\ogroup$-representation $U$.
\end{notation}

We define the $\ogroup$ partial order on $\ParN$ as follows.

\begin{definition}\label{D:O-dominanceorder}
Suppose $\lambda,\mu\in\ParN$ correspond 
to $(\wta,\epsilon)$ and $(\wtb,\epsilon^{\prime})$, 
respectively. Let 
\begin{gather*}
\mu\ltO\lambda\text{ if }\wtb<\wta,\text{ or }\sigma(\wtb)<\wta,
\end{gather*}
where $<$ on $X^{+}\big(\sogroup[\lsym]\big)$ is the restriction of the usual partial order on $X(\mathfrak{so}_{\lsym})$.
\end{definition}

\begin{lemma}
\cref{D:O-dominanceorder} defines a partial order.
\end{lemma}

\begin{proof}
This follows from \cite[Lemma 3.1]{AcHaRi-disconnected-groups}.
\end{proof}

\begin{example}
Consider $\lambda=(0^{\lsym})$ and $\mu=(1^{\lsym})$. 
In this case, both $\lambda$ and $\mu$ correspond 
via $\xp$ to $0\in X^{+}\big(\sogroup[\lsym]\big)$. 
However, since $\lambda \neq\mu$, the 
partitions are not comparable with respect 
to $\ltO$. On the other hand, we have 
$(0^{\lsym})\ltO (1^{2},0^{\lsym-2})$ and $(1^{\lsym})\ltO 
(1^{2}, 0^{\lsym-2})$. 
\end{example}

\begin{lemma}\label{R:characterizing-ON-dominance-order}
We have $\mu<\lambda$ if and only if one of the following holds:
\begin{gather*}
\wtb<\wta,
\text{ or }
\wtb<\sigma(\wta),
\text{ or }
\sigma(\wtb)<\wta,
\text{ or }\sigma(\wtb)<\sigma(\wta). 
\end{gather*}
\end{lemma}

\begin{proof}
Since $\sigma$ is a Dynkin diagram automorphism, it preserves the usual partial order on the set $X(\mathfrak{so}_{\lsym})$, and therefore 
\begin{gather*}
\wtb<\wta\text{ if and only if }\sigma(\wtb)<\sigma(\wta).
\end{gather*}
Since $\sigma$ is an involution, we have 
\begin{gather*}
\sigma(\wtb)<\wta\text{ if and only if }\wtb<\sigma(\wta).
\end{gather*}
Thus, the claim follows.
\end{proof}

\begin{definition}\label{D:dom-order-on-partitions}
Given two partitions $\xp$ and $\xp^{\prime}$, we say that $\xp\led \xp^{\prime}$ if $\sum_{i=1}^{k}\xp_{i}\leq\sum_{i=1}^{k} \xp^{\prime}_{i}$ for all $k\geq 0$.
\end{definition}

Note that $(1^{\lsym})\ltO (1^{2}, 0^{\lsym-2})$, but $(1^{2}, 0^{\lsym-2})\ltd(1^{\lsym})$. Thus, the order on partitions from \cref{D:dom-order-on-partitions} is not adapted to $\ParN$. It is however useful to use $\led$ to compare $\wta$ and $\wtb$ in $X^{+}\big(\sogroup[\lsym]\big)$ when considering $\xp(\wta)$ and $\xp(\wtb)$. 

\begin{lemma}\label{L:SON-aleb-implies-YaleYb}
Let $\wta,\wtb\in X^{+}\big(\sogroup[\lsym]\big)$. 
If $\wta\leq\wtb$, then $\xp(\wta)\led\xp(\wtb)$. 
\end{lemma}

\begin{proof}
We prove this for $\llsym=4$, the general case when $\lsym=2\llsym$ is an exercise and when $\lsym=2\llsym+1$ is an easier exercise.

Suppose $\wta,\wtb\in X^{+}\big(\sogroup[8]\big)\subset\oplus_{i=1}^{4}\Z\epsilon_{i}$ and that $\wta\leq\wtb$. Since a $\Z_{\geq 0}$-linear combination of positive roots is a 
$\Z_{\geq 0}$-linear combination of simple 
roots, we have $\wtb-\wta=w\alpha_{1}+x\alpha_{2}+y\alpha_{3}
+z\alpha_{4}$, where $w,x,y,z\in\Z_{\geq 0}$. Thus,
$\wtb=(\wta_{1}+w,\wta_{2}-w+x,\wta_{3}-x+y+z,\wta_{4}-y+z)$
and
\begin{itemize}

\item $\wtb_{1}-\wta_{1}=w\geq 0$,

\item $(\wtb_{1}+\wtb_{2})-(\wta_{1}+\wta_{2})=x \geq 0 $,

\item $(\wtb_{1}+\wtb_{2}+\wtb_{3})-(\wta_{1}+\wta_{2}+\wta_{3})
=y+z\geq 0$, \quad \text{and}

\item $(\wtb_{1}+\wtb_{2}+\wtb_{3}+|\wtb_{4}|)-
(\wta_{1}+\wta_{2}+\wta_{3}+|\wta_{4}|)
=y+z+|\wta_{4}-y+z|-|\wta_{4}|\geq 0$.
\end{itemize}
The last inequality follows from noticing that $|z-y|\leq|z|+|y|=z+y$, so
\begin{gather*}
|\wta_{4}|=|\wta_{4}+(z-y)-(z-y)|\leq
|\wta_{4}+z-y|+|z-y|
\leq|\wta_{4}+z-y|+z+y.
\end{gather*}
The proof is complete.
\end{proof}

\begin{definition}
Let $\lambda,\mu\in\ParN$ correspond to $(\wta,\epsilon)$ and $(\wtb, \epsilon^{\prime})$, respectively. We define a partial order on $\ParN$ by declaring $\lambda<\mu$ if $\xp(\wta)\ltd\xp(\wtb)$. 
\end{definition}

\begin{lemma}\label{L:leO-implies-le-partition}
If $\lambda\leO\mu$, then $\lambda\leq\mu$.
\end{lemma}

\begin{proof}
Follows from \cref{D:O-dominanceorder}, \cref{L:SON-aleb-implies-YaleYb}, and \cref{R:xp-SON-sigma-compatibility}.
\end{proof}

%%%%%%%%%%%%%%%%

\subsubsection{Standard representations for $\ogroup$}

%%%%%%%%%%%%%%%%  

Since $U$ is a finite dimensional $\ogroup[\lsym]$-representa\-tion. 
It in particular has a weight space decomposition as a representation of $\sogroup[\lsym]$. 

\begin{lemma}
As $\sogroup[\lsym]$-representations we have
\begin{gather*}
\ResN\IndN\big(\weyl[\wta]\big)
\cong
\weyl[\wta]\oplus\weyl[\sigma(\wta)].
\end{gather*}
\end{lemma}

\begin{proof}
It is easy to see from \cref{E:ON-smash-product-reln} that if $u\in U_{\lambda}$, then $\sigma(u)\in U_{\sigma(\lambda)}$, and that if $u\in U$ is annihilated by $\uplain(\mathfrak{so}_{\lsym})^{+}$, then so is $\sigma(u)$. It follows that
\begin{gather*}
1\hcirc v_{\wta}^{+}\mapsto(v_{\wta}^{+},0),
\quad
\sigma\hcirc v_{\wta}^{+}\mapsto(0,v_{\sigma(\wta)}^{+}),
\end{gather*}
is the desired isomorphism.
\end{proof}

Suppose that $\sigma(\wta)=\wta$. Then we let 
\begin{gather*}
v_{(\wta,+1)}^{+}=\tfrac{1}{2}
\left(1\hcirc v_{\wta}^{+}+\sigma\hcirc v_{\wta}^{+}\right)
,\quad
v_{(\wta,-1)}^{+}=\tfrac{1}{2}\left(1\hcirc v_{\wta}^{+}-\sigma\hcirc v_{\wta}^{+}\right).
\end{gather*}

\begin{lemma}
Each $v_{(\wta,\pm 1)}$ generates an $\ogroup[\lsym]$-subrepresentation of $\IndN\big(\weyl[\wta]\big)$.
Moreover,
\begin{gather*}
\uplain(\mathfrak{so}_{\lsym})\cdot v_{(\wta,+1)}
\oplus\uplain(\mathfrak{so}_{\lsym})\cdot v_{(\wta,-1)} 
=\ResN\IndN\big(\weyl[\wta]\big)\xrightarrow{\cong} 
\weyl[\wta]\oplus\weyl[\wta]
\end{gather*}
induces isomorphisms of $\uplain(\mathfrak{so}_{\lsym})$-representations, $\uplain(\mathfrak{so}_{\lsym})\cdot v_{(\wta,\epsilon)}\cong \weyl[\wta]$, for $\epsilon\in\{\pm 1\}$.
\end{lemma}

\begin{proof}
The vectors $v_{(\wta,\pm 1)}$ are highest weight vectors for $\uplain(\mathfrak{so}_{\lsym})$ and eigenvectors, with eigenvalues $\pm 1$ respectively, for $\sigma$. It follows from \cref{E:ON-smash-product-reln}, that for $\epsilon\in\{\pm 1\}$, the $\ogroup[\lsym]$-subrepresenta\-tion generated by $v_{(\wta,\epsilon)}^{+}$ is the $\epsilon$ eigenspace of $\IndN\big(\weyl[\wta]\big)$.

The second claim can then easily be checked.
\end{proof}

Using these lemmas we can make the following definition.

\begin{definition}
We define the \emph{Weyl representation} for $\ogroup$ with highest weight $\lambda\in\ParN$ as
\begin{gather*}
\weyl[\lambda]=
\begin{cases}
\uplain(\mathfrak{so}_{\lsym})\cdot v_{(\wta, +1)} & 
\text{if $\sigma(\wta)=\wta$ and $\lambda=\xp(\wta)$}, 
\\
\uplain(\mathfrak{so}_{\lsym})\cdot v_{(\wta, -1)} & 
\text{if $\sigma(\wta)=\wta$ and $\lambda=\xp(\wta)^{tw}$}
\\
\IndN(\weyl[\wta]) & \text{if $\sigma(\wta)\neq\wta$}.
\end{cases}
\end{gather*}
We define the \emph{dual Weyl representation} with highest weight $\lambda\in\ParN$ as the dual space $\coweyl[\lambda]=\weyl[\lambda]^{\pivo}$.
\end{definition}

There is another natural definition of dual Weyl representation, paralleling our definition of Weyl representation. That is as a summand of a dual Weyl representation for $\sogroup$ induced to $\ogroup$. But in fact one arrives at the same definition this way. 

\begin{lemma}\label{L:dual-Weyl-two-defs}
For $\wta\in X^{+}\big(\sogroup[\lsym]\big)$, we have
\begin{align*}
\IndN(\weyl[\wta])\cong&
\begin{cases}
\weyl[\xp(\wta)]\oplus\weyl[\xp(\wta)^{tw}]
& \text{if $\sigma(\wta)=\wta$}
\\
\weyl[\xp(\wta)] & \text{if $\sigma(\wta)\neq\wta$}.
\end{cases}
\\
\IndN(\coweyl[\wta])\cong&
\begin{cases}
\coweyl[\xp(\wta)]\oplus\coweyl[\xp(\wta)^{tw}] & 
\text{if $\sigma(\wta)=\wta$}
\\
\coweyl[\xp(\wta)] & \text{if $\sigma(\wta)\neq\wta$.}
\end{cases}
\end{align*}
\end{lemma}

\begin{proof}
This is essentially immediate from definitions.
\end{proof}

\begin{lemma}\label{L:Res-standards}
Let $M$ denote either a Weyl or a dual Weyl representation, 
and let $\lambda,\mu\in\ParN$ with $\lambda\neq\mu$.
\begin{gather*}
\ResN\big(M(\lambda)\big)
\cong
\begin{cases}
M(\wta) & \text{if $\xp(\wta)\in\{\lambda,\lambda^{tw}\}$ 
and $\lambda\neq\lambda^{tw}$},
\\
M(\wta)\oplus M\big(\sigma(\wta)\big) & 
\text{if $\xp(\wta)=\lambda$ and $\lambda=\lambda^{tw}$}.
\end{cases} 
\end{gather*}
Moreover, if $M(\lambda)[\mu]\neq 0$, then $\mu\ltO\lambda$.
\end{lemma}

\begin{proof}
Because of the discussion above, it suffices to analyze the decomposition of the $\sogroup[\lsym]$-representation $\ResN\IndN\big(M(\wta)\big)$, for $\wta\in X^{+}\big(\sogroup[\lsym]\big)$. This is analogous to Mackey theory for finite groups, and we leave it to the reader to fill in the details.

The second claim is \cite[Proposition 3.4]{AcHaRi-disconnected-groups}.
\end{proof}

\begin{definition}
For $\lambda\in \ParN$ define a map of $\ogroup$-representations $\hs_{\lambda}\colon\Delta_{\A}(\lambda)\to\nabla_{\A}(\lambda)$ as follows. Suppose $\mathcal{Y}(\wta)=\lambda$. If $\lambda\neq\lambda^{tw}$, then $\hs_{\lambda}:=\hs_{\wta}$, and if $\lambda\neq\lambda^{tw}$, then $\hs_{\lambda}=\IndN(\hs_{\wta})$. 
\end{definition} 

\begin{lemma}\label{L:Ohs-spans}
For $\lambda\in\ParN$, we have $\Hom_{\ogroup}(\Delta_{\wta}(\lambda), \nabla_{\A}(\lambda))=\A\cdot\hs_{\lambda}$.
\end{lemma}

\begin{proof}
Note that if $\lambda\neq\lambda^{tw}$, then $\hs_{\lambda} = \ResN(\hs_{\wta})$, and if $\lambda\neq\lambda^{tw}$, then $\ResN(\hs_{\lambda})=\hs_{\wta}\oplus\hs_{\sigma(\wta)}$. Observing that we have $\sigma(v_{\wta}^{+})=v_{\sigma(\wta)}^{+}$, the claim follows from
the fact that $\hs_{\wta}$ is spanning $\Hom_{\sogroup}(\Delta_{\A}(\wta),\nabla_{\A}(\wta))$.
\end{proof}

\subsubsection{Simple representations for $\ogroup$}

We can view finite dimensional representations 
of $\sogroup[\lsym]$ over $\F$ as $\uplain[\F](\mathfrak{so}_{\lsym})$-representations, with weight space decompositions, such that the weight spaces are contained in $X^{+}\big(\sogroup[\lsym]\big)$. Similarly, we view representations of $\ogroup[\lsym]$ as such $\uplain[\F](\mathfrak{so}_{\lsym})$-representations, with a compatible action of $\sigma$. 

\begin{lemma}
The $\ogroup[\lsym]$-representation $\fweyl[\lambda]$ has a unique maximal $\ogroup[\lsym]$-subrepresenta\-tion, consisting of the sum of subrepresentations $U$ such that $U\cap\fweyl[\lambda]_{\lambda}=0$.
\end{lemma}

\begin{proof}
The usual Yoga.
\end{proof}

Using the previous lemma, we define $\fsimple[\lambda]$ as the unique simple quotient of $\fweyl[\lambda]$. It then follows by duality that $\fcoweyl[\lambda]$ has a simple socle which is isomorphic to $\fsimple[\lambda]$. 

\begin{lemma}
The set $\{\fsimple[\lambda]\}_{\lambda\in \ParN}$ is a complete and irredundant list of the finite dimensional simple $\ogroup[\lsym]$-representations.
\end{lemma}

\begin{proof}
Let $S$ be an simple $\ogroup[\lsym]$-representation. A standard argument coming from Clifford theory shows that $\ResN(S)$ is completely reducible. Thus, it is a direct sum of simple representations for $\sogroup[\lsym]$. Choosing a direct sum decomposition into simple subrepresentations, we then obtain a map to a direct sum of dual Weyl representations. By \cref{L:dual-Weyl-two-defs}, Frobenius reciprocity yields a non-zero map from $S$ to a direct sum of dual Weyl representations for $\ogroup[\lsym]$. Since $S$ is simple, it follows that $S$ is isomorphic to a summand of the socle of this direct sum of dual Weyl representations. Hence, $S\cong\fsimple[\lambda]$ for some $\lambda\in\ParN$. We leave it as an exercise, using highest weights and the action of $\sigma$, to argue that $\fsimple[\lambda]\cong\fsimple[\mu]$ implies $\lambda=\mu$. 
\end{proof}

\begin{lemma}\label{L:res-simples}
Let $\lambda\in\ParN$, then
\begin{gather*}
\ResN\big(\fsimple[\lambda]\big)\cong 
\begin{cases}
\fsimple[\wta] & 
\text{if $\xp(\wta)\in\{\lambda,\lambda^{tw}\}$, and $\lambda\neq\lambda^{tw}$}, 
\\
\fsimple[\wta]\oplus\fsimple[\sigma(\wta)] & 
\text{if $\xi(\wta)=\lambda=\lambda^{tw}$}.
\end{cases}
\end{gather*}
\end{lemma}

\begin{proof}
Since $\fsimple[\lambda]$ is isomorphic to the socle 
of $\fcoweyl[\lambda]$, we get an injective map
\begin{gather*}
\ResN\big(\fsimple[\lambda]\big)\hookrightarrow\ResN\big(\fcoweyl[\lambda]\big).
\end{gather*}
Another standard Clifford theory argument shows that $\ResN\big(\fsimple[\lambda]\big)$ is a completely reducible finite dimensional $\sogroup[\lsym]$-representation. It follows that $\ResN\big(\fsimple[\lambda]\big)$ is a non-zero subrepresentation of the socle of $\ResN\big(\fcoweyl[\lambda]\big)$, which by \cref{L:res-simples} is isomorphic to $\fsimple[\wta]$, if $\xp(\wta)\in\{\lambda,\lambda^{tw}\}$ and $\lambda\neq\lambda^{tw}$, or $\fsimple[\wta]\oplus\fsimple[\sigma(\wta)]$, if 
$\xp(\wta)=\lambda=\lambda^{tw}$. If $\lambda\neq\lambda^{tw}$, the desired result is immediate. If $\lambda=\lambda^{tw}$, so $\sigma(\wta) \neq\sigma(\wta)$, then the desired result follows by looking at the $\sigma$ action on weight spaces. 
\end{proof}

\subsubsection{The orthogonal highest weight category}

Recall that by, for example, \cite[Section 6.4]{BrSt-semi-infinite-highest-weight}, the category $\Rep(\sogroup[\lsym])$, equipped with the poset $(X^{+}\big(\sogroup[\lsym]\big),\leq)$, is an upper finite highest weight category.

\begin{lemma}\label{L:O-delta/nabla-Ext-vanishing}
We have \emph{Ext-vanishing}, i.e.:
\begin{gather*}
\Ext^{i}_{\ogroup[\lsym]}\big(\weyl[\lambda],\coweyl[\mu]\big)
\cong
\begin{cases} 
0 & \text{if $\lambda\neq\mu$ or $i>0$},
\\
\A\cdot\hs & \text{if $\lambda=\mu$ and $i=0$}.
\end{cases}
\end{gather*}
\end{lemma}

\begin{proof}
Since $2$ is invertible in $\A$, one can argue, see \cite[Corollary 3.6.18]{Be-rep-cohomology}, that restriction induces an injective map of $\A$-modules
\begin{gather*}
\Ext^{i}_{\ogroup[\lsym]}\big(\weyl[\lambda],\coweyl[\mu]\big)
\to 
\Ext^{i}_{\sogroup[\lsym]}\bigg(\ResN\big(\weyl[\lambda]\big),
\ResN\big(\coweyl[\mu]\big)\bigg).
\end{gather*}
The usual Ext-vanishing implies that
$\Ext_{\sogroup[\lsym]}^{>0}(\placeholder,\placeholder)=0$ whenever the first entry has a Weyl representation filtration and the second entry has a dual Weyl filtration. It follows from the statement \cref{L:Res-standards} that 
$\Ext_{\ogroup}^{i}(\weyl[\lambda],\coweyl[\mu])\cong 0$ 
for all $\lambda,\mu\in\ParN$. Since 
$\Ext^{0}=\Hom$, and also $\Hom_{\sogroup[\lsym]}(\weyl[\wta],\coweyl[\wtb])\cong 0$, whenever $\wta\neq\wtb$, it suffices to show that
\begin{gather*}
\Hom_{\ogroup[\lsym]}\big(\weyl[\lambda],\coweyl[\lambda^{tw}]\big)\cong0
\text{ when $\lambda\neq\lambda^{tw}$},
\\
\Hom_{\ogroup[\lsym]}\big(\weyl[\lambda],\coweyl[\lambda]\big)
\cong\A\cdot\hs_{\lambda}\text{ when $\lambda=\lambda^{tw}$}.
\end{gather*}
The first equality follows from noting that $\hs_{\wta}$ spans the space of $\sogroup[\lsym]$-homomorphisms over $\A$, here $\xp(\wta)\in\{\lambda,\lambda^{tw}\}$, and $\hs_{\wta}$ does not commute with $\sigma$. The second equality is \cref{L:Ohs-spans}.
\end{proof}

\begin{lemma}
In $\Rep_{\F}\big(\ogroup[\lsym]\big)_{\leq\lambda}$:
The $\ogroup$-representation $\fweyl[\lambda]$, respectively $\fcoweyl[\lambda]$, is the projective cover, respectively the injective hull, of $\fsimple[\lambda]$. 
\end{lemma}

\begin{proof}
The second claim follows from the first by duality. To show the first claim, it suffices to show that $\Ext^{>0}_{\ogroup}(\fweyl[\lambda], \fsimple[\nu]) = 0$ for all $\nu\le \lambda$. Again, noting that $\F$ is a field over $\A$, so $2\in \F^{\times}$, we can use \cite[Corollary 3.6.18]{Be-rep-cohomology} to observe that restriction induces an injection 
\begin{gather*}
\Ext^{i}_{\ogroup[\lsym]}\big(\fweyl[\lambda],\fsimple[\nu]\big)
\hookrightarrow
\Ext^{i}_{\sogroup[\lsym]}\big(\ResN(\fweyl[\lambda]),\ResN(\fsimple[\nu])\big).
\end{gather*}
Since $\Rep\big(\sogroup[\lsym]\big)$ is well-known to be a highest weight category, we can observe that if $\wta\in X^{+}\big(\sogroup[\lsym]\big)$, then $\ffweyl[\wta]$ is a projective cover of $\fsimple[\wta]$ in $\Rep\big(\sogroup[\lsym]\big)_{\leq\wta}$. Thus, $\Ext_{\sogroup[\lsym]}^{>0}\big(\fweyl[\wta],\fsimple[\wtb]\big)\cong0$ for all $\wtb\leq\wta$. The claim then follows from \cref{L:Res-standards} and \cref{R:characterizing-ON-dominance-order}.
\end{proof}

By \cite[Lemma 4.1 and Theorem 4.2]{BrSt-semi-infinite-highest-weight}, it follows that 
\begin{gather*}
\Tilt_{\F}\big(\ogroup[\lsym]\big)=\Wcat_{\F}\big(\ogroup[\lsym]\big)\cap\DWcat_{\F}\big(\ogroup[\lsym]\big)
\end{gather*}
is an additive category, with isomorphism classes of indecomposable in bijection with $\ParN$. For $\lambda\in\ParN$, we write $\ftilting[\lambda]$ for the indecomposable tilting representations with subrepresentation $\fweyl[\lambda]$. 

\begin{proposition}
The category $\Rep_{\F}\big(\ogroup[\lsym]\big)$ 
equipped with $(\ParN,\leq)$ is a(n upper finite) highest weight category.
\end{proposition}

\begin{proof}
We use \cite[Corollary 3.64]{BrSt-semi-infinite-highest-weight}
and the above discussion.
\end{proof}

\subsubsection{Combinatorial orthogonal category}

The following is an important property:

\begin{lemma}\label{L:o-tensor-weyl-is-weyl}
The tensor product of two Weyl representations (respectively dual Weyl representations) in 
$\Rep_{\F}\big(\ogroup[\lsym]\big)$ has a filtration by Weyl representations 
(respectively dual Weyl representations).
\end{lemma}

\begin{proof}
Because of compatibility of $\hcirc$ and $(\placeholder)^{\pivo}$, along with exactness of $(\placeholder)^{\pivo}$, it suffices to prove the result for 
Weyl representations. Let $\lambda,\mu\in\ParN$. There are 
$\wta,\wtb\in X^{+}\big(\sogroup[\lsym]\big)$ such that 
$\fweyl[\lambda]$, respectively $\fweyl[\mu]$, is a 
direct summand of $\IndN(\fweyl[\wta])$, respectively of
$\IndN(\fweyl[\wtb])$. Then $\fweyl[\lambda]\hcirc\fweyl[\mu]$ 
is a direct summand of
\begin{gather*}
\IndN\big(\fweyl[\wta]\big)\hcirc\IndN\big(\fweyl[\wtb]\big) 
\cong 
\IndN\Big(\fweyl[\wta]\hcirc\ResN\big(\IndN(\fweyl[\wtb])\big)\Big).
\end{gather*}
We use \cite[Proposition A2.2(vi)]{Do-q-schur} as follows.
Since a summand of Weyl filtered representations is Weyl
filtered, 
inductions and restrictions of Weyl representations are Weyl 
filtered, and tensor products of Weyl 
$\sogroup$-representations are Weyl filtered, the claim follows from exactness of $\IndN$ and $\ResN$.
\end{proof}

Note that in the proof above we work over $\F$. This is because we are using that summand of a Weyl filtered representation is Weyl filtered, which follows from highest weight category theory, and is therefore not present over $\A$.

\begin{proposition}
The category $\Tilt_{\F}(\ogroup[\lsym])$ is a symmetric 
ribbon $\F$-linear category.
\end{proposition}

\begin{proof}
\cref{L:o-tensor-weyl-is-weyl} implies that $\Wcat_{\F}\big(\ogroup[\lsym]\big)$ and $\DWcat_{\F}\big(\ogroup[\lsym]\big)$ are closed under tensor products. Being pivotal and symmetric is inherited from $\Rep_{\F}\big(\ogroup[\lsym]\big)$, and so is the ribbon property. 
\end{proof}

\begin{lemma}\label{L:SimpleTilting}
The Weyl representations $\weyl[1^{i}]$ are isomorphic to $\ext[i]=\ext[i](\A^{\lsym})$ and are simple tilting representations, for $i=0,1,\dots,\lsym$.
\end{lemma}

\begin{proof}
Recall that if $\lsym$ is even, then
\begin{gather*}
\weyl[\varpi_{i}]\cong\ResN(\ext[i])\text{ for $i\in[1,\llsym-1]$},
\quad
\weyl[2\varpi_{\lsym}]\cong\ResN(\ext[\lsym]).
\end{gather*}
Moreover, if $\lsym$ is odd, then
\begin{gather*}
\weyl[\varpi_{i}]\cong\ResN(\ext[i])
\text{ for $i\in[1,\llsym-2]$}
,\quad 
\weyl[\varpi_{\llsym-1}+\varpi_{\lsym}]\cong\ResN(\ext[\llsym-1]),
\\
\weyl[2\varpi_{\llsym-1}]\oplus\weyl[2\varpi_{\lsym}]\cong
\ResN(\ext[\llsym]).
\end{gather*}
Note that each highest weight above is fixed by $\sigma$ except for $2\varpi_{\llsym-1}$ and $2\varpi_{\lsym}$ which are permuted by $\sigma$.

It then follows from \cite[Sections 3.6.2 and 3.6.4]{JuMaWi-parity-tilting}
that each Weyl $\sogroup[\lsym]$-representa\-tion appearing above is tilting. In particular, each of these Weyl representations is isomorphic to its dual over $\A$. 

One then easily argues that independent of whether $\lsym$ is even or odd, we have $\ext[i]\cong\weyl[1^{i}]$ for $i\in[0,\llsym]$, and 
$\ext[i]\F^{\lsym}$ is in $\Tilt_{\F}(\ogroup[\lsym])$ for $i\in[0,\llsym]$.
\end{proof}

Let $\FundF\subset\Rep_{\F}\big(\ogroup[\lsym]\big)$ be the full subcategory spanned by the representations 
as in \cref{L:SimpleTilting}.

\begin{proposition}\label{P:Kar-Fund-is-Tilt}
There is an equivalence of pivotal symmetric ribbon categories 
\begin{gather*}
\Kar[\FundF]\to\TiltF.
\end{gather*}
\end{proposition}

\begin{proof}
Recall from \cref{D:dominant-orthogonal-wt} that dominant weights for $\ogroup$ are given by pairs of a dominant $\sogroup$ weight and 
$\epsilon\in\{\pm 1,0\}$.
Tensoring with the determinant $\ogroup$-representation corresponds to swapping the sign in $\epsilon$. Thus, since the determinant $\ogroup$-representation is an exterior power, the results follows 
from the same statement about $\sogroup$-representations.
\end{proof}

\begin{remark}
One could expect that there is an equivalence of symmetric ribbon $\F$-linear categories $\F\hcirc_{\A}\Fund\to\FundF$.
\end{remark}

\appendix

%%%%%%%%%%%%%%%%

\section{Some combinatorial facts for Howe's duality}\label{S:BackgroundWeights}

%%%%%%%%%%%%%%%%

We now fill in some details regarding \cref{S:Semisimple}.

\subsection{A short overview}

\begin{remark}
Nothing in this section is new. And since it can be pieced together from the literature, we will be very brief.
\end{remark}

Recall that $\ext(V\hcirc\A^{\rsym})$ is an $\ogroup$-$\uplain(\gso)^{op}$-birepresentation by \cref{P:actions-commute}, 
and we will only use the associated two actions.

There is a basis, different than the $w_{S}$ basis, which is a weight basis for both actions. 
These vectors will necessarily also be indexed by $S\subset\NBoxm$, so for each box and dot diagram, we will associate a $\ogroup$-weight and an $\gso$-weight. 

In the case of $\ogroup$ and $\gso$, dominant weights can be encoded three different ways: 
\begin{enumerate}[label=(\roman*)]

\item As a \emph{sequence} of positive integers which are immediately read off of the box and dot diagram, see \cref{D:ParN} and \cref{D:RepsLambdaWtS}.

\item As an simple representation of a \emph{maximal torus} (that is in terms of the usual notion of weight, generalized to disconnected groups), see \cref{D:ParN}.

\item As a pair of a \emph{Young diagram} (a.k.a. integer partitions) and an integer in 
$\{-1,0,1\}$, see \cref{R:three-ways-on-wts} and \cref{N:youngdiagram-for-gso-wts}.

\end{enumerate}

%%%%%%%%%%%%%%%%%%%%%%%%

\subsection{$\uplain(\gso)$-weights in $\ext(V\hcirc\A^{\rsym})$}

%%%%%%%%%%%%%%%%%%%%%%%

It turns out that there is a simpler convention for writing the
weight vectors in \cref{L:RepsMoreOActions}
instead of using the $\epsilon_{j}$ basis. Note that the 
operators $e_{j,j}^{0}=x_{j}\partial_{j}$ acting on 
$\ext(\A^{\rsym})$ can be made to act on 
$\ext(\A^{\rsym})^{\hcirc\lsym}$ by derivations, 
and then we can transport this action to 
$\ext(V\hcirc\A^{\rsym})$ with the isomorphism $\hreading$. 
If we instead think of the weight of $w_{S}^{h}$ 
in terms of the eigenvalues of the elements $e_{j,j}^0$, 
then it is possible to recover how the $h_{j}$ act from \cref{D:RepsHoweChevalleyGens}.
We will give the details now.

\begin{notation}\label{N:RepsNotationForWeightStuff}
Let $\Comp$ be the set of \emph{compositions} of length $\rsym$, i.e. tuples $\wtk=(\wtk_{1},\dots,\wtk_{\rsym})$ such 
that $\wtk_{j}\in\N$. We write $\Comp^{\lsym}$ to 
denote the subset of $\Comp$ consisting of 
$\wtk$ such that $\wtk_{j}\leq\lsym$ for $j\in[1,\rsym]$, or 
in other words, the compositions that fit into a $\lsym$-$\rsym$ rectangle. 
\end{notation}

\begin{definition}\label{D:RepsLambdaWtS}
Given $S\subset\NBoxm$, we define $\wt(w^{h}_{S})=\big(|S_{1}|,\dots,|S_{\rsym}|\big)\in\CompN$.
\end{definition}

\begin{example}\label{E:RepsLambdaWtS}
Here is an example:
\begin{gather*}
\ytableausetup{boxsize=0.62cm}
S=
\begin{ytableau}
*(magenta!50)\bullet & *(magenta!50)\bullet & \phantom{a} 
& \phantom{a} & \phantom{a} & \phantom{a}
\\
*(magenta!50)\bullet & \phantom{a} & \phantom{a} 
& \phantom{a} & \phantom{a} & *(magenta!50)\bullet
\\
*(magenta!50)\bullet & \phantom{a} & \phantom{a} 
& \phantom{a} & *(magenta!50)\bullet & \phantom{a}
\end{ytableau}
\rightsquigarrow
\wt(w^{h}_{S})=(3,1,0,0,1,1).
\end{gather*}
Indeed, in terms of dot diagrams \cref{D:RepsLambdaWtS} counts the number of dots in columns.
\end{example}

Recall that $\CompN$ denotes the set of $\gso$-weights appearing in $\ext(V\hcirc\A^{\rsym})$.

\begin{lemma}\label{L:RepsLambdaWtS}
We have $\{\wt(w^{h}_{S})|S\subset\NBoxm\}=\CompN$.
\end{lemma}

\begin{proof}
Note that $\wt(w^{h}_{S})\in\CompN$. Given $\wtk\in\CompN$, let
\begin{gather*}
S^{\wtk}=\cup_{1\leq j\leq\rsym}\{(1,j),(2,j),\dots,(\wtk_{j},j)\}.
\end{gather*}
Since $1\leq\wtk_{j}\leq\lsym$ for $j\in[1,\rsym]$, we 
have $S\subset\NBoxm$. Then from 
\cref{D:RepsLambdaWtS} we see $\wt w^{h}_{S^{\wtk}}=\wtk$.
\end{proof}

\begin{lemma}\label{L:RepsWeightSame}
Suppose $\wt(w^{h}_{S})=(\wtk_{1},\dots,\wtk_{\rsym})$. Then
we have $e_{j,j}^{o}\acts w^{h}_{S}=|S_{j}|\cdot w^{h}_{S}$ and therefore
\begin{gather*}
h_{j}\acts w^{h}_{S}= 
\big(|S_{j}|-|S_{j+1}|\big)\cdot w^{h}_{S}
,\quad
h_{\rsym}\acts w^{h}_{S}=
\big(|S_{\rsym-1}|+|S_{\rsym}|-N\big)\cdot w^{h}_{S},
\end{gather*}
where $j\in[1,\rsym-1]$.
\end{lemma}

\begin{proof}
Since $d_{j}(S) = |S_{j}| - (N-|S_{j}|)$, we have 
\begin{gather*}
\frac{1}{2}\big(d_{j}(S)-d_{j+1}(S)\big) 
=|S_{j}|-|S_{j+1}|,\text{ for $j\in[1,\rsym-1]$,}
\\
\frac{1}{2}\big(d_{\rsym-1}(S)+d_{\rsym}(S)\big) 
=|S_{\rsym-1}|+|S_{\rsym}|-\lsym.
\end{gather*}
Since $\alpha_{j}=\epsilon_{j}-\epsilon_{j+1}$ and 
$\alpha_{\rsym}=\epsilon_{\rsym-1}+\epsilon_{\rsym}$, 
the claim then follows from \cref{L:RepsMoreOActions}.
\end{proof}

It follows from \cref{L:RepsWeightSame} that
\begin{gather}\label{E:convert-wtk-to-epsilon}
(\wtk_{1},\dots,\wtk_{\rsym})
\mapsto
\sum_{i=1}^{\rsym}\big(\wtk_{i}-\tfrac{\lsym}{2}\big)\epsilon_{i}
\end{gather}
converts between weights $\wt(w^{h}_{S})\in\CompN$ 
and $\wt_{\gso}(w^{h}_{S})\in X(\gso)\subset\oplus_{i=1}^{\rsym}\Z\frac{\epsilon_{i}}{2}$ (the notation $X(\gso)$ was specified in \cref{N:RepsDynkinSO}). Note that the $\lsym$ in the $-\lsym/2$ factor depends on $S\subset\NBoxm$. Thus, we will not refer to $\Comp$ in what follows, only $\CompN$.

\begin{example}\label{E:RepsWeightSame}
For $\rsym=6$, $\lsym=1$ and $S=\{1,3,4\}$ we have
\begin{gather*}
\begin{ytableau}
*(magenta!50)\bullet & \phantom{a} & *(magenta!50)\bullet 
&*(magenta!50)\bullet & \phantom{a} & \phantom{a}
\end{ytableau}
\rightsquigarrow
\wt(w^{h}_{S})=(1,0,1,1,0,0)
\mapsto
\wt_{\gso}(x_{S})=\tfrac{1}{2}(1,-1,1,1,-1,-1),
\end{gather*}
See also \cref{E:RepsIntegralSOAction}.
\end{example}

\begin{lemma}\label{L:RepsSameWeightLemma}
Let $S\subset\NBoxm$. We have
\begin{gather}\label{E:wtk-dominant}
\begin{gathered}
\wt_{\gso}(w^{h}_{S})\in X^{+}(\gso)
\\
\Longleftrightarrow
\\
\wt(w^{h}_{S})=\wtk\in\CompN\text{ is such that }
\wtk_{1}-\tfrac{\lsym}{2}\geq\dots\geq
\wtk_{\rsym-1}-\tfrac{\lsym}{2}\geq
|\wtk_{\rsym}-\tfrac{\lsym}{2}|\geq 0,
\end{gathered}
\end{gather}
\begin{gather}\label{E:wtk-anti-dominant}
\begin{gathered}
\wt_{\gso}(w^{h}_{S})\in X^{-}(\gso)
\\
\Longleftrightarrow
\\
\wt(w^{h}_{S})=\wtk\in\CompN\text{ is such that }
\wtk_{1}-\tfrac{\lsym}{2}\leq\dots\leq
\wtk_{\rsym-1}-\tfrac{\lsym}{2}\leq-|\wtk_{\rsym}-\tfrac{\lsym}{2}|\leq 0,
\end{gathered}
\end{gather}
where $X^{-}(\gso)$ means antidominant $\gso$-weights.
\end{lemma}

\begin{proof}
Use that $\wta\in X^{+}(\gso)$ if and only if $\alpha_{i}^{\vee}(\wta)\in\Z_{\geq 0}$, for $i\in[1,\rsym]$, to deduce 
\begin{gather*}
\wta\in X^{+}(\gso)
\Leftrightarrow
\wta_{1}\geq\wta_{2}\geq\dots\geq\wta_{\rsym-1}\ge|\wta_{\rsym}|.
\end{gather*}
Then apply \cref{E:convert-wtk-to-epsilon}.
\end{proof}

\begin{notation}\label{N:RepsSameWeightLemma}
Write $\DCompN$ for the set of $\wtk$ such that \cref{E:wtk-dominant} holds and $\ACompN$ for the set of $\wtk$ such that \cref{E:wtk-anti-dominant} holds.
\end{notation}

For $\wtk\in\CompN$, we write $\ext[\wtk]=\ext[k_{1}](V)\hcirc\dots\hcirc\ext[k_{\rsym}](V)\subset \ext(V)^{\hcirc\rsym}$.
Note that since each $\ext[k_{i}]$ is a direct summand of $\ext(V)$, $\ext[\wtk]$ is a summand of $\ext(V)^{\hcirc\rsym}$. 

\begin{lemma}\label{L:wtk-space-is-Lambda-k}
Under the isomorphism $\vreading$ from \cref{L:RepsFixIso}, the summand $\ext[\wtk]\subset\ext(V)^{\hcirc\rsym}$ corresponds to the the $\wtk$-weight space of the $\uplain(\gso)$-representation $\ext(V\hcirc\A^{\rsym})$, i.e.
\begin{gather*}
\vreading^{-1}(\ext[\wtk])=\ext(V\hcirc\A^{\rsym})[\wtk].
\end{gather*}
\end{lemma}

\begin{proof}
Follows from description of $\vreading$ 
in \cref{L:RepsFixIso}, and the observation that 
$\wt x_{S}^{h}=\wtk$ if and only if 
$k_{j}=|S_{j}|$ for $j\in[1,\rsym]$.
\end{proof}

We have thus seen two ways to encode $\gso$-weights.
There is a third way that appears in the literature to encode 
the data of a $\gso$-weight in $X^{+}(\gso)$. This third way is 
by Young diagrams, ordered by the usual \emph{dominance order} $\led$
(we also write $\ltd$ etc. having the evident meaning), with our notation specified by:

\begin{example}\label{E:Dominance}
For partitions of six we have:
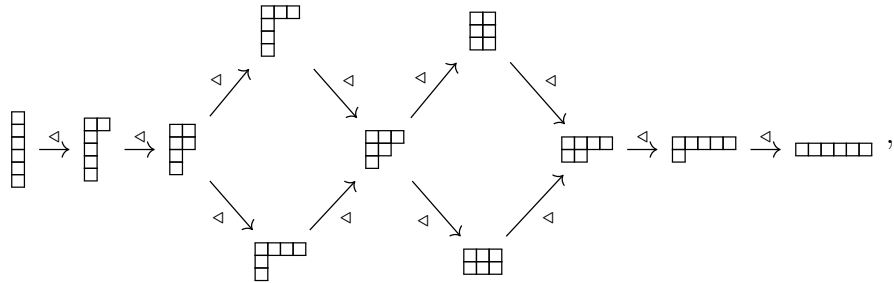
\begin{gather*}
\ytableausetup{boxsize=0.15cm}
\begin{tikzcd}[ampersand replacement=\&,column sep=1em,row sep=1em]
\&
\&
\& \ydiagram{3,1,1,1}\ar[rd,"\ltd"]
\&
\& \ydiagram{2,2,2}\ar[rd,"\ltd"]
\&
\&
\&
\\
\ydiagram{1,1,1,1,1,1}\ar[r,"\ltd"]
\&\ydiagram{2,1,1,1,1}\ar[r,"\ltd"]
\&\ydiagram{2,2,1,1}\ar[ru,"\ltd"]
\ar[rd,"\ltd"']
\&
\&\ydiagram{3,2,1}\ar[ru,"\ltd"]
\ar[rd,"\ltd"']
\&
\&\ydiagram{4,2}\ar[r,"\ltd"]
\&\ydiagram{5,1}\ar[r,"\ltd"]
\&\ydiagram{6}
\\
\&
\&
\& \ydiagram{4,1,1}\ar[ru,"\ltd"']
\&
\& \ydiagram{3,3}\ar[ru,"\ltd"']
\&
\&
\&
\end{tikzcd}
,
\end{gather*}
(this is the \emph{English convention})
with the order increases when reading left-to-right.
\end{example}

Given $\wta\in X^{+}(\gso)$, we have 
$\wta=A_{1}\frac{\epsilon_{1}}{2}
+A_{2}\frac{\epsilon_{2}}{2}+
\dots+A_{\rsym}\frac{\epsilon_{\rsym}}{2}$ such that 
$A_{i}\in\Z$, for $i\in[1,\rsym]$, 
and $A_{1}\geq\dots\geq A_{\rsym-1}
\geq|A_{\rsym}|\geq 0$. 
Moreover, either $\frac{A_{i}}{2}\in\Z$ 
for $i\in[1,\rsym]$, or $\frac{A_{i}}{2}\in\frac{1}{2}+\Z$ for $i\in[1,\rsym]$. We use this as follows:

\begin{definition}\label{N:youngdiagram-for-gso-wts}
We associate a Young diagram to $\wta$, denoted $\xp(\wta)$, with $i$th row of length
\begin{gather*}
\xp(\wta)_{i}
=
\begin{cases} 
\frac{|A_{i}|}{2} & \text{if $A_{i}$ is even} 
\\
\frac{|A_{i}|-1}{2} & \text{if $A_{i}$ is odd}.
\end{cases}
\end{gather*}
We also associate an element $\epsilon(\wta)\in\{0,\pm 1\}$ by $\epsilon(\wta)=0$, if $\wta_{\rsym}=0$, 
and $\epsilon(\wta)=\pm 1$ if 
$\wta_{\rsym}=\pm|\wta_{\rsym}|\neq 0$.
\end{definition}

Note that one can recover $\wta$ uniquely from the pair 
$\big(\xp(\wta),\epsilon(\wta)\big)$. 

\begin{example}\label{E:RepsLambdaWtSMore}
Let us take $\rsym=6$, $\lsym=2$ and the following dot diagram:
\begin{gather*}
\ytableausetup{boxsize=0.62cm}
S=
\begin{ytableau}
*(magenta!50)\bullet & *(magenta!50)\bullet & \phantom{a} 
& *(magenta!50)\bullet & *(magenta!50)\bullet & *(magenta!50)\bullet
\\
*(magenta!50)\bullet & *(magenta!50)\bullet & *(magenta!50)\bullet 
& \phantom{a} & \phantom{a} & \phantom{a}
\end{ytableau}
\rightsquigarrow
\begin{gathered}
\wt(w^{h}_{S})=(2,2,1,1,1,1),
\\
\wt_{\gso}(w^{h}_{S})=(1,1,0,0,0,0),
\\
\xp=\big((1,1),0\big).
\end{gathered}
\end{gather*}
We leave it to the reader to draw the Young diagram.
\end{example}

\begin{lemma}\label{L:so2m-dominanceorder-vs-partitionorder}
If $\wta\leq\wtb$, then $\xp(\wta)\led\xp(\wtb)$.
\end{lemma}

\begin{proof}
Standard, see the proof of \cref{L:SON-aleb-implies-YaleYb} for the $\sogroup$ version.
\end{proof}

\begin{notation}
Suppose that $\wtk,\wtl\in\DCompN$ correspond to $\wta,\wtb\in X^{+}(\gso)$. If $\wta\leq\wtb$, then we write $\wtk\legso\wtl$.
\end{notation} 

\begin{definition}
Let $\wta,\wtb\in X^{+}(\gso)$. Define a partial order 
$\big(\xp(\wta),\epsilon(\wta)\big)<\big(\xp(\wtb),\epsilon(\wtb)\big)$ 
if and only if $\xp(\wta)\ltd\xp(\wtb)$.
\end{definition}

\begin{notation}\label{N:wtk-le-wtl-partition}
Suppose that $\wtk,\wtl\in\DCompN$ correspond to $\wta,\wtb\in X^{+}(\gso)$. If 
$\big(\xp(\wta),\epsilon(\wta)\big)<\big(\xp(\wtb),\epsilon(\wtb)\big)$, then we write $\wtk<\wtl$.
\end{notation}

\begin{lemma}
Let $\wtk,\wtl\in\DCompN$. If $\wtk\legso\wtl$, then $\wtk\leq \wtl$. 
\end{lemma}

\begin{proof}
Follows from \cref{L:so2m-dominanceorder-vs-partitionorder}.
\end{proof}

%%%%%%%%%%%%%%%%%%%%%%%%%

\subsection{$\ogroup$-weights in $\ext(V\hcirc\A^{\rsym})$}

%%%%%%%%%%%%%%%%%%%%%%%%%

We follow the conventions and notation of \cref{SS:highestweights-for-on}. We use that $\sogroup\subset\ogroup$ acts on
$\ext(V\hcirc\A^{\rsym})$ by restriction.

As we observed in \cref{R:v-not-wt-basis}, the basis $v_{i}$ is not a weight basis with respect to our choice of $T\subset\sogroup$ from \cref{D:split-torus}. Thus, neither is the basis $w^{v}_{S}$ for $\ext(V\hcirc\A^{\rsym})$.

\begin{definition}
We write 
\begin{gather*}
z_{ij}= a_{i}\hcirc x_{j},
z_{N-i+1, j}=b_{i}\hcirc x_{j}
\text{ for $i\in[1,n]$}, 
\\
z_{n+1, j}=u\hcirc x_{j}\text{ if $\lsym=2\llsym+1$}.
\end{gather*}
For $S\subset\NBoxm$, we also write $z_{S}$ for the product of $z_{ij}$ such that $(i,j)\in S$ ordered by the \emph{vertical reading}, see \cref{D:vertical-horizontal-reading}, of $S$. We do not consider the horizontal reading for this basis.
\end{definition}

\begin{lemma}
The set $\{z_{S}|S\subset\NBoxm\}$ is an $\A$-basis of 
$\sogroup$-weight vectors.
\end{lemma}

\begin{proof}
Since $\sogroup\subset\ogroup$ acts on $\ext(V\hcirc\A^{\rsym})\cong \ext(V)^{\hcirc m}$ by the usual tensor product rule, this is easy to
check and omitted.
\end{proof}

\begin{lemma}
Let $S\subset\NBoxm$. Then 
\begin{gather*}
\wt_{SO}(z_{S})=\sum_{i=1}^{\lsym}
\left(|S_{i}|-|S_{N-i+1}|\right)\epsilon_{i}.
\end{gather*}
\end{lemma}

\begin{proof}
Immediate from description of the 
$T$-action in \cref{D:split-torus}.
\end{proof}

We now will describe the $\sigma$-action on the basis $z_{S}^{v}$. 

\begin{notation}
Suppose $\lsym=2\llsym$. For $S\subset\NBoxm$, write 
$\sigma(S)\subset \NBoxm$ to denote the set determined by the conditions:
\begin{gather*}
\sigma(S)_{i}=S_{i}\text{ for $i\neq n,n+1$}, 
\\
\sigma(S)_{\lsym}=\{(t, n)|(t,n+1)\in S_{n+1}\}\text{ and } \sigma(S)_{n+1}=\{(t,n+1)|(t,n)\in S_{\lsym}\}.
\end{gather*}
Here $\sigma$ is as in \cref{D:OThingy}.
\end{notation}

\begin{lemma}
We have
\begin{gather*}
\sigma\acts z_{S}=
\begin{cases}
(-1)^{|S_{n+1}|}z_{S} & \text{if $\lsym=2\llsym+1$},
\\
(-1)^{|\{x|(x,n),(x,n+1)\in S\}|}z_{\sigma(S)} 
& \text{if $\lsym=2\llsym$}.
\end{cases}
\end{gather*}
\end{lemma}

\begin{proof}
Boring and omitted.
\end{proof}

Thus, if $\lsym=2\llsym+1$, then $z_{S}$ is always a $\sigma$-eigenvector. While if $\lsym=2\llsym$, then $z_{S}$ is a $\sigma$-eigenvector if and only if $\sigma(S)=S$, i.e. the $\lsym$, $n+1$ horizontal strip of $S$ only contains empty columns or double dot columns.

\begin{remark}
Two different subsets of $\NBoxm$ can give rise to two distinct basis vectors with the same $\ogroup$-weight. For example, suppose that $\lsym=2$ and $m=2$. Then
\begin{gather*}
\wt_{O}(z_{\emptyset})=(0, 1)=\wt_{O}(z_{[1,N]\times[1,m]}),
\end{gather*}
as one easily checks.
\end{remark}

\begin{lemma}
If $(\wta,\epsilon)\in X\big(\ogroup[\lsym]\big)$, then 
$[(\wta,\epsilon)]\ext(V\hcirc\A^{\rsym})$ is an $\uplain(\gso)$-direct summand of $\ext(V\hcirc\A^{\rsym})$. 
\end{lemma}

\begin{proof}
This follows from \cref{P:actions-commute}.
\end{proof}

%%%%%%%%%%%%%%%%%%%%%%%%

\subsection{$\ogroup$-$\uplain(\gso)$-weights in $\ext(V\hcirc\A^{\rsym})$}

%%%%%%%%%%%%%%%%%%%%%%%

The formulas for the action of the raising and lowering operators in $\uplain(\gso)$ in terms of the $z_{S}$ basis are different than in the $w_{S}$ basis. However, it is still the case that $z_{S}$ is a weight vector for $\uplain(\gso)$.

\begin{lemma}\label{L:wt-of-farright-is=anitdom}
We have the following.
\begin{enumerate}[label=\arabic*.]

\item Let $S\subset\NBoxm$, then
$\wt_{\gso}(z_{S})=\wt_{\gso}(w_{S}^{h})$.

\item Let $\lambda\in\ParNm$ correspond to 
$(\wta,\epsilon)\in X^{+}\big(\ogroup[\lsym]\big)$. Consider a set 
$S\subset\NBoxm$ which has box and dot diagram with $\lambda_{i}$ dots in the $i$th row. Then $\wt_{O}(z_{S})=(\wta,\epsilon)$.

\item Let $\lambda\in\ParNm$ and consider a set $S\subset\NBoxm$ which has box and dot diagram with $\lambda_{i}$ dots in the $i$th row, and so that all dots are as far right as possible. Then $\wt z_{S}\in\ACompN$.

\end{enumerate}
\end{lemma}

\begin{proof}
A calculation.
\end{proof}

\begin{notation}
Let $\lambda\in \ParNm$. We write $S_{\lambda}$ to denote the set $S$ as in \cref{L:wt-of-farright-is=anitdom}.(c).
\end{notation}

\begin{lemma}\label{L:z-lambda-is-joint-high-low-wt-vector}
Let $\lambda\in\ParNm$. The vector $z_{S_{\lambda}}$ is a highest weight vector for $\sogroup$ and a lowest weight vector for $\uplain(\gso)$.
\end{lemma}

\begin{proof}
Another boring calculation.
\end{proof}

We prefer to not label representations by their lowest weight.

\begin{remark}\label{R:z-generates-L(lambda)-L(w0(wtb))}
Recall that the simple with lowest weight $\wtb$ has highest weight $w_{0}(\wtb)$. Since $\wtb\in X_{-}(\gso)$, it follows that $z_{S_{\lambda}}$ generates a subrepresentation of $\K\hcirc_{\A}\ext(V\hcirc\A^{\rsym})$ isomorphic to 
$\ksimple[\wta, \epsilon]\boxtimes\ksimple[w_{0}(\wtb)]$.
\end{remark}

\begin{definition}
Let $\lambda\in\ParNm$ correspond 
to $(\wta,\epsilon)$. Let 
$\wtb=\wt_{\gso}(z_{S_{\lambda}})\in X_{-}(\gso)$. 
Define $(\wta,\epsilon)^{\dagger}=w_{0}(\wtb)$.
\end{definition}

We now set out to understand the dagger operation combinatorially. 

\begin{definition}
Write $\lsym=2\llsym+1$, if $\lsym$ is odd, and $\lsym=2\llsym$, if $\lsym$ is even. Given a Young diagram $\mathcal{Y}$ with at most $\lsym$ rows and at most $\rsym$ columns, i.e. the diagram fits on an $\lsym$ by $\rsym$ checkerboard. The complement of $\mathcal{Y}$, is defined as the Young diagram with $m-\mathcal{Y}_{i}$ boxes in the $n+1-i$th row. Define $\mathcal{Y}^{\ct}$ to be the Young diagram obtained by taking the transpose of the complement of $\mathcal{Y}$.
\end{definition}

\begin{lemma}\label{L:ct-reverses-dominance-order}
If $\mathcal{Y}_{1}\led \mathcal{Y}_{2}$, then $\mathcal{Y}_{2}^{\ct}\led \mathcal{Y}_{1}^{\ct}$
\end{lemma}

\begin{proof}
Taking complements preserves the partial order on Young diagrams. The claim follows from noting that transpose reverses the partial order.
\end{proof}

\begin{lemma}
Let $\lambda\in\ParNm$ correspond to $(\wta,\epsilon)$ 
and let $\wtb=\wt_{\gso}(z_{S_{\lambda}})$. Then $-\wtb\in X^{+}(\gso)$ corresponds to $(\xp(\wta)^{\ct},\epsilon)$. 
\end{lemma}

\begin{proof}
Omitted.
\end{proof}

Note that $w_{0}(\wtb)=-w_{0}\acts (-\wtb)$. If $\rsym$ is even, then $w_{0}$ acts by $-1$, so $-w_{0}$ is the identity, while if $\rsym$ is odd, then $-w_{0}$ is the Dynkin diagram automorphism, which multiplies the $\epsilon_{\rsym}$ coordinate by $-1$. Since $-\wtb$ corresponds to $(\xp(\wta)^{\ct}, \epsilon)$, it follows that
\begin{gather*}
w_{0}(\wtb)
\leftrightsquigarrow
\begin{cases} 
(\xp(\wta)^{\ct}, \epsilon) & \text{if $\rsym$ is even}, 
\\
(\xp(\wta)^{\ct}, -\epsilon) & \text{if $\rsym$ is odd}.
\end{cases}
\end{gather*}

\begin{definition}\label{D:lambda-dagger}
Let $\lambda\in\ParNm$ correspond to $(\wta,\epsilon)$. Define $\lambda^{\dagger}\in\DCompN$ as the weight corresponding to
\begin{gather*}
(\xp(\wta),\epsilon)^{\dagger}= 
\begin{cases} 
(\xp(\wta)^{\ct}, \epsilon) & \text{if $\rsym$ is even}, 
\\
(\xp(\wta)^{\ct}, -\epsilon) & \text{if $\rsym$ is odd},
\end{cases}
\end{gather*}
depending on the parity of $\rsym$.
\end{definition}

Let $\lambda,\mu\in\ParNm$ correspond to $(\xp(\wta),\epsilon)$ and $(\xp(\wtb),\epsilon^{\prime})$. If $(\xp(\wta),\epsilon)<(\xp(\wtb),\epsilon^{\prime})$, then $\xp(\wta)\ltd\xp(\wtb)$, so by \cref{L:ct-reverses-dominance-order}, $\xp(\wtb)^{\ct}\ltd \xp(\wta)^{\ct}$, and therefore $(\xp(\wtb),\epsilon^{\prime})^{\dagger}<(\xp(\wta),\epsilon)^{\dagger}$.

\begin{proposition}\label{P:order-reversing-bijection}
The map $\dagger\colon\ParNm\rightarrow\DCompN$ is an order reversing bijection.
\end{proposition}

\begin{proof}
A calculation.
\end{proof}

\begin{lemma}\label{L:non-zero-wt-space-implies-le}
Let $\lambda\in\ParNm$ and 
$\wtk\in \DCompN$. If $[\lambda]\ext(V\hcirc\A^{\rsym})[\wtk]\neq 0$, then $\wtk\leq\lambda^{\dagger}$. 
\end{lemma}

\begin{proof}
Since $\ksimple[\mu]_{\lambda}\neq 0$ implies $\lambda\ltO\mu$ by \cref{L:Res-standards}, \cref{L:leO-implies-le-partition} gives
\begin{gather*}
[\lambda]\K\hcirc_{\A}\ext(V\hcirc\A^{\rsym})[\wtk]\subset  \bigoplus_{\lambda\leq\mu}\ksimple[\mu]\boxtimes\ksimple[\mu^{\dagger}].
\end{gather*}
It follows from $[\lambda]\K\hcirc_{\A}\ext(V\hcirc\A^{\rsym})[\wtk]\neq 0$, that there is some $\mu\in \ParNm$, such that $\lambda\leq\mu$ and $\ksimple(\mu^{\dagger})[\wtk]\neq 0$. Thus, $\wtk\leq\mu^{\dagger}\leq\lambda^{\dagger}$.
\end{proof}

\begin{lemma}\label{L:T(wtk)-is-summand}
If $\lambda\in\ParNm$, then $\ftilting[\lambda^{\dagger}]$ is a direct summand of $[\lambda]\F\hcirc_{\A}\ext(V\hcirc\A^{\rsym})$.
\end{lemma}

\begin{proof}
By \cref{L:has-Weyl-filt-for-som}, $\F\hcirc_{\A}\ext(V\hcirc\A^{\rsym})$ is a tilting representations for $\uplain[\F](\gso)$. The weight space $[\lambda]\F\hcirc_{\A}\ext(V\hcirc\A^{\rsym})$ is a $\uplain[\F](\gso)$ direct summand of $\F\hcirc_{\A}\ext(V\hcirc\A^{\rsym})$, so is a tilting representations for $\uplain[\F](\gso)$. Thus, $[\lambda]\F\hcirc_{\A}\ext(V\hcirc\A^{\rsym})$ is a direct sum of tilting representations of the form $\ftilting(\wtk)$, and by \cref{L:non-zero-wt-space-implies-le} $\wtk \le \lambda^{\dagger}$. Since $z_{S_{\lambda}}\in [\lambda]\F\hcirc_{\A}\ext(V\hcirc\A^{\rsym})[\lambda^{\dagger}]\neq 0$, and since $\ftilting(\wtk)[\wtl]\neq0$ implies $\wtl\leq\wtk$, we conclude that $\ftilting(\lambda^{\dagger})$ is a summand of $[\lambda]\F\hcirc_{\A}\ext(V\hcirc\A^{\rsym})$. 
\end{proof}

\end{document}